\documentclass[twoside, 10pt]{article}

\usepackage[hmargin= 0.75in, height=9.4in]{geometry}
\usepackage{amssymb,amsthm, amsfonts,amsmath,mathrsfs, amsxtra, url, hyperref}

\usepackage{fancyhdr}
\fancypagestyle{plain}{
 \fancyhf{}

}

\renewcommand{\theequation}{\thesection.\arabic{equation}}
 \numberwithin{equation}{section}

\newtheorem {thm}{Theorem}[section]
\newtheorem {prop}{Proposition}[section]
\newtheorem {lem}{Lemma}[section]
\newtheorem {deff}{Definition}[section]
\newtheorem {cor}{Corollary}[section]
\newtheorem {rem}{Remark}[section]

\def\ba{\begin{array}}
\def\ea{\end{array}}
\def\bea{\begin{eqnarray}}
\def\eea{\end{eqnarray}}
\def\beas{\begin{eqnarray*}}
\def\eeas{\end{eqnarray*}}
\def\bi{\begin{itemize}}
\def\ei{\end{itemize}}
\def\bc{\begin{cases}}
\def\ec{\end{cases}}

\def\a{\alpha}
\def\ga{\gamma}
\def\d{\delta}
\def\e{\varepsilon}
\def\z{\zeta}
\def\k{\kappa}
\def\l{\lambda}
\def\vr{\varrho}
\def\si{\sigma}

\def\t{\tau}
\def\th{\theta}
\def\o{\omega}
\def\vf{\varphi}
\def\vth{\vartheta}
\def\p{\psi}
\def\f{\phi}

\def\D{\Delta}
\def\G{\Gamma}
\def\L{\Lambda}
\def\O{\Omega}

\def\U{\Upsilon}

\def\bF{{\bf F}}

\def\cB{{\cal B}}
\def\cC{{\cal C}}

\def\cE{{\cal E}}
\def\cF{{\cal F}}
\def\cG{{\cal G}}

\def\cM{{\cal M}}
\def\cN{{\cal N}}

\def\cT{{\cal T}}
\def\cU{{\cal U}}

\def\cX{{\cal X}}
\def\cY{{\cal Y}}
\def\cZ{{\cal Z}}

\def\hD{\mathbb{D}}
\def\hE{\mathbb{E}}

\def\hL{\mathbb{L}}
\def\hM{\mathbb{M}}
\def\hN{\mathbb{N}}

\def\hR{\mathbb{R}}
\def\hS{\mathbb{S}}

\def\hU{\mathbb{U}}

\def\hZ{\mathbb{Z}}

\def\sB{\mathscr{B}}

\def\sD{\mathscr{D}}

\def\sJ{\mathscr{J}}

\def\sM{\mathscr{M}}

\def\sP{\mathscr{P}}

\def\sU{\mathscr{U}}

\def\sY{\mathscr{Y}}
\def\sZ{\mathscr{Z}}

\def\fC{\mathfrak{C}}

\def\fL{\mathfrak{L}}
\def\fM{\mathfrak{M}}
\def\fN{\mathfrak{N}}

\def\fS{\mathfrak{S}}

\def\fp{\mathfrak{p}}

\def\ff{\mathfrak{f}}
\def\fm{\mathfrak{m}}

\def\fy{\mathfrak{y}}
\def\fz{\mathfrak{z}}

\def\fk{\mathfrak{k}}
\def\fn{\mathfrak{n}}

\def\fe{\mathfrak{e}}
\def\fc{\mathfrak{c}}

\def\tnp{\wt{N}_\fp}

\def\ti{\n \times \n}
\def\oti{\n \otimes \n}
\def\df{\n := \n}
\def\ls{\n \le \n}
\def\gs{\n \ge \n}
\def\={\n = \n}
\def\+{\n + \n}
\def\-{\n - \n}
\def\ins{\n \in \n}
\def\ld{\n \land \n}
\def\ve{\n \vee \n}
\def\sb{\n \subset \n}
\def\>{\n > \n}
\def\<{\n < \n}

\def\({\textnormal{(}}
\def\){\textnormal{)}}
\def\[{[\n[}
\def\]{]\n]}
\def\lan{\langle}
\def\ran{\rangle}

\def\no{\noindent}

\def\ss{\smallskip}

\def\q{\quad}
\def\qq{\qquad}
\def\n{\negthinspace}
\def\dn{\n \n}
\def\tn{\n \n \n}

\def\ol{\overline}

\def\ua{\mathop{\uparrow}}
\def\da{\mathop{\downarrow}}

\def\wt{\widetilde}
\def\wh{\widehat}

\def\dtp{{\hbox{$dt  \ti  d P-$a.s.}}}
\def\dsp{{\hbox{$ds  \ti  d P-$a.s.}}}
\def\pas{{\hbox{$ P-$a.s.}}}

\def\hb{\hbox}
\def\dis{\displaystyle}
\def\cd{\cdot}
\def\cds{\cdots}

\def\fa{\,\forall \,}

\def\es{\emptyset}

\def\b1{{\bf 1}}
\def\qed{\hfill $\Box$ \medskip}

\def\esssup{\mathop{\rm esssup}}

\def\limsup{\mathop{\ol{\rm lim}}}

\newcommand{\lsup}[1]{ \underset{#1}{\limsup}}

\newcommand{\lmt}[1]{ \underset{#1}{\lim}}
\newcommand{\lmtu}[1]{ \underset{#1}{\lim} \n \ua \,}
\newcommand{\lmtd}[1]{ \underset{#1}{\lim} \n \da \,}

\begin{document}

\title{\bf $\hL^p$ Solutions of Backward Stochastic \\ Differential Equations with Jumps}

\author{
Song Yao\thanks{
\noindent Department of
  Mathematics, University of Pittsburgh, Pittsburgh, PA 15260; email: {\tt songyao@pitt.edu}. } }

\date{}

\maketitle

\begin{abstract}

Given $p \ins (1, 2)$, we study   $\hL^p$ solutions   of a multi-dimensional backward stochastic
differential equation with  jumps (BSDEJ) whose generator may not be Lipschitz continuous in $(y,z )-$variables.
  We show that such a BSDEJ with   $p-$integrable terminal data  admits a unique  $\hL^p$  solution
  by approximating the monotonic generator by  a sequence of Lipschitz generators via convolution with mollifiers
and using a stability result.

\end{abstract}



\section{Introduction}
 \label{sec:introduction}

 Let $p \ins (1, 2)$ and $T \ins (0,\infty)$. In this paper,  we study $\hL^p$ solutions of
a multi-dimensional  backward stochastic differential equation with jumps (BSDEJ)
 \bea \label{BSDEJ}
 Y_t=\xi+\int_t^T f(s, Y_s, Z_s, U_s) ds -\int_t^T Z_s \, dB_s-\int_{(t,T]} \int_\cX U_s(x)\tnp(ds,dx), \q t\in  [0,T]
 \eea
 over a probability space $(\O,\cF,P)$ on which $B$ is a Brownian motion and
 $\fp$ is  an $\cX-$valued Poisson point process independent of $B$.
    Practically speaking,   if the Brownian motion stands for the noise from the financial market,
 then the   Poisson random measure   can be interpreted as the randomness of   insurance claims.
 In the BSDEJ \eqref{BSDEJ}   with generator $f$ and terminal data $\xi$,
 a solution consists of an adapted c\`adl\`ag process $Y$, a locally square-integrable
 predictable process $Z$ and a locally $p-$integrable predictable random field  $U$.


  The backward stochastic equation (BSDE)   was introduced 
 by Bismut \cite{Bismut-73} as  the adjoint equation  for   the   Pontryagin maximum principle in
 stochastic control theory.
 Later, Pardoux and Peng  \cite{PP-90} commenced a systematical  research of    BSDEs.
 Since then, the BSDE theory  has grown rapidly    and has    been applied to
    various areas  such as   mathematical finance, theoretical economics, stochastic control and optimization,   partial differential equations,    differential geometry
  and etc, (see the references in \cite{EPQ-97, Karatzas_Soner_1998}).

   Li and Tang \cite{Tang_Li_94} introduced into the BSDE a jump term that is driven
   by a Poisson random measure  independent of the Brownian motion.
   These authors obtained the existence   of a unique solution to a BSDEJ with a Lipschitz generator and
   square-integrable terminal data.
   Then Barles, Buckdahn and Pardoux \cite{Buckdahn_Pardoux_1994,BBP_1997} showed that
   the wellposedness of  BSDEJs gives rise to a viscosity solution of a  semilinear parabolic
   partial integro-differential equation (PIDE) and thus provides a  probabilistic interpretation of  such a PIDE.
   Later, Pardoux \cite{Pardoux_BSDEJ_1995}  relaxed the  Lipschitz condition of the generator on  variable $y$
   by  assuming a monotonicity condition 
   on variable $y$ instead.
   Situ \cite{Situ_BSDEJ_1997} and Mao and Yin \cite{YM_BSDEJ_2008}  even  degenerated the monotonicity condition
   of the generator  to a weaker version so as to remove the Lipschitz condition on variable $z$.

    During the development of the  BSDE theory, some efforts were made  in relaxing  the square integrability
  on the terminal data so as to be compatible with the fact that linear BSDEs are well-posed for integrable terminal data
  or that   linear expectations have $\hL^1$ domains:
  El Karoui et al. \cite{EPQ-97} showed that for any $p-$integrable terminal data, 
  the BSDE with a Lipschitz   generator   admits a unique $\hL^p-$solution.
  Then Briand and Carmona \cite{Briand_Carmona_2000} reduced the  Lipschitz condition  of the generator
   on  variable $y$ by a strong monotonicity condition as well as a polynomial growth condition on variable $y$.
   Later, Briand et al.  \cite{BH_Lp_2003} found that the polynomial growth condition is not necessary if one uses the  monotonicity condition similar to that of   \cite{Pardoux_BSDEJ_1995}.

  In the present paper, assuming that the generator $f$ satisfies   monotonicity conditions  (H6) and (H3)  on $(y,z)$;
  that $f$ has a general growth condition (H2) on $y$, a linear growth condition (H4) on $z$;
  and that $f$ is Lipschitz continuous in $u$,  we show in Theorem \ref{thm_BSDEJ1}  that for any $p-$integrable
  terminal data $\xi$,   the BSDEJ \eqref{BSDEJ} admits a unique $\hL^p-$solution
  $(Y,Z,U)   $ (see the notations in Subsection \ref{subsect:notation}).  Consequently, we obtain
    a general martingale representation theorem for $p-$integrable martingales in the jump case (Corollary \ref{cor_martingale}).

To demonstrate  Theorem \ref{thm_BSDEJ1}, we start  with an inequality \eqref{eq:b143} about the difference of two local
$p-$integrable solutions to BSDEJs with different parameters under a general monotonicity condition \eqref{nonlip-cond2}.
   The basic inequality \eqref{eq:b143} gives rise to   an  {\it a priori} estimate \eqref{eqn-apropri} of
   the $\hL^p-$norm  of a solution $(Y,Z,U)$ of a BSDEJ with parameter $(\xi,f)$
   in terms  of the $L^p$ norms of $|\xi| \+ \int_0^T \n |f(t,0,0,0)| dt$.
       The inequality \eqref{eq:b143} also  leads to  a stability result of
     $\hL^p-$solutions of BSDEJs (Proposition \ref{prop_stab}),
 which claims that a sequence of solutions to BSDEJs is a Cauchy sequence under the $L^p-$norm
 if   their terminal data is a Cauchy sequence under the $L^p-$norm
  and if the solutions satisfy an asymptotic monotonicity condition \eqref{nonlip-cond3}.
  Then the uniqueness of the $\hL^p-$solution to BSDEJ \eqref{BSDEJ} immediately   follows.

    For the existence of an $\hL^p-$solution to BSDEJ \eqref{BSDEJ},
    we first deal with the case when the monotonic generator $f$ has linear growth (H2') in $y$  and
    when the random variable  $|\xi| \+ \int_0^T \n |f(t,0,0,0)| dt$ is bounded.
    In Proposition \ref{prop_exist_bdd}, we exploit {\it convolution with mollifiers}
    to  approach the monotonic generator $f$ by     a sequence of Lipschitz  generators,
    and utilize  the stability result (Proposition \ref{prop_stab}) to show  that
    the $\hL^2-$solutions of the BSDEJs with the approximating Lipschitz  generators and the bounded
    terminal data is actually
    a Cauchy sequence   in $\hS^p$ whose limit    solves the BSDEJ \eqref{BSDEJ}.
    Then by truncating the generator $f$ and the terminal data $\xi$ respectively,
     we employ the stability result again to obtain the general existence result in Theorem \ref{thm_BSDEJ1}.

    When     the generator $f$ is Lipschitz in $(y,z,u)$, we can use the classic {\it fixed-point} argument to
    demonstrate the existence of a unique $\hL^p-$solution   of BSDEJ \eqref{BSDEJ}
    with $p-$integrable terminal data $\xi$, see Theorem \ref{thm_BSDEJ2}.  In Proposition \ref{gen_repre},
     we  further     represent a Lipschitz generator as the   limit of the difference quotients of
 $\hL^p $ solutions to the corresponding BSDEJs.

 \ss \no {\bf Main Contributions.}

  Given $ U \ins \hU^2 $,
  unlike the case of Brownian stochastic integrals,
  the    Burkholder-Davis-Gundy inequality
  is not applicable for the $p/2-$th power of the Poisson stochastic integral
  $ \int_{(0,t]} \n \int_\cX U_s(x) \tnp(ds,dx) $, $t \ins [0,T]$: i.e.
  $E \bigg[ \underset{t \in [0,T]}{\sup} \Big( \int_{(0,t]} \n \int_\cX U_s(x) \tnp(ds,dx) \Big)^{\frac{p}{2}} \bigg] $
  can not be dominated by
    $E \Big[  \big(  \int_{(0,T]}   \n   \int_\cX  \n    |U_t(x)|^2  N_\fp(dt,dx)  \big)^{\frac{p}{2}}  \Big]$.
    Moreover, one may not be able to compare
    $E \Big[  \big(  \int_{(0,T]}   \n   \int_\cX  \n    |U_t(x)|^2  N_\fp(dt,dx)  \big)^{\frac{p}{2}}  \Big]$
    with $ E \Big[  \big(  \int_0^T   \n   \int_\cX  \n    |U_t(x)|^2  \nu (dx) dt   \big)^{\frac{p}{2}}  \Big] $.
  Then we could not follow  the classical argument in the proof of   \cite[Proposition 3.2]{BH_Lp_2003},
     neither could we employ the space  $ \hU^{2,p}  \df \Big\{ U \n  :
  E \Big[  \big(   \int_0^T  \n \int_\cX \n  | U_t (x) |^2  \nu (dx) dt  \big)^{\frac{p}{2}}  \Big]
     \<  \infty \Big\}  $ or the space
     $  \wt{\hU}^{2,p} \df \Big\{ U \n  :
  E \Big[  \big(  \int_{(0,T]}   \n   \int_\cX  \n    |U_t(x)|^2  N_\fp(dt,dx)  \big)^{\frac{p}{2}}  \Big]
     \<  \infty \Big\}  $.

  \ss   To address these technical difficulties, we first    generalize the Poisson stochastic integral
  for a random field  $U \ins \hU^p$
 by constructing in Lemma \ref{lem_stoch_integr_Lp}
   a c\`adl\`ag uniformly integrable martingale $M^U_t \df \int_{(0,t]} \int_\cX U_s(x)\tnp(ds,dx)$, $t \ins [0,T]$,
    whose   quadratic variation $[M^U  \n , M^U]$ is still
 $  \int_{(0,t]}   \n   \int_\cX  \n    |U_s(x)|^2  N_\fp(ds,dx) $, $   t  \ins  [0,T] $.
 Our inequality \eqref{eq:a017} shows that
 \bea \label{eq:a481}
 E \Big[ [M^U  \n , M^U]^{\frac{p}{2}} \Big] \ls E \int_{(0,T]} \n \int_\cX \n    |U_t(x)|^p N_\fp(dt,dx)
 \= E \int_0^T \n \int_\cX \n    |U_t(x)|^p \nu(dx) dt
 \eea
 In deriving the key $\hL^p -$type   inequality   \eqref{eq:b143} about   the difference $Y^1 \- Y^2$ of
  two  local $p-$integrable solutions  to   BSDEJs    with different parameters,
 our delicate analysis  shows that   the variational jump part of
   the dynamics   of $|Y^1\-Y^2|^p$  
 will eventually boil down  to
 the term   $E \n  \int_0^T  \n \int_\cX \n   |U^1_t (x) \- U^2_t (x)  |^p \nu (dx) dt$,
 which justifies our choice of $\hU^p$ over $\hU^{2,p}$ or $\wt{\hU}^{2,p}$ as the space for jump diffusion.
 The   estimation course of the variational jump is full of  analytical subtleties,
 but we manage to overcome them by leveraging Taylor's expansion, \eqref{eq:a481} and   some new techniques
 \big(see \eqref{eq:a451}$-$\eqref{eq:a159} for details\big).

    It is also worthy mentioning that
    although our   ``convolution with mollifiers" approach seems similar to that of \cite{Situ_BSDEJ_1997},
    some special treatments   are
    necessary along the way to overcome various technical hurdles raising in the $\hL^p-$jump case;
    and some auxiliary results,   like Lemma \ref{lem_a02}
      and  \ref{lem_a04}, are interesting in their own right.

 The  financial  significance of the present paper lies in the fact that it allows us to study
 many  mathematical finance problems for a large class of $p-$integrable financial positions
 (which may not be square-integrable)
 under nonlinear evaluation criteria or  risk measurement in a  market with jumps.
 In particular, the paper   provides a solid technical ground for our accompanying articles
 \cite{Lp_gexp,NE_Lp,RBSDEJ_Lp}:

 Given a real-valued  $p-$integrable $\xi$,
 the wellposedness result (Theorem   \ref{thm_BSDEJ1} or  \ref{thm_BSDEJ2})  shows that
 the BSDEJ with a   generator $g$ and the  terminal data $\xi$ admits a unique solution, whose
 $Y-$component $Y^\xi$ can be regarded as the   so-called ``(conditional)  $g-$expectation" of $\xi$:
   $\cE_g[\xi|\cF_t] \df Y^\xi_t$, $t \ins [0,T]$.
   In   \cite{Lp_gexp},  we first   derive a strict comparison theorem for $\hL^p-$solutions of BSDEJs
  with Lipschitz generators $g$ that satisfy an additional condition on $u-$variable (see (A3) therein).
  Then we demonstrate that the $g-$expectations, as nonlinear expectations with $\hL^p$ domains  under jump filtration,
  inherit many basic 
   properties   from the classic linear expectations
  and are  closely related to  axiom-based coherent and convex risk measures (see \cite{ADEH,FS-02,Rosazza_2006})   in mathematical finance.

   In \cite{NE_Lp},  we study  a general class of jump-filtration consistent nonlinear expectations $\cE$ with $\hL^p-$domains,
 which includes  many coherent or convex time-consistent risk measures $\rho \= \{\rho_t\}_{t \in [0,T]}$.
 Under certain domination condition, we demonstrate that the nonlinear expectation $\cE$ preserves
 many important (martingale) properties of  linear expectations
 (including optional sampling and  Doob-Meyer decomposition),
   and thus can be   represented by some $g-$expectation. 
 Consequently,  one can  utilize the BSDEJ theory to systematically analyze the risk measure $\rho$
 with $\hL^p-$domains  and  employ  numerical schemes of BSDEJs to run simulation
 for financial problems involving  $\rho$  in a  financial  market with jumps.

 Moreover,  we analyze in \cite{RBSDEJ_Lp} a   BSDEJ   with a $p-$integrable reflecting barrier $\fL$
 whose generator $g$ is Lipschitz continuous in $(y,z,u)$.
 We show that such a reflected BSDEJ with $p-$integrable parameters  admits a unique  $\hL^p$  solution,
 and thus  solves   the corresponding optimal stopping problem
  under the  $g-$expectation  or some dominated risk measure  with $\hL^p-$domain:
   the $Y-$component of the unique solution is exactly  the {\it Snell} envelope of process $\fL$ under the $g-$expectation
   and the first time it meets $\fL$ is an optimal stopping time  for maximizing  the $g-$expectation of reward $\fL$
   or minimizing the risk measure of financial position $\fL$.

\no {\bf Relevant Literature.}

 Besides the aforementioned works, we would like to make a synopsis of some recent research on BSDEJs.

\no 1)  Kruse and Popier \cite{Kruse_Popier_2016} lately studied a similar $\hL^p-$solution problem of BSDE under a
   right-continuous filtration which may be larger than the jump filtration:
\bea \label{BSDE_gf}
Y_t \= \xi \+ \int_t^T f(s, Y_s, Z_s, U_s) ds \- \int_t^T Z_s \, dB_s \- \int_{(t,T]} \int_\cX U_s(x)\tnp(ds,dx)
\- \int_t^T dM_s , \q t \ins  [0,T] ,
\eea
where $M$ is a local martingale   orthogonal to the jump filtration.
However, their wellposedness result requires  a relatively stronger monotone condition  and
   Lipschitz continuity  of $f$ in $z$ \big(see (H1) and (H3) therein\big).

      Klimsiak  
       studied $\hL^p$ solutions of reflected BSDEs under a   general right-continuous filtration in \cite{Klimsiak_2015},
      and analyzed $\hL^p$  solutions
      to  BSDEs with monotone generators and two irregular reflecting barriers in \cite{Klimsiak_2013}.


\ss \no 2)    The researches on BSDEs over   general filtered probability spaces
     have recently attracted more and more attention.
    A series of works  \cite{Buckdahn_preprint_1993,Karoui_Huang_1995,EPQ-97,BDM_2002,CFS_2008,LLQ_2011,CCR_2014}
    are dedicated to the theory of  BSDEs \eqref{BSDE_gf} but driven by a c\`adl\`ag martingale
    under a right-continuous filtration     that is also     quasi-left continuous.
    Lately,    \cite{BPTZ_2015,PPS_2016} removed  the quasi-left continuity assumption from the filtration
    so that the quadratic variation of the driving martingale
   does not need to be  absolutely continuous.
   On the other hand,   based on  a  general martingale representation result  due to
   Davis and Varaiya \cite{Davis_Varaiya_1974},
   Cohen and Elliott \cite{Cohen_Elliott_2012,CEP_2010}
   discussed the case where the driving martingales   are not {\it a priori} chosen but imposed by the filtration;
   see   Hassani and Ouknine \cite{Hassani_Ouknine_2002} for a  similar   approach
   on a BSDE in form of  a generic map  from a space of semimartingales to the spaces
   of   martingales    and those of finite-variation processes.
   Also,  Mania and Tevzadze \cite{Mania_Tevzadze_2003} and
   Jeanblanc et al. \cite{JMSS_2012} studied BSDEs for   semimartingales
   and their applications to mean-variance hedging.

 As to BSDEs driven by other discontinuous random sources,
 Xia \cite{Xia_2000} and Bandini \cite{Bandini_2015} studied  BSDEs driven by a random measure;
 Confortola et al.  \cite{Confortola_Fuhrman_2013,CFJ_2016}  considered BSDEs driven by a marked point process;
 \cite{Nualart_Schoutens_2001,BEE_2003,Ren_Fan_2009,Geiss_Steinicke_2016} analyzed  BSDEs driven by L\'evy processes;
 \cite{ABE_2010,Shen_Elliott_2011,Kharroubi_Lim_2014} discussed   BSDEs driven by a  process
 with a finite number of marked jumps.

\ss \no 3) There are also plenty of researches on quadratic BSDEJs:

To study the exponential utility maximization problem with an additional liability,  Becherer \cite{Becherer_2006} extended Kobylanski \cite{Kobylanski_2000}'s
 monotone stability approach to a jump-diffusion model  and obtained  a unique bounded solution to
 a related BSDE  driven by a random measure   whose generator may not be   Lipschitz continuous  in $u$.
   Becherer et al. \cite{BBK_2016} recently generalized this result for
   random measures of infinite activity with a non-deterministic compensator.
 Meanwhile, Morlais \cite{Morlais_2009} utilized  a similar monotone stability approach and   dynamic programming to
 show that   a special quadratic BSDEJ  with   bounded terminal data has a unique solution, whose $Y$ component
 is the value process of   an exponential utility maximization problem   with jumps.
 Morlais \cite{Morlais_2010} even obtained an existence result for such   quadratic BSDEJs  with
 exponentially integrable  terminal data.

 For   general quadratic BSDEJs with unbounded terminal data,
 Ngoupeyou \cite{Ngoupeyou_2010} and  El Karoui et al. \cite{EMN_2016}   extended
 Barrieu and El Karoui \cite{Barrieu_Karoui_2013}'s  quadratic semimartingales approach  to the jump case.
 They  managed to obtain an existence result for   quadratic-exponential BSDEJs
 (i.e.   quadratic  BSDEJs   whose generators have a exponential growth in $u$) with  unbounded terminal data.
  Also,  Jeanblanc et al. \cite{JMN_2016} described the value process of a utility optimization problem
 under Knightian-uncertainty in a jump setting as a class of quadratic-exponential BSDEJs.
 When   generators  of   quadratic-exponential BSDEJs are allowed to be locally-Lipschitz,
   Fujii and Takahashi  \cite{Fujii_Takahashi_2015}
 provided a sufficient condition  for the Malliavin's differentiability
 of such   BSDEJs  with  bounded terminal data
 while  \cite{Antonelli_Mancini_2016} could still employ   \cite{Kobylanski_2000}'s
 monotone stability approach to  show the wellposedness of such BSDEJs.

  As to different methods on quadratic BSDEJs,
   Kazi-Tani et al. \cite{KPZ_2015,KPZ_2016}
 exploited the fixed-point approach as in Tevzadze \cite{Tevzadze_2008} and
   an exquisite splitting technique  to demonstrate the wellposedness of
     quadratic-exponential BSDEJs with bounded terminal data
     and applied this result to study the related nonlinear expectations;
 Laeven and Stadje \cite{Laeven_Stadje_2014} took a duality approach to characterize
  the value of an optimal portfolio valuation problem    as the unique solution to a BSDEJ  with a convex generator
  which has at most quadratic growth   in $z$.

\ss \no 4)  For results on BSDEJs in  other interesting directions,
 see   \cite{KPZ_2015b,KPZ_2015c} for second-order BSDEs with jumps and  the related fully-nonlinear PIDEs;
 see Cohen and Elliott \cite{Cohen_Elliott_2008,Cohen_Elliott_2010,Cohen_Szpruch_2012}
for BSDEs driven by Markov chains;
see Kharroubi et al. \cite{KMPZ_2010} for (minimal) solutions to  BSDEs with constrained jumps and related
quasi-variational inequalities; see Aazizi and Ouknine \cite{Aazizi_Ouknine_2011} for a class of constrained BSDEJs
 and its application in pricing and hedging  American options;
 see Klimsiak and   Rozkosz \cite{Klimsiak_Rozkosz_2013,Klimsiak_Rozkosz_2016} for
   a general (non-Markovian)  BSDE and  a related semilinear elliptic equation with measure data
whose operator is associated with a regular semi-Dirichlet form;
  see  \cite{Kruse_Popier_2016b,GHQ_2015} for BSDEJs with singular terminal data
  and their applications to optimal position targeting and a non-Markovian liquidation problem respectively;
  see also \cite{Geiss_Labart_2016} for  numerical simulation of BSDEJs by Wiener chaos expansion  among other.

   The rest of the paper is organized as follows: In Section \ref{sec:introduction},
   we list  necessary notations, and we generalize the Poisson stochastic integral for   $U \ins \hU^p$
   so as to define   BSDEJs in $\hL^p$ sense.   After making some assumptions
   on generator $f$ \big(including the monotonicity conditions in $(y,z)$\big),
   we present in Section \ref{sec:main}, the   main result of our paper,
   the existence and uniqueness of an $\hL^p-$solution to a BSDEJ
   with $p-$terminal data, which gives rise to
   a general martingale representation theorem for $p-$integrable martingales in the jump case.
   In Section \ref{sec:prepare}, we give an  inequality  about   the difference of
  two  local $p-$integrable solutions  to   BSDEJs  as well as two consequences of it:
  an {\it a priori} estimate and a stability result of $\hL^p-$solutions of  BSDEJs,
 both are important to prove Theorem \ref{thm_BSDEJ1}.  Section \ref{sec:prepare} also
 includes a basic existence result of $\hL^p-$solutions to  BSDEJs with bounded parameters,
 which is also crucial for Theorem \ref{thm_BSDEJ1}.
 Section \ref{sec:Lip}  further discusses the wellposedness
 of  BSDEJs with Lipschitz generators  in $\hL^p$ sense and the corresponding generator representation.
 The proofs of our results are deferred to Section \ref{sec:proof}, and the Appendix contains
 some necessary technical lemmata.

\subsection{Notation and Preliminaries} 

\label{subsect:notation}

 Throughout this paper, we fix a time horizon $T \ins (0,\infty)$ and consider a  complete probability
 space $(\O,\cF, P )$ on which a $d-$dimensional  Brownian motion $B$ is defined.

  For a generic c\`adl\`ag process $X$, we denote  its corresponding jump process by
  $\D X_t \df  X_t   \-   X_{t-} $, $t  \ins  [0,T]$  with $X_{0-} \df X_0$.
 Given a measurable space $(\cX, \cF_\cX)$,
 let $\fp$ be an $\cX   -$valued Poisson point process
 on $(\O,\cF, P)$ that is independent of $B$.
 For any scenario $\o \in \O$,
 let $D_{\fp(\o)}$ collect all jump times of the path $\fp(\o)$, which is a countable subset of $(0,T]$
 (see e.g. Section 1.9 of \cite{SDE-diffu}). We assume that for some  finite measure $\nu$  on $\big(\cX,\cF_\cX\big)$,
the counting measure $N_\fp(dt, dx)$ of $\fp$
on $[0,T] \ti  \cX$  has  compensator $E\big[ N_\fp(dt, dx)\big] \=  \nu(dx)   dt$.
   The corresponding compensated Poisson random measure $\tnp$ is
   $\tnp(dt, dx)  \df  N_\fp(dt, dx)  \-  \nu(dx)   dt$.

   For any $t  \ins  [0, T]$, we define sigma-fields
 \beas
    \cF^B_t := \si\big\{B_s ;   s \le t\big\},  \q
   \cF^N_t := \si\big\{N_\fp \big((0,s],A \big) ;  s \le t, A \in \cF_\cX \big\},
    \q \cF_t  \df  \si\big( \cF^B_t  \n \cup \n  \cF^N_t \big)
 \eeas
 and augment   them by all  $P-$null sets in $\cF$.
Clearly, the jump filtration $\bF  \= \{\cF_t \}_{t  \in [0, T]}$
is complete and right-continuous (i.e. satisfies the {\it usual hypotheses}, see e.g.,  \cite{Pr-90}).
   Let $\sP$ \big(resp. $\wh{\sP}$\,\big)  denote the  $\bF-$progressively measurable
 (resp. $\bF-$predictable) sigma-field on $[0,T]  \ti  \O  $,
 and let    $\cT$   collect  all   $\bF-$stopping times.

 \if{0}
   For any $i \in \hN$,  as $ \big\{  N_\fp ((0,t],\cX) \big\}_{t \in [0,T]}$ is
  an $\bF -$adapted c\`adl\`ag process,
 its $i-$th jump time
  \beas
  \tau^N_i  : =  \inf\big\{t  \in  (0,T]  :  N_\fp ((0,t],\cX)  \ge  i \big\}
  \eeas
 is an $\bF-$stopping time.
 For each $\o  \ins  \O$, $\tau^N_i (\o)$ is   the $i-$th smallest element in $D_{\fp(\o)}$
 or $T$. Put in another way, $\tau^N_i (\o)  \df  \min\big\{ t  \ins  D_{\fp(\o)}  \n :
 t  \>  \tau^N_{i-1} (\o)   \big\} \ld T$, with $\tau^N_0 (\o)  \df  0$
 and $\min \es  \df  \infty$.

 \fi

  For  a generic   Euclidean space $\hE$ with norm $\|\cd \|$,   we define:
 \beas
 \sD(x) := \b1_{\{x \ne 0\}} \frac{1}{\|x\|}  x\q \hb{and} \q
     \pi_r(x)   := \frac{r   }{r \vee \|x\|} x,      \qq      \fa  x \in \hE , \q \fa r \in (0,T].
  \eeas
   See  Lemma \ref{lem_pi} and Lemma \ref{lem_a04}    for the  properties of these two functions.

   Given $l \in \hN$, the following spaces  of functions will be used in the sequel:

 \no 1)   For any $p \ins [1,\infty)$, let $L^p_+[0,T]$ be the space of all measurable functions
$\p \n : [0, T]  \n \to \n  [0, \infty)$ with $ 
\int_0^T \n \big(\p(t)\big)^p dt  \<  \infty$.

     \no 2)  For $p \ins (1,2]$, 
     let  $L^p_\nu  \df  L^p (\cX, \cF_\cX, \nu; \hR^l) $ be
       the space of all $\hR^l-$valued, $\cF_\cX-$measurable functions $u$
    with     $  \|u\|_{L^p_\nu} \df \big( \int_\cX   \n   |u(x)|^p \nu(dx) \big)^{\frac{1}{p}}   \<  \infty$.
    For any $u_1,u_2 \ins L^p_\nu$, we say $u_1 \= u_2$ if $u_1 (x) \= u_2  (x)  $ for
    $\nu-$a.s. $x \in  \n \cX \n $.

  \no 3) For any sub-sigma-field $\cG$ of $\cF$, let

  \no  $\bullet$  $L^0_+(\cG )$ be the space of all real-valued non-negative $\cG-$measurable random variables;

   \no  $\bullet$   $L^p_+(\cG )  := \Big\{ \xi \in L^0_+(\cG ) :  \|\xi\|_{L^p_+(\cG )} := \big\{E\big[\xi^p\big]\big\}^{\frac{1}{p}}<\infty  \Big\}$ for all $p \in [1,2)$;

    \no   $\bullet$   $L^\infty_+(\cG ) := \Big\{ \xi \in L^0_+(\cG ) : \|\xi\|_{L^\infty_+(\cG )} := \underset{\o \in   \O}{\esssup} \, \xi(\o)  <\infty \Big\}$;

   \no  $\bullet$  $L^0(\cG )$ be the space of all $\hR^l-$valued,
$\cG-$measurable random variables;

   \no  $\bullet$   $L^p(\cG )  := \Big\{ \xi \in L^0(\cG ) :  \|\xi\|_{L^p(\cG )} := \big\{E\big[|\xi|^p\big]\big\}^{\frac{1}{p}}<\infty  \Big\}$ for all $p \in [1,2)$;

   \no   $\bullet$    $L^\infty(\cG ) := \Big\{ \xi \in L^0(\cG ) : \|\xi\|_{L^\infty(\cG )} := \underset{\o \in   \O}{\esssup} \, |\xi(\o)|  <\infty \Big\}$.

 \no 4) Let    $\hD^0$ be the space of all $\hR^l-$valued, $\bF-$adapted c\`adl\`ag processes,
 and let $\hD^\infty$ be the space of all $\hR^l-$valued, $\bF-$adapted c\`adl\`ag  processes $X$ with
 $ \|X\|_{\hD^\infty} \df  \underset{(t,\o) \in [0,T] \times \O}{\esssup} |X_t(\o)|
  \= \underset{ \o  \in   \O}{\esssup} \, X_* (\o)  \< \infty $,
  where $ X_* (\o) \df  \underset{t \in [0,T]  }{\sup} \big| X_t  (\o)  \big|  $.

 \no 5)  Set    $\hZ^2_{\tiny \rm loc}  \df  L^2_{\hb{\tiny \rm loc}}
 \big( [0,T] \ti \O , \wh{\sP} , dt \ti dP; \hR^{l \times d} \big) $,
      the space of all $\hR^{l \times d}-$valued, $\bF -$predictable processes $Z  $ with
 $ \int_0^T \n  | Z_t |^2   \, dt   \\ <  \n  \infty $, \pas ~

 \no 6)  For any $p \in [1,2]$, we   let

 \no $\bullet$
     $\hD^p \df \Big\{ X  \ins  \hD^0 \n :  \| X \|_{\hD^p } \df
  \big\{ E  [   X^p_*   ] \big\}^{ \frac{1}{p}}  \<  \infty \Big\} $.

 \no $\bullet$
   $\hZ^{2,p}   \df
   \Big\{Z  \ins  \hZ^2_{\tiny \rm loc}  \n : \| Z \|_{\hZ^{2,p}  }
     \df  \Big\{ E \Big[  \big(   \int_0^T \n  | Z_t |^2   \, dt  \big)^{\frac{p}{2}}
    \Big] \Big\}^{\frac{1}{p}}   \<  \infty \Big\}$.
    We will   simply denote     $\hZ^{2,2}$   by $\hZ^2$.
           For any $Z \ins \hZ^{2,p}$,  the Burkholder-Davis-Gundy inequality    implies  that
   \bea \label{eq:f347}
       E   \bigg[ \, \underset{t \in [0, T]  }{\sup} \Big|\int_0^t \n Z_s dB_s  \Big|^p
        \bigg]   \ls  c_{p,l}  E   \bigg[    \Big( \int_0^T  \n
       \big|  Z_s   \big|^2    ds \Big)^{\n \frac{p}{2}}
      \bigg]            \<  \infty
      \eea
   for some constant $c_{p,l} \> 0$ depending on $p$ and $l$.
 So   $ \big\{ \int_0^t \n Z_s dB_s \big\}_{t \in [0,T]} $ is a  uniformly integrable    martingale.

  \no $\bullet$
 $\hU^p_{\tiny \rm loc}  \df  L^p_{\hb{\tiny \rm loc}} \big( [0,T]   \ti    \O    \ti    \cX  ,
 \wh{\sP}    \oti    \cF_\cX, dt  \ti  dP   \ti  \nu(dx); \hR^l \big) $
 be   the space of  all   $\wh{\sP}    \oti    \cF_\cX -$measurable
 random fields $ U \n : [0,T]   \ti    \O    \ti    \cX  \n \to \n  \hR^l $ such that
 $ \int_0^T   \n   \int_\cX    \n   |U_t (x) |^p \nu (dx) dt
  \=   \int_0^T  \n  \|U_t \|^p_{L^p_\nu}   dt \<  \infty $,    \pas ~
 For any $U \ins \hU^p_{\tiny \rm loc}$, it is clear that
 $U (t , \o) \ins L^p_\nu$ for \dtp ~ $(t,\o) \ins [0,T] \ti \O$.

 \no $\bullet$
        $\hU^p \df   \Big\{U  \ins  \hU^p_{\tiny \rm loc}  \n :  \| U \|_{\hU^p  }
     \df  \big\{ E \n  \int_0^T  \n \int_\cX \n   |U_t (x) |^p \nu (dx) dt  \big\}^{\frac{1}{p}}   \<  \infty \Big\}
     \=   L^p  \big( [0,T]   \ti    \O    \ti    \cX  ,
 \wh{\sP}    \oti    \cF_\cX, dt  \ti  dP   \ti  \nu(dx); \hR^l \big) $.

 \no $\bullet$ Let us simply denote   $ \hD^p \ti \hZ^{2,p} \ti \hU^p $ by $\hS^p$.

        In this paper,        we   use the convention   $ \inf \es := \infty$
      and  let $c_{p,l}$ denote  a generic constant depending only on $p$
       and $l$ (in particular,  $c_l$ stands for a generic constant depending only on $l$),
       whose form may vary from line to line.

 \subsection{Generalization of Poisson Stochastic Integrals} 
\label{subsection:PSI_Lp}

 The stochastic integral  with respect to the compensated Poisson random measure $\tnp (dt, dx)$
 (or simply ``Poisson stochastic integral")
 is usually defined for locally square integrable random fields $U \ins \hU^2_{\tiny \rm loc}$.
In this subsection, we will generalize such kind of stochastic integral   for
random fields in $ \underset{p \in [1,2)}{\cup} \hU^p_{\tiny \rm loc} $   in spirit of \cite[VIII.75]{Proba_pot_2}.

 Let $ \hM^1  $ be the space of all  c\`adl\`ag  local martingales
 $M \= \{M_t\}_{t \in [0,T]}$ 
 with $\|M\|_{\hM^1} \df    E \big\{\big[M,M\big]^{\frac12}_T   \big\}   \< \infty$.
 According to \cite[VII.81-VII.92]{Proba_pot_2}, $\|\cd \|_{\hM^1}$ is a norm on $\hM^1  $ that  is equivalent to
 $\|\cd\|_{\hD^1} $, 
 thus $\big( \hM^1 , \|\cd \|_{\hM^1}  \big)$ is a Banach space.

  Let $p \ins [1,2)$ and   $U  \ins  \hU^p$. For any $n  \ins  \hN$,
since $E \n \int_0^T \dn \int_\cX \n  \b1_{\{|U_s(x)|\le n\}}
|U_s(x)|^2 \nu(dx) ds \ls n^{2-p} E \n \int_0^T \dn \int_\cX   |U_s(x)|^p \nu(dx) ds \\ < \n \infty $,
 $ M^{U,n}_t  \df   \int_{(0,t]}  \n \int_\cX \n  \b1_{\{|U_s(x)|\le n\}}
U_s(x) \tnp(ds,dx)    $, $t  \ins  [0,T]  $    defines a square integrable martingale.

\begin{lem} \label{lem_stoch_integr_Lp}
 Let $p \ins [1,2)$. For any    $U  \ins  \hU^p$,
 $ \{M^{U,n}\}_{n \in \hN} $ is a Cauchy sequence in $\big( \hM^1 , \|\cd \|_{\hM^1}  \big)$, whose limit
   $M^U$  is a c\`adl\`ag uniformly integrable martingale with   quadratic variation
 $[M^U  \n , M^U]_t \= \int_{(0,t]}   \n   \int_\cX  \n    |U_s(x)|^2  N_\fp(ds,dx) $, $   t  \ins  [0,T] $.
The jump process of $M^U$ satisfies that for \pas ~ $\o \ins \O$,
 \bea \label{eq:a021}
 \D M^U_t (\o)  \=  \b1_{\{t \in D_{\fp(\o)} \}}   U \big(t,\o,\fp_t(\o)\big) , \q \fa t \ins (0,T] .
\eea
 Moreover,  $U  \n \to \n  M^U$ is a linear   mapping on $\hU^p$.

\end{lem}

 We shall assign $M^U$ as the Poisson stochastic   integral
  \bea  \label{def_Poisson_integr}
    \int_{(0,t]} \n \int_\cX U_s(x) \tnp(ds,dx) , \q  t \in [0,T]
  \eea
 of $U \ins \hU^p $. Analogous to the classic extension
 of Poisson stochastic integrals from $\hU^2$ to  $\hU^2_{\tiny \rm loc}$, one can define
 the   stochastic   integral \eqref{def_Poisson_integr}
(or simply $M^U$) for any $U  \ins  \hU^p_{\tiny \rm loc} $,
 which is a c\`adl\`ag local martingale with quadratic variation
 $ \int_0^t  \n  \int_\cX  \n  |U_s(x)|^2  N_\fp(ds,dx) $, $   t  \ins  [0,T] $
and whose jump process satisfies \eqref{eq:a021} also.
 This generalized Poisson stochastic integral is still linear in $U \ins \hU^p_{\tiny \rm loc}$.

\subsection{BSDEs with Jumps}

 From now on, let us fix    $p \ins (1, 2)$.
 A mapping $f  \n :    [0, T]  \ti  \O  \ti  \hR^l  \ti  \hR^{l \times d}
  \ti   L^p_\nu   \n \to \n  \hR^l $ is called a {\it $p-$generator} if
 it is  $\sP  \oti  \sB  ( \hR^l   )   \oti  \sB \big( \hR^{l\times d} \big)
  \oti  \sB\big(L^p_\nu   \big)/\sB(\hR^l ) -$measurable.   For any $\tau \ins \cT$,
  \beas
   f_\tau (t,\o,y,z,u) \df \b1_{\{t < \tau (\o)\}} \, f (t,\o,y,z,u) , \q
 \fa (t,\o,y,z,u) \ins [0, T] \ti  \O  \ti  \hR^l  \ti  \hR^{l \times d}  \ti  L^p_\nu
  \eeas
  is also $\sP  \oti  \sB  ( \hR^l    )   \oti  \sB  ( \hR^{l \times d}  )
  \oti  \sB\big(L^p_\nu   \big)/\sB(\hR^l  ) -$measurable.

 \begin{deff} \label{def_BSDEJ}
 Given $p \ins (1,2) $,  let $\xi \ins L^0(\cF_T)$ and $f$ be a $p-$generator.
 A triplet of processes   $(Y,Z,U) \ins \hD^0  \ti  \hZ^2_{\tiny \rm loc}   \ti  \hU^p_{\tiny \rm loc}$
 is called a solution of a backward stochastic differential equation with jumps that has terminal data $\xi$
 and  generator $f$ \big(BSDEJ\,$(\xi,f)$ for short\big) if
 $   \int_0^T  \n    |f (s, Y_s, Z_s, U_s ) |   ds   \<  \infty $,   \pas ~
 and if \eqref{BSDEJ} holds \pas ~

 \end{deff}

 \begin{rem} \label{rem_integral} Let $p \ins (1,2) $.

\no \(1\) Let $U  \ins    \hU^p_{\tiny \rm loc}$. For any   $\tau \ins \cT$, since
   $\big\{\b1_{\{t \le \tau  \}}\big\}_{t \in [0,T]}$
   is   an $\bF-$adapted c\`agl\`ad   process \(and thus $\bF-$predictable\),
   the process  $\big\{\b1_{\{t \le \tau  \}}U_t\big\}_{t \in [0,T]}$   also belongs to $ \hU^p_{\tiny \rm loc}$.
 By Subsection \ref{subsection:PSI_Lp}, 
 the stochastic integral  $\int_{(0, \tau]} \n \int_\cX \n U_s (x)   \tnp (ds, dx)
   \=  \int_{(0, T]} \n \int_\cX \n \b1_{\{s \le \tau  \}}  U_s (x) \tnp (ds, dx)$
 is well defined.      More general, the stochastic integral
   $\int_{(\tau, \ga]} \n \int_\cX \n U_s (x) \tnp (ds, dx)$ is valid for any
   $\tau ,\ga \ins \cT$ with   $\tau  \ls  \ga $, \pas

  \no \(2\) Given $\xi \ins L^0(\cF_T)$ and    a $p-$generator $f$,
 let $(Y,Z,U)$ be a solution of   BSDEJ\,$(\xi,f)$ as described in Definition \ref{def_BSDEJ}.
 For \pas ~ $\o \ins \O$,  we see from  \eqref{BSDEJ} and \eqref{eq:a021}   that
   \bea    \label{eq:b617}
   \D Y_t(\o) =   \D M^U_t (\o)
   =   \b1_{\{t \in D_{\fp(\o)}\}} U \big(t,\o, \fp_t (\o)\big), \q \fa t \in  [0,T] ,
   \eea
   which implies that
  \bea \label{Y_jump}
  \big\{t \in  [0,T]:   Y_{t-} (\o)  \ne Y_t (\o)      \big\} \subset D_{\fp(\o)}
  \hb{ is a countable subset of $[0,T]$. }
  \eea

 \end{rem}

 \section{Main Result}
 \label{sec:main}

 In the rest of this paper,  we set $q \df \frac{p}{p-1} \> 2$
 and let   $  \beta $ be a $[0,\infty)-$valued,
 $\bF-$progressively measurable process with    $  \int_0^T \beta_t dt  \ins L^\infty_+(\cF_T)   $.
  We   make the following assumptions on $p-$generators $f$:

   \no {\bf (H1)} For each $(t,\o, u)  \ins  [0,T]  \ti  \O \times L^p_\nu$,
   the mapping $ (y,z) \n \to \n  f(t,\o,y,z,u)$ is continuous.

   \no {\bf (H2)} For any $\d \>0$, there exists a $[0,\infty)-$valued, $\bF -$progressively measurable process
   $ \phi^\d$ with   $ E \int_0^T \phi^\d_t dt  \< \infty$ such that
   $   \underset{|y| \le \d}{\sup} \big|  f (t,  y ,0,0   ) \- f (t,  0 ,0,0   )  \big|     \ls \phi^\d_t  $, \dtp

   \no {\bf (H3)} It holds for   $\dtp$ $(t,\o)  \ins  [0,T]  \ti  \O$ that
 \beas
 && \hspace{-2cm}   \big\lan   y    ,  f (t,\o,y ,0,0   ) \- f (t,\o,0 ,0,0   )  \big\ran
 \ls    \beta (t,\o) |y |^2    , \q \fa  y \ins  \hR^l   .
 \eeas

   \no {\bf (H4)} For some $c_1 (\cd)  \ins     L^2_+ [0,T]$, it holds for \dtp~$(t,\o) \in [0,T]  \ti  \O$ that
   \beas
     |f(t, \o, y,z,0) \- f(t, \o, y,0,0)|  \le  \beta (t,\o) \+  c_1(t)|z| ,  \q \fa (y,z) \in  \hR^l \times \hR^{l \times d} .
     \eeas

     \no {\bf (H5)} For some $c_2(\cd)  \ins    L^q_+ [0,T]$, it holds for   \dtp ~ $(t,\o) \in [0,T]  \ti  \O$ that
 \beas 
 \big| f(t,\o,y,z,u_1) \- f(t,\o,y,z,u_2) \big|
 \ls c_2 (t)   \|  u_1  \- u_2 \|_{L^p_\nu}  \,  ,
 \q \fa (y,z,u_1,u_2)  \ins  \hR^l  \ti  \hR^{l \times d}  \ti   L^p_\nu  \ti   L^p_\nu \,.
 \eeas

  \no {\bf (H6)}
  It holds for   $\dtp$ $(t,\o)  \ins  [0,T]  \ti  \O$ that
 \beas
 && \hspace{-2cm} |y_1 \- y_2|^{p-1} \big\lan \sD ( y_1 \- y_2 ) ,  f (t,\o,y_1,z_1,u   )
 \- f (t,\o,y_2,z_2,u   )   \big\ran
 \ls \l(t) \, \th \big( |y_1 \- y_2|^p \big) \+  \Phi (t,\o) |y_1 \- y_2|^p \\
&&     +    \L (t,\o)  |y_1 \- y_2|^{p-1}   |z_1 \- z_2|   , \q
 \fa  (y_1,z_1 ) , (y_2,z_2 ) \ins  \hR^l  \ti  \hR^{l \times d} , ~ \fa u   \ins   L^p_\nu   ,
 \eeas
 where $\l(\cd) \ins L^1_+ [0,T]$;
     $\th \n : [0, \infty)  \n \to \n  [0, \infty) $ is an increasing  concave function
   satisfying $ \int_{0+}^1  \n   \frac{1}{ \th(t) } dt  \=  \infty$; and
 $\Phi,\L$ are two $[0,\infty)-$valued, $\sB [0,T] \oti  \cF_T-$measurable process such that
   $  \int_0^T \n (\Phi_t \ve \L_t^2) dt  \ins   L^\infty_+(\cF_T)$  and
   $   E \n \int_0^T \n   \L^{2+\fe}_t  dt  \<  \infty$ for some $\fe  \ins   (0, 1  )$.

\begin{rem}  \label{rem_drift}
Given $p \ins (1,2)$,  let $f$ be a $p-$generator  satisfying \(H2\), \(H4\), \(H5\) and that
   $\int_0^T \n |f (t, 0, 0, 0 )| dt \< \infty$, \pas ~  Then
 it holds for any $ (Y,Z,U) \ins \hD^1 \ti \hZ^2_{\tiny \rm loc} \ti \hU^p_{\tiny \rm loc}$
 that  $ \int_0^T  \n    |f (t, Y_t, Z_t, U_t ) |   dt   \<  \infty $,   \pas

\end{rem}

  For simplicity,   set $ \ol{C} \df   \Big( \int_0^T \n  \big(c_1(t)\big)^2  dt \Big)
    \vee    \Big( \int_0^T \n  \big(c_2(t)\big)^q  dt \Big) $,
    $C_\beta \df  \big\| \int_0^T \beta_t dt \big\|_{L^\infty_+(\cF_T) }   $,
    $C_\Phi \df \big\| \int_0^T \Phi_t dt \big\|_{L^\infty_+(\cF_T) }  $
    and $C_{\L} \df  \big\| \int_0^T \L^2_t dt \big\|_{L^\infty_+(\cF_T) }   $.

  Our main goal  is the following existence and uniqueness result of BSDEJs for  case ``$p \in (1,2)$".

 \begin{thm} \label{thm_BSDEJ1}
 Given $p \ins  (1,2)$,  let
  $\xi  \ins  L^p(\cF_T)$ and let  $f$ be a $p-$generator satisfying  \(H1\)$-$\(H6\)
 such that  $\int_0^T \n |f(t,0,0,0)| dt \ins L^p_+ (\cF_T)$ and that
 the parameter $c_2(\cd) \ins   L^{q'}_+ [0,T]$ for some $q' \ins (q, \infty)$.
  Then the BSDEJ\,$(\xi,f)$ admits a unique solution
  $(Y,Z,U)  \ins  \hS^p $.

 \end{thm}

 \if{0}
      Now, let $\tau$ be an  $\bF-$stopping time   that may take the infinite value $\infty$.
   Thanks to Theorem \ref{thm_BSDEJ1}, the BSDEJ with random time horizon $\tau$ is also well-posed
   for any terminal data $\xi \in L^p(\cF_\tau)$ under hypotheses (H1)-(H6).

\begin{cor}   \label{cor-case1}
 Given $p \ins  (1,2)$ and $\tau \ins \cT$,
   let    $\xi  \ins  L^p(\cF_\tau)$ and $f$ be a $p-$generator.
If the mapping
\bea  \label{def_f_tau}
f_\tau (t,\o,y,z,u) : =  \b1_{\{t \le \tau(\o)\}}f(t,\o,y,z,u) , \q \fa (t,\o,y,z,u)
\ins [0, T]  \ti  \O  \ti  \hR^l  \ti  \hR^{l \times d}  \ti   L^p_\nu
\eea
satisfies   \(H1\)$-$\(H6\),  the following BSDEJ
    \bea   \label{BSDEJ_random}
   Y_{\tau \land t}=\xi+\int_{\tau \land t}^\tau f(s, Y_s, Z_s, U_s) ds -\int_{\tau \land t}^\tau Z_s \, dB_s-\int_{( \tau \land t,\tau]} \int_\cX U_s(x)\tnp(ds,dx), ~\,  \fa t\in  [0,T]; \q  \pas
    \eea
admits a unique solution $\Big\{\big(Y_t(\o), Z_t(\o), U_t(\o)\big)\Big\}_{(t,\o) \in \[0, \tau\]}   $ such that  $ \Big\{\Big(Y_{\tau \land t}, \b1_{\{t \le \tau\}} Z_t, \b1_{\{t \le \tau\}} U_t\Big)\Big\}_{t \in [0,T]} \ins \hS^p  $.

\end{cor}

 \fi

 This wellposedness gives rise to  a general martingale representation theorem  in the jump case as follows:

 \begin{cor}   \label{cor_martingale}
 Let $p \ins  (1,2)$.
 For any $\xi \in L^p (\cF_T)$, there exists a unique pair
 $(Z,U) \in   \hZ^{2,p}      \ti  \hU^p $  such that \pas
 \bea   \label{eq:a235}
  E[\xi|\cF_t]   =     E[\xi]    +    \int_0^t    Z_s dB_s
    +    \int_{(0,t]}   \int_\cX   U_s (x) \tnp (ds,dx) , \q t  \ins  [0,T] .
 \eea
 \end{cor}

\section{A priori Estimate and Stability Result}

\label{sec:prepare}

  To prove Theorem \ref{thm_BSDEJ1}, we started with       an inequality   about   the difference of
  two  local $p-$integrable solutions  to   BSDEJs    with different parameters   under a general monotonicity condition.

\begin{lem}   \label{lem-a-priori}

Let $  p  \ins (1, 2 ) $.
    For $i \= 1,2$,    let $\xi_i \ins L^0(\cF_T)$,
        let $  f_i   $   be  a $p-$generator,   and let
    $(Y^i, Z^i, U^i) \ins \hD^0  \ti  \hZ^2_{\tiny \rm loc}   \ti
     \hU^p_{\tiny \rm loc}  $ be a solution of BSDEJ\,$(\xi_i,  f_i)$
  such that $ Y^1 \- Y^2 \ins \hD^p$.   Assume that \dsp
  \bea
     && \hspace{-1.5cm}    | Y^1_s \-  Y^2_s|^{p-1}   \big\lan \sD (  Y^1_s \-  Y^2_s ) ,
   f_1  (s,Y^1_s, Z^1_s, U^1_s  )  \-  f_2  (s,Y^2_s,  Z^2_s, U^2_s  ) \big\ran
     \nonumber    \\
  &&   \ls    |Y^1_s  \-  Y^2_s|^{p-1} \Big[ g_s +  \Phi_s  |Y^1_s  \-  Y^2_s|
   \+  \L_s      |Z^1_s  \-  Z^2_s  |
  \+ \G_s \big\| U^1_s   \-  U^2_s   \big\|_{L^p_\nu}   \Big]
     \+  \U_s    ,  \q   \label{nonlip-cond2}
 \eea
 where
   $g$, $\Phi$,  $\L$, $\U$, $\G$ are five
 $[0,\infty)-$valued,  $\sB [0,T] \oti  \cF_T-$measurable processes satisfying
 $    \int_0^T \n  ( \Phi_t \ve   \L^2_t \ve \G^q_t ) dt   \ins  L^\infty_+(\cF_T)$
 and $ E \big[    ( \int_0^T \n  g_s ds  )^p
     \+   \int_0^T \n  \U_s ds \big] \< \infty $.
 Then for some constant $\fC$   depending on  $T$,  $\nu(\cX)$, $p$,
    $C_\Phi    $,    $C_{\L}  $    and $C_\G \df  \big\| \int_0^T \G^q_t dt \big\|_{L^\infty_+(\cF_T) }$,
 \bea
   && \hspace{-2cm}  E\bigg[ \, \underset{s \in [t , T]}{\sup}   | Y^1_s \- Y^2_s|^p
  \+  \Big(\int_t^T   \big| Z^1_s \- Z^2_s \big|^2 ds \Big)^{\frac{p}{2}}
   \+     \int_t^T   \n  \int_\cX \n
    |U^1_s(x) \- U^2_s(x)|^p \nu(dx) ds   \bigg] \nonumber \\
 &&   \ls   \fC  \,   E \bigg[  |\xi_1 \- \xi_2|^p  \+  \Big( \int_t^T \n  g_s ds \Big)^p
     \+   \int_t^T \n  \U_s ds \bigg]    , \q \fa t \in [0,T]   .   \label{eq:b143}
 \eea

 \end{lem}

 This basic inequality gives rise to an  {\it a priori} estimate and a stability result of $\hL^p-$solutions of  BSDEJs,
 both of which will play important roles in the demonstration of Theorem \ref{thm_BSDEJ1}.

 \begin{prop} \label{prop-a-priori}
 Given $  p  \ins  (1,  2)  $, let $\xi \ins L^p (\cF_T)$
 and     $  f    $   be  a $p-$generator  satisfying \(H3\)$-$\(H5\)  and
   $\int_0^T \n |f(t,0,0,0)| dt          \in \n L^p_+ (\cF_T)$.
 If  $(Y , Z , U )  \ins  \hD^p  \ti  \hZ^2_{\tiny \rm loc}   \ti
     \hU^p_{\tiny \rm loc}  $ solves    BSDEJ\,$(\xi ,  f )$, then
  \bea
  \|Y\|^p_{\hD^p} + \big\| Z \big\|^p_{\hZ^{2,p}} + \big\| U \big\|^p_{\hU^p}
   \ls    \cC  E \bigg[  1 \+  |\xi  |^p
   \+  \Big( \int_0^T \n   |f(t,0,0,0)|    dt \Big)^p \, \bigg] < \infty
   \q  \label{eqn-apropri}
  \eea
 for some constant $\cC$   depending on $T$, $\nu(\cX)$, $p$,   $ \ol{C} $  and $C_\beta$.

\end{prop}

\begin{prop}   \label{prop_stab}
 Given  $p \ins  (1,2)$,    let   $\{\xi_n\}_{n \in \hN}$ be a Cauchy sequence 
  in $L^p(\cF_T)$.
  For each $n  \ins  \hN$,  let $ f_n   $   be a $p-$generator  and let
  $(Y^n, Z^n, U^n)  \ins  \hS^p $ be a solution of BSDEJ\,$(\xi_n,  f_n)$.
  Assume that   for any $m,n  \ins  \hN$ with $m  \>  n$,
   $(Y^{m,n},Z^{m,n}, U^{m,n})  \df  (Y^m \- Y^n, Z^m \- Z^n,  U^m \- U^n)$ satisfies that \dsp
  \bea
      && \hspace{-1cm}  |Y^{m,n}_s|^{p-1}
     \big\lan  \sD ( Y^{m,n}_s )  ,  f_m \big(s,Y^m_s, Z^m_s, U^m_s \big) - f_n\big(s,Y^n_s,
  Z^n_s, U^n_s \big) \big\ran    \nonumber    \\
  &&   \le  \l(s) \, \th \big(|Y^{m,n}_s|^p  + \eta_n \big) \+  \Phi_s      |Y^{m,n}_s|^p
   \+  |Y^{m,n}_s|^{p-1} \Big[ \, \L_s |Z^{m,n}_s|   \+ c(s) \big\| U^{m,n}_s \big\|_{L^p_\nu}    \Big]
     \+  \U^{m,n}_s,  \qq   \label{nonlip-cond3}
 \eea
 where

  \no \(i\) $\l(\cd) \ins L^1_+ [0,T]$ and
   $\th: [0, \infty) \to [0, \infty)$ is an increasing  concave function
   satisfying $ \int_{0+}^1 \frac{1}{ \th(t) } dt = \infty$;

  \no \(ii\)   $c(\cd)  \ins  L^q_+ [0,T]$
 and   $ \Phi $, $\L$  are two  $[0,\infty)-$valued,  $\sB [0,T] \oti  \cF_T-$measurable processes  with
 $    \int_0^T \n  ( \Phi_t \ve \L_t^2  ) dt   \ins  L^\infty_+(\cF_T)$;

  \no \(iii\)   $ \eta_n \ins L^1_+ (\cF_T) $  and
   $ \U^{m,n}    $ is a  $[0,\infty)-$valued, $\sB [0,T] \oti  \cF_T-$measurable process    such that
  \bea \label{eqn-b215}
    \lmt{ n \to \infty} E [\eta_n ] \=
  \lmt{ n \to \infty}  \; \underset{m>n}{\sup} \, E \n \int_0^T \U^{m,n}_t dt   = 0  .
  \eea
  If   $\int_0^T \l(t) dt  \> 0$, we further assume that
  \bea \label{eqn-b210}
     \underset{n \in \hN}{\sup} \Big( \| Y^n\|_{\hD^p}
   \+    \|Z^n\|_{\hZ^{2,p}}  \+  \| U^n\|_{\hU^p} \Big)    \<  \infty    .
   \eea
  Then $\big\{(Y^n, Z^n, U^n)\big\}_{n \in \hN} $ is a Cauchy sequence in $ \hS^p   $.
 \end{prop}

  The following result shows that  a BSDEJ with bounded terminal data has
  a   solution, which will also play a key role in the proof of   Theorem \ref{thm_BSDEJ1}.

\begin{prop}  \label{prop_exist_bdd}
 Given $p \ins (1,2) $, let $\xi \ins L^\infty(\cF_T)$ and $f$ be a $p-$generator satisfying  \(H1\),
 \(H3\)$-$\(H6\) and that

 \ss \no {\bf (H2')} For some $ \k_0 \ins (0,\infty)  $,
 it holds for $\dtp$   $(t,\o)  \ins  [0,T]  \ti  \O$ that
 \beas
 \big| f(t, \o, y, 0, 0) \- f(t, \o, 0, 0, 0) \big| \ls   \k_0 \big( 1 \+ |y| \big)  , \q \fa  y \ins  \hR^l  .
 \eeas
 If $\int_0^T \n  |f(t,0,0,0)| dt   \in \n  L^\infty_+(\cF_T)$,
 then      the BSDEJ\,$(\xi,f)$ has a solution
    $(Y,Z,U)  \ins  \hD^\infty   \ti   \hZ^{2,p}      \ti  \hU^p $.
\end{prop}

\section{Wellposedness with Lipschitz Generators}

\label{sec:Lip}

 When the $p-$generator   is   Lipschitz continuous in $(y,z,u)$,
 the condition (H1) is not necessary to derive a unique solution for   the corresponding BSDE with jump.
 We shall demonstrate this using  a fixed-point argument as well as Theorem \ref{thm_BSDEJ1}.

\begin{thm} \label{thm_BSDEJ2}

Given $p \in (1,2)$,   let
$  \xi   \ins  L^p(\cF_T)    $
and let  $f$ be a $p-$generator with $ \int_0^T \n |f(t,0,0,0)| dt \ins L^p_+ (\cF_T)$.
If there exists
 two    $[0,\infty)-$valued, $\sB [0,T]    \oti       \cF_T-$measurable  processes $\wt{\beta}$,
 $\L$ with
 $    \int_0^T \n   ( \wt{\beta}^q_t \ve \L^2_t)  dt   \ins  L^\infty_+(\cF_T)$
 such that   for   \dtp ~ $(t,\o)  \ins  [0,T]  \ti  \O$
 \bea
 \hspace{-5mm}
  \big| f(t,\o,y_1,z_1 ,u_1) \-  f(t,\o,y_2,z_2 ,u_2)\big|
  & \tn \le  & \tn  \wt{\beta} (t , \o) \big( |y_1 \- y_2|   \+     \|  u_1  \- u_2 \|_{L^p_\nu} \big)
    \nonumber \\
   & \tn  & \tn      +     \L (t , \o)  |z_1 \- z_2| , \q
 \fa   (y_i,z_i,u_i)     \ins  \hR^l  \ti  \hR^{l \times d}   \ti    L^p_\nu , \; i  \= 1,2 . \q
 \label{eq:f173}
 \eea
Then    BSDEJ $(\xi,f)$  admits a unique solution $(Y,Z,U )
   \ins     \hS^p    $.

\end{thm}

 As a consequence of Theorem \ref{thm_BSDEJ1} and Theorem \ref{thm_BSDEJ2},
 we   have the following result on BSDEJs whose generator $f $ is null after some stopping time $\tau$.

\begin{cor} \label{cor_f_tau}
 Given $p \ins (1,2) $,   let   $f$ be a $p-$generator with
  $\int_0^T \n |f(t,0,0,0)| dt \ins L^p_+ (\cF_T)$
 such that   either  \(H1\)$-$\(H6\) or \eqref{eq:f173} holds.
 For any $\tau \ins \cT$ and $\xi \ins L^p (\cF_\tau)$,
 the unique solution $  \big( Y  , Z  , U  \big) $
 of the BSDEJ\,$(\xi,f_\tau)$ in $\hS^p$ satisfies that
   $P \big\{Y_t \= Y_{\tau \land t}, \, t \ins [0,T] \big\} \= 1$ and   that
 $\big(Z_t,U_t\big) \= \b1_{\{t \le \tau\}} \big(Z_t,U_t\big)$, \dtp

\end{cor}

   We can even  represent a Lipschitz $p-$generator as the \pas ~ limit of the difference quotients of
 $\hL^p $ solutions to the corresponding BSDEJs:

\begin{prop} \label{gen_repre}
 Given $p \ins (1,2) $ and $\fc \> 0$,  let   $f$ be a $p-$generator with
  $\int_0^T \n |f(t,0,0,0)| dt \ins L^p_+ (\cF_T)$  such that
 \bea
  \big| f(t,  y_1,z_1 ,u_1) \-  f(t,  y_2,z_2 ,u_2)\big|
     \ls      \wt{\beta}_t      |y_1 \- y_2|
     \+      \fc \big(   |z_1 \- z_2|  \+     \|  u_1  \- u_2 \|_{L^p_\nu} \big) , ~
 \fa   (y_i,z_i,u_i)     \ins  \hR^l  \ti  \hR^{l \times d}   \ti    L^p_\nu , \; i  \= 1,2   \q
 \label{eq:f173c}
 \eea
 holds \dtp ~ for some     $[0,\infty)-$valued, $\sB [0,T]       \otimes          \cF_T-$measurable  process  $\wt{\beta}$
   with  $    \int_0^T \n    \wt{\beta}^q_t    dt   \ins  L^\infty_+(\cF_T)$.
 Let $ (  t,y,z,u)  \ins    [0,T) \ti  \hR^l \ti  \hR^{l \times d} \ti L^p_\nu$  such that
 \bea \label{eq:x231}
 \lmt{s \to t+} f(s,y,z,u) \= f(t,y,z,u) , ~ \pas  \hb{ and }
  E \Big[ \, \underset{s \in [t,t+\d ]}{\sup} |f(s,y,0,0)|^p \Big] \< \infty
  \hb{ for some }  \d \= \d (t,y ) \ins (0,T \- t] .
 \eea
 For any $s \ins (t,T]$, let $(Y^{s,y,z,u},Z^{s,y,z,u},U^{s,y,z,u})$ denote  the unique $\hL^p$  solution to
 BSDEJ\,$\big(y\+V(t,s,z,u), f_s \big)$
 with $V (t,s,z,u) \df z(B_s  \- B_t)  \+ \int_{r \in (t,s]} \n \int_\cX \n u(x) \tnp (dr,dx) \ins \cF_s$
 and  $f_s  (r,\o,y,z,u)  \df   \b1_{\{r \le s\}}f(r,\o,y,z,u) $, $ \fa (r,\o,y,z,u)
\ins [0, T]  \ti  \O  \ti  \hR^l  \ti  \hR^{l \times d}  \ti   L^p_\nu$.
   Then it holds \pas ~ that  $  f(t,y,z,u) \= \underset{\e \to 0+}{\lim} \frac{1}{\e}
  \big(Y^{t+\e,y,z,u}_t    \- y \big)  $.

\end{prop}

\section{Proofs}

\label{sec:proof}


 \no {\bf Proof of Lemma \ref{lem_stoch_integr_Lp}:}
 {\bf 1)} Let   $U  \ins  \hU^p $. Given $\o \ins \O$,
we denote the countable    set $D_{\fp(\o)}$ by $\{t_i(\o)\}_{i \in \hN}$. For any $j \ins \hN$,
  Lemma \ref{lem_lp_esti} shows that
\beas
\hspace{-3mm}
   \Big(   \sum^j_{i=1}   \big|U \big(t_i (\o),\o, \fp_{t_i (\o)}(\o)\big)\big|^2 \Big)^{\n \frac{p}{2}}
  \ls  \sum^j_{i=1}   \big|U \big(t_i (\o),\o, \fp_{t_i (\o)}(\o)\big)\big|^p
     \ls \underset{t \in D_{\fp(\o)}}{\sum} \dn \big|U \big(t,\o,\fp_t(\o)\big)\big|^p
 \= \Big( \int_{(0,T]} \n \int_\cX \n    |U_t(x)|^p N_\fp(dt,dx) \Big)  (\o) .
\eeas
Letting $j \n \to \n \infty$ on the left-hand-side yields that
\bea \label{eq:a017}
\hspace{-5mm}
\Big( \int_{(0,T]} \n \int_\cX \n    |U_t(x)|^2 N_\fp(dt,dx) \Big)^{\n \frac{p}{2}} (\o)
    \=   \bigg( \underset{t \in D_{\fp(\o)}}{\sum} \dn \big|U \big(t,\o,\fp_t(\o)\big)\big|^2 \bigg)^{\n \frac{p}{2}}    \ls   \Big( \int_{(0,T]} \n \int_\cX \n    |U_t(x)|^p N_\fp(dt,dx) \Big)  (\o) .
\eea
It follows that
 \bea \label{eq:a011}
  E \bigg[ \Big( \int_{(0,T]} \n \int_\cX \n  |U_t(x)|^2 N_\fp(dt,dx) \Big)^{\n \frac12} \bigg]
   \ls  1  \+  E    \bigg[ \Big( \int_{(0,T]} \n \int_\cX \n
   |U_t(x)|^2 N_\fp(dt,dx) \Big)^{\n \frac{p}{2}} \bigg]
    \ls  1  \+  E  \n   \int_0^T \dn \int_\cX \n  |U_t(x)|^p \nu(dx) dt   \<  \infty  ,  ~ \;\;
    \eea
    which implies  that  $  \int_{(0,T]} \n \int_\cX \n  |U_t(x)|^2 N_\fp(dt,dx) \< \infty  $,   \pas

  For any $k,n \ins  \hN$ with $k  \>  n $,
 since  $\big[ M^{U,k} \dn - \n M^{U,n}, M^{U,k} \dn - \n M^{U,n} \big]_T
 \=  \int_0^T \n \int_\cX \n  \b1_{\{ n < |U_s(x)|\le k\}} |U_s(x)|^2 N_\fp(ds,dx)$,
 one has
 \beas
 \underset{k \ge n}{\sup} \,  E \Big\{  \big[ M^{U,k} \- M^{U,n}, M^{U,k} \- M^{U,n}
  \big]^{\frac12}_T \Big\}   \ls
   E \bigg[ \Big( \int_{(0,T]} \n \int_\cX \n  \b1_{\{   |U_t(x)| > n \}}
   |U_t(x)|^2 N_\fp(dt,dx) \Big)^{\n \frac12} \bigg] .
 \eeas
 As $n \n \to \n \infty$,    \eqref{eq:a011} and the monotone convergence theorem show  that
 $ \{M^{U,n}\}_{n \in \hN} $ is a Cauchy sequence in $\big( \hM^1 , \|\cd \|_{\hM^1}  \big)$.
 Let $M^U$ be its limit.

 \no  {\bf 2)}   By Kunita-Watanabe inequality,
\beas
 && \hspace{-0.7cm}
 \big| [M^{U,n}, M^{U,n}]_t  \- [M^U, M^U]_t  \big|
  \=  \big| [M^{U,n} \- M^U, M^{U,n} \- M^U]_t  \-  2 [M^{U,n} \- M^U, M^{U,n}]_t \big|   \\
&&  \ls  [M^{U,n} \- M^U, M^{U,n} \- M^U]_t
  \+  2 \Big(   [M^{U,n} \- M^U, M^{U,n} \- M^U]_t   \Big)^{\n \frac12}
 \Big(   [ M^{U,n},  M^{U,n}]_t   \Big)^{\n \frac12} \\
 &&  \=  [M^{U,n} \- M^U, M^{U,n} \- M^U]_t
  \+  2 \Big(   [M^{U,n} \- M^U, M^{U,n} \- M^U]_t   \Big)^{\n \frac12}
 \Big( \int_{(0,t]} \n  \int_\cX \b1_{\{|U_s(x)| \le n\}} |U_s(x)|^2  N_\fp(ds,dx) \Big)^{\n \frac12} ,
 ~ \fa  t \ins [0,T] .
\eeas
 Then Lemma \ref{lem_lp_esti} and H\"older's inequality imply that
\beas
&& \hspace{-1cm}
E \bigg[ \, \underset{t \in [0,T]}{\sup} \big| [M^{U,n}, M^{U,n}]_t \- [M^U, M^U]_t  \big|^{\frac12} \bigg] \\
&& \ls  E \Big\{ \big[M^{U,n} \- M^U, M^{U,n} \- M^U \big]^{\frac12}_T \Big\}
 \+  \sqrt{2} E \bigg[ \big(   [M^{U,n} \- M^U, M^{U,n} \- M^U]_T   \big)^{\n \frac14}
 \Big(  \int_{(0,T]} \n  \int_\cX   |U_t(x)|^2  N_\fp(dt,dx)  \Big)^{\n \frac14} \bigg]  \\
&& \ls \big\| M^{U,n} \- M^U \big\|_{\hM^1}
 \+  \sqrt{2} \big\| M^{U,n} \- M^U \big\|_{\hM^1}^{\n \frac12}
 \bigg( E \bigg[ \Big(  \int_{(0,T]} \n  \int_\cX  \n
  |U_t(x)|^2  N_\fp(dt,dx)  \Big)^{\n \frac12} \bigg] \bigg)^{\n \frac12} .
\eeas
 Letting $n \n \to \n  \infty$ yields that $ \lmt{n \to \infty}
 E \bigg[ \, \underset{t \in [0,T]}{\sup} \big| [M^{U,n}, M^{U,n}]_t \- [M^U, M^U]_t  \big|^{\frac12} \bigg]
 \= 0 $.
 So there exists a subsequence of $ \{M^{U,n}\}_{n \in \hN} $
 \big(we still denote it by $\{M^{U,n}\}_{n \in \hN}$\big) such that
 $ \lmt{n \to \infty} \, \underset{t \in [0,T]}{\sup} \big| [M^{U,n}, M^{U,n}]_t \- [M^U, M^U]_t  \big|
 \= 0 $, \pas , which together with   the monotone convergence theorem yields that for \pas~ $\o \ins \O$
  \beas
[M^U, M^U]_t (\o) & \tn =  & \tn  \lmt{n \to \infty} [M^{U,n}, M^{U,n}]_t (\o)
 \=  \lmtu{n \to \infty} \Big( \int_{(0,t]} \n \int_\cX \n \b1_{\{|U_s(x)| \le n\}} |U_s(x)|^2  N_\fp(ds,dx) \Big)(\o)\\
 & \tn  =  & \tn   \lmtu{n \to \infty}  \underset{s \in D_{\fp(\o)} \cap (0,t]}{\sum} \dn  \b1_{\{|U  (s,\o,\fp_s(\o) )|
 \le n\}}   \big|U \big(s,\o,\fp_s(\o)\big)\big|^2
  \n = \dn    \underset{s \in D_{\fp(\o)} \cap (0,t]}{\sum} \dn \big|U \big(s,\o,\fp_s(\o)\big)\big|^2 \\
  & \tn  =  & \tn  \Big( \int_{(0,t]} \int_\cX   |U_s(x)|^2  N_\fp(ds,dx) \Big)(\o) , ~ \; \fa t  \ins  [0,T].
\eeas

 Then the Burkholder-Davis-Gundy inequality and \eqref{eq:a011} show that
 \beas
  E   \bigg[ \underset{t \in [0, T]  }{\sup} \big|M^U_t\big|^p  \bigg]
  \le c_{p,l} \Big[ \big[M^U, M^U\big]^{\frac{p}{2}}_T  \Big]
  = c_{p,l} E   \bigg[ \Big( \int_{(0,T]} \n \int_\cX
 |U_t(x)|^2  N_\fp(dt,dx) \Big)^{\frac{p}{2}}  \bigg]  \< \infty ,
 \eeas
 which implies that $ M^U $ is a uniformly integrable martingale.

  \no {\bf 3)} As $\|\cd \|_{\hM^1}$    is equivalent to
 $\|\cd\|_{\hD^1} $ on $\hM^1  $, we see that $\lmt{n \to \infty} E \bigg[ \underset{t \in [0,T]}{\sup}
 \big| M^{U,n}_t \- M^U_t \big| \bigg] = 0$.  So there exists a subsequence of $ \{M^{U,n}\}_{n \in \hN} $
 \big(we still denote it by $\{M^{U,n}\}_{n \in \hN}$\big) such that
 $ \lmt{n \to \infty} \, \underset{t \in [0,T]}{\sup} \big| M^{U,n}_t \- M^U_t  \big|  \= 0 $
 except on a $P-$null set $\cN$.
 We also assume that for any $\o \ins \cN^c$, the paths $ M^U(\o)  $ and $ M^{U,n} (\o)$, $n \ins \hN$
 are c\`adl\`ag.

   Let $\o \ins \cN^c$, $t \ins (0,T]$ and $\e \> 0$.
 One can find $N \= N (\o) \ins \hN$
 such that $\underset{t \in [0,T]}{\sup} \big| M^{U,n}_t \- M^U_t  \big| \< \e/2$
 for any $n \ge N$.  Also,  there   exists  $\d \= \d (t,\o) \ins (0,t)$
 such that $ \big|M^U_s (\o) \- M^U_{t-} (\o)\big|  \< \e/2$ for any $s \in (t-\d,t)$.
 Then   for any $ n \gs N$,  we have  $\big|M^{U,n}_s (\o)  \-  M^U_{t-} (\o) \big|
\le \big|M^{U,n}_s (\o)  \-  M^U_s (\o) \big| \+  \big|M^U_s (\o) \- M^U_{t-} (\o)\big| \< \e$,
$\fa s \ins (t-\d,t) $. Letting $s \nearrow t$ yields that
$ \big|M^{U,n}_{t-} (\o)  \-  M^U_{t-} (\o) \big| \ls \e $,
which  shows that $ \lmt{n \to \infty} M^{U,n}_{t-} (\o) \= M^U_{t-} (\o) $. It follows that
\beas
 \D M^U_t (\o) & \tn =  & \tn  M^U_t (\o) \- M^U_{t-} (\o)
 \= \lmt{n \to \infty} \big( M^{U,n}_t (\o) \- M^{U,n}_{t-} (\o) \big)
 \= \lmt{n \to \infty}  \D M^{U,n}_t (\o) \\
  & \tn  =  & \tn  \lmt{n \to \infty}  \b1_{\{t \in D_{\fp(\o)}\}}
 \b1_{\{|U (t,\o,\fp_t(\o))| \le n\}} U \big(t,\o,\fp_t(\o)\big)
 \=    \b1_{\{t \in D_{\fp(\o)}\}}   U \big(t,\o,\fp_t(\o)\big) .
\eeas

  \no {\bf 4)} Let $ U^1,  U^2 \ins  \hU^p $ and $ n \ins  \hN$.
 For $i \= 1,2$,   define
 \beas
 \hspace{-3mm}
  X^{i,n}_t  \df  \int_{(0,t]} \n \int_\cX \n \b1_{\{|U^1_s(x)+U^2_s(x)| \le n \}}
 U^i_s (x) \tnp (ds,dx)
 \hb{ and }      \wt{X}^{i,n}_t  \df  \int_{(0,t]} \n \int_\cX \n \b1_{\{|U^1_s(x)+U^2_s(x)| \le n, |U^i_s(x)| \le n\}}
 U^i_s (x) \tnp (ds,dx)  , ~ \; t  \ins  [0,T] .
 \eeas
 We can deduce that
\beas
&& \hspace{-1.5cm} \big\|M^{U^1+U^2 \n , \,n } \-  M^{U^1 \n ,\,n}  \-  M^{U^2 \n ,\,n} \big\|_{\hM^1}
 \=  \big\|X^{1,n} \+ X^{2,n} \-  M^{U^1 \n ,\,n}  \-  M^{U^2 \n ,\,n} \big\|_{\hM^1}
 \ls  \sum_{i=1,2} \Big( \big\|X^{i,n} \- \wt{X}^{i,n}   \big\|_{\hM^1}
\+ \big\|\wt{X}^{i,n} \-  M^{U^i \n ,\,n}    \big\|_{\hM^1} \Big) \\
&& =   \sum_{i=1,2}   E \bigg[ \Big( \int_{(0,T]} \n \int_\cX \n
  \b1_{\{|U^1_t(x)+U^2_t(x)| \le n, |U^i_t(x)| > n\}} |U^i_t (x)|^2 N_\fp (dt,dx) \Big)^{\n \frac12} \bigg] \\
&& \q    +   \sum_{i=1,2}   E \bigg[ \Big( \int_{(0,T]} \n \int_\cX \n
  \b1_{\{|U^1_t(x)+U^2_t(x)| > n, |U^i_t(x)| \le n\}} |U^i_t (x)|^2 N_\fp (dt,dx) \Big)^{\n \frac12} \bigg] \\
&& \le \sum_{i=1,2}  E \bigg[ \Big( \int_{(0,T]} \n \int_\cX \n
  \b1_{\{  |U^i_t(x)| > n\}} |U^i_t (x)|^2 N_\fp (dt,dx) \Big)^{\n \frac12}
   \+   \Big( \int_{(0,T]} \n \int_\cX \n
  \b1_{\{|U^1_t(x)+U^2_t(x)| > n \}} |U^i_t (x)|^2 N_\fp (dt,dx) \Big)^{\n \frac12} \bigg] .
\eeas
 As $n \n \to \n \infty$,   \eqref{eq:a011} and the monotone convergence theorem show  that
 $ \lmt{n \to \infty} \big\|M^{U^1+U^2 \n , \,n } \-  M^{U^1 \n , \,n }  \-  M^{U^2 \n , \,n } \big\|_{\hM^1}
 \= 0$, which implies that   $M^{U^1+U^2}  \=  M^{U^1 }  \+  M^{U^2 } $.

   Next,   let   $ U  \ins  \hU^p $, $\a \ins  \hR$ and $ n \ins  \hN$. One has
  \beas
  \big\|M^{\a U,n } \- \a M^{U ,n}   \big\|_{\hM^1}
 & \tn =  & \tn     E \bigg[ \Big( \int_{(0,T]} \n \int_\cX \n
  \b1_{\{(1 \land |\a|^{-1}) n < |U_s(x) | \le (1 \vee |\a|^{-1}) n \}} |\a U_s (x)|^2 N_\fp (dt,dx) \Big)^{\n \frac12} \bigg] \\
 & \tn  \le  & \tn |\a|   E \bigg[ \Big( \int_{(0,T]} \n \int_\cX \n
  \b1_{\{   |U_s(x) |  > (1 \land |\a|^{-1}) n \}} |  U_s (x)|^2 N_\fp (dt,dx) \Big)^{\n \frac12} \bigg] .
\eeas
 Letting $n \n \to \n \infty$, using  \eqref{eq:a011} and the monotone convergence theorem again yield  that
 $ \lmt{n \to \infty} \big\|M^{\a U,n } \- \a M^{U,n } \big\|_{\hM^1} \= 0$,
 which implies that   $M^{\a U}  \= \a M^{U } $.
 Therefore $U  \n \to \n  M^U$ is a linear   mapping on $\hU^p$. \qed

  \no {\bf Proof of Remark \ref{rem_drift}:  }   Let $ (Y,Z,U) \ins \hD^1 \ti \hZ^2_{\tiny \rm loc} \ti \hU^p_{\tiny \rm loc}$.
 Fix $n,k \ins \hN$. Define
 \beas 
     \tau_n  \df  \inf\bigg\{ t \in [0, T]:    \int_0^t \n |f (s, 0, 0, 0 )| ds  \+
  \int_0^t  \n  \big|Z_s\big|^2      ds
      \+  \int_0^t  \n \int_\cX   |U_s (x)  |^p  \nu(dx) ds   \> n \bigg\}  \ld  T  \ins  \cT
  \eeas
  and $A_k \df \{Y_* \ls k\} \ins \cF_T$.

  Since $|Y_t| \ls k$, $\fa t \ins [0,T]$ on $A_k$,   (H2), (H4), (H5) and  H\"older's inequality imply that
 \beas
   && \hspace{-1.5cm} E   \Big[ \b1_{A_k} \n  \int_0^{\tau_n} \n   \big|f (t, Y_t, Z_t, U_t )\big|   dt \Big]
       \ls       E \dn  \int_0^{\tau_n} \dn
 \big(  |f (t, 0, 0, 0 )|   \+  \phi^k_t \+ \beta_t
   \+  c_1(t) |Z_t|   \+ c_2(t)   \|U_t\|_{L^p_\nu}      \big) dt \\
 && \hspace{-0.3cm}  \ls  n \+ C_\beta \+    E \n     \int_0^T \n   \phi^k_t   dt   \+
     \Big( E \n \int_0^{\tau_n} \n    \big(c_1(t)\big)^2   dt \Big)^{\frac12} \Big( E  \n \int_0^{\tau_n} \n
       |Z_t |^2   dt \Big)^{\frac12}
  \+    \Big( E \n  \int_0^{\tau_n} \n    \big( c_2 (t)\big)^q   dt \Big)^{\frac{1}{q}}
  \Big(  E \int_0^{\tau_n} \n   \|U_t\|^p_{L^p_\nu}    dt \Big)^{\frac{1}{p}}     \\
  && \hspace{-0.3cm}  \ls  n \+ C_\beta \+     E \n     \int_0^T \n  \phi^k_t   dt        \+
     \sqrt{n} \Big(   \int_0^T \n    \big(c_1(t)\big)^2   dt \Big)^{\frac12}
    \+   n^{\frac{1}{p}}   \Big(    \int_0^T \n
  \big( c_2 (t)\big)^q   dt \Big)^{\frac{1}{q}}     \<  \infty ,
 \eeas
 which shows that  $ \b1_{A_k} \int_0^{\tau_n} \n   \big|f (t, Y_t, Z_t, U_t )\big|   dt  \<  \infty $, \pas ~
 As $Y_* \< \infty$, \pas, letting $k \n \to \n \infty$, we see that
 $   \int_0^{\tau_n} \n   \big|f (t, Y_t, Z_t, \\  U_t )\big|   dt  \<  \infty $
  except on a $P-$null set $\cN_n$.
 Since  $\int_0^T \n |f (t, 0, 0, 0 )| dt \< \infty$, \pas ~ and since
 $(Z,U) \ins  \hZ^2_{\tiny \rm loc}   \ti  \hU^p_{\tiny \rm loc}$, there exists a
 $P-$null set $\cN_0$ such that for any $\o \ins \cN^c_0$,
  $ \tau_\fn(\o) \= T $ for some   $ \fn  \=  \fn (\o)  \ins  \hN$.
  Now, for any $ \o  \ins      \underset{n \in    \hN \cup \{0\}}{\cap} \cN^c_n$, one  can deduce that
 $ \int_0^T  \n   \big|f (t,\o, Y_t (\o) , Z_t (\o), U_t (\o) )\big|   dt
  \= \int_0^{\tau_\fn  (\o)}  \n   \big|f (t,\o, Y_t (\o) , Z_t (\o), U_t (\o) )\big|   dt  \<  \infty $.
   \qed

 \no {\bf Proof of Lemma \ref{lem-a-priori}:}
  Set   $\wp \df  \big(2^{p-4} \, p(p \- 1)\big)^{\frac{1}{p}}$    and   define processes
 \beas
 a_t :=   \Phi_t  \+    \frac{\L^2_t}{p-1}
   +   \frac{p-1}{p} \wp^{-q} \G^q_s  + \frac{1}{p} \wp^p \nu(\cX)
  \q \hb{and} \q           A_t  \df  p \n \int_0^t  a_s   ds  , \q
    t  \ins  [0,T] \, .
 \eeas
 Then  $ C_A  \df  \|A_T\|_{L^\infty_+(\cF_T) }
  \ls  p \, C_\Phi  \n + q C_{\L}
   \+   (p \- 1) \wp^{-q}    C_\G
   \+ \wp^p \nu(\cX) T   $.
   In this proof, we let   $\fC$ denote a generic constant   depending on
    $T$, $\nu(\cX)$, $p$,     $C_\Phi$, $C_{\L}$ and $C_\G$,
    whose form may vary from line to line.

   \no {\bf  1)}  {\it Denote $(Y,Z,U) \df (Y^1 \- Y^2, Z^1 \- Z^2,U^1 \- U^2 )$.
   We first apply It\^o's formula to derive the dynamics of the approximate $p-$th power of process $Y$:
   $\vf_\e( Y_t )  \df  \big( |Y_t|^2 \+  \e   \big)^{\n \frac12}$}.

 \ss  Let us    fix  $t_0 \ins [0,T]$,  $n  \ins  \hN$  and define 
       \bea    \label{sp_tau_n}
     \tau_n := \inf\bigg\{ t \in [0, T]:     \int_0^t  |Z_s|^2      ds
     + \int_0^t  \n \int_\cX   |U_s (x)  |^p  \nu(dx) ds  > n  \bigg\} \land T    \in   \cT . 
  \eea
 For any   $\e  \ins  (0,1]$,
 the function $  \vf_\e(  x)  \df  \big( |x|^2 \+  \e   \big)^{\n \frac12}$,
  $   x   \ins     \hR^l$
   has the following derivatives of its $p-$th power:
   \bea \label{eqn-b150}
      D_{i}  \, \vf^p_\e( x)  = p \, \vf^{p-2}_\e( x) \, x_i
  \q \hb{and} \q   D^2_{i j}  \, \vf^p_\e( x) = p \, \vf^{p-2}_\e( x)
  \, \d_{ij} +  p(p \- 2)\vf^{p-4}_\e( x) \, x_i x_j   , \q \fa i, j \in \{1, \cds, l \} .
  \eea
 We also set   $\fS^\e_t \=  \fS^{n,\e}_t  \df  \underset{s \in [\tau_n \land  t ,  \tau_n  ]}{\sup}
    \vf_\e(  Y_s ) $,    $t  \ins  [t_0,T]$.   By Lemma \ref{lem_lp_esti},
    \bea \label{eq:a143}
     E \Big[ \big(\fS^\e_{t_0}\big)^{p  } \Big]
    \ls   E \bigg[ \, \underset{s \in [0,T]}{\sup}  \vf^p_\e ( Y_s) \bigg]  \ls
    E \bigg[  \, \underset{t \in [0, T] }{\sup}    |Y_t|^p \bigg] \+   \e^{\frac{p}{2}}
    \=    \|Y\|^p_{\hD^p}  \+    \e^{\frac{p}{2}}    \<  \infty .
    \eea

 Now, let us fix $(t,\e) \ins [t_0,T] \ti (0,1]$.
 Applying It\^o's formula  \big(see e.g. \cite[Theorem VIII.27]{Proba_pot_2}   or \cite[Theorem II.32]{Pr-90}\big)
 to process $e^{A_s} \vf^p_\e( Y_s)$ over the interval $[\tau_n \ld t, \tau_n  ]$
 and using \eqref{Y_jump} yield  that
 \bea
   && \hspace{-1.5 cm}
   e^{A_{\tau_n \land t}} \vf^p_\e (   Y_{\tau_n \land t}) \+  \frac{1}{2} \n \int_{\tau_n \land t}^{\tau_n} \n e^{A_s}    \mbox{trace}\big( Z_s
  Z^T_s D^2   \vf^p_\e( Y_s)  \big) ds
     \n + \dn  \underset{s \in (\tau_n \land t,\tau_n] }{\sum}  e^{A_s}   \Big( \vf^p_\e( Y_s)
      \- \vf^p_\e( Y_{s-})
      \-  \big\lan D  \vf^p_\e( Y_{s-}),   \D Y_s \big\ran \n \Big)   \nonumber \\
  && \hspace{-0.8 cm} =  \n e^{A_{\tau_n}} \vf^p_\e (  Y_{\tau_n})
  \+  p \n  \int_{\tau_n \land t}^{\tau_n}  \n  e^{A_s}  \big[   \vf^{p-2}_\e( Y_s)
 \big\lan Y_s ,  f_1  (s,Y^1_s, Z^1_s, U^1_s  )
 \-  f_2  (s,Y^2_s,  Z^2_s, U^2_s  ) \big\ran \-  a_s \vf^p_\e( Y_s)\big]  ds   \nonumber \\
 & & \hspace{-0.8 cm} \q   - \, p  (M_T \- M_t  \+   \cM_T  \-   \cM_t  )  , \qq \pas,
  \label{eqn-b190}
 \eea
 where  $   M_s    \df  M^\e_s \=  \int_0^{\tau_n \land s } \n \b1_{\{r > t_0\}}
        e^{A_r}  \vf^{p-2}_\e ( Y_{r-}) \lan Y_{r-} , Z_r   dB_r \ran $
       and $  \cM_s   \df  \cM^\e_s \=     \int_{(0, \tau_n  \land   s]} \n    \int_\cX  \n
         \b1_{\{r > t_0\}} e^{A_r} \vf^{p-2}_\e ( Y_{r-})  \lan  Y_{r-}, U_r(x)\ran  \\  \tnp(dr,dx)   $,
           $\fa  s \ins [0, T]$. 
    \if{0}
    Since $M^U$ is a c\`adl\`ag local martingale, the stochastic integral
    \beas
    \int_0^t   \b1_{\{s \le \tau_n\}} \b1_{\{s > t_0\}}   e^{A_s} \vf^{p-2}_\e ( Y_{s-})   Y_{s-} d M^U_s
   =  \int_{(0,t]} \n \int_\cX  \b1_{\{s \le \tau_n\}} \b1_{\{s > t_0\}}   e^{A_s} \vf^{p-2}_\e ( Y_{s-})
   \lan  Y_{s-}, U_s(x)\ran    \tnp(ds,dx)   , \q  t \in [0,T]
    \eeas
    is well-defined.

   \fi
    Since an analogy to \eqref{eq:a017} shows that for any $ t \ins [0,T] $
    \bea \label{eq:a433}
    E \bigg[  \Big( \int_{( \tau_n \land t, \tau_n]}  \n  \int_\cX
         |U_s (x) |^2     N_\fp (  ds, dx) \Big)^{\frac{p}{2}} \bigg]
        \ls  E \n \int_{( \tau_n \land t, \tau_n]}  \n  \int_\cX   |U_s (x) |^p     N_\fp (  ds, dx)
          \=  E \n  \int_{\tau_n \land t}^{ \tau_n}  \n  \int_\cX  |U_s (x) |^p   \nu(dx) ds  \ls  n  ,
    \eea
   we can deduce from the Burkholder-Davis-Gundy inequality, Young's inequality, \eqref{eq:a143}
and \eqref{sp_tau_n}    that
   \bea
   && \hspace{-1.4cm}   E   \bigg[ \, \underset{s \in [0, T]  }{\sup}|M_s|
    \+ \underset{s \in [0, T]}{\sup}|\cM_s|    \bigg]
    \ls  c_l e^{   C_A} E   \bigg[
  \big(\fS^\e_{t_0}\big)^{p-1   } \Big( \int_0^{\tau_n}  \n
       |Z_s|^2    ds \Big)^{\n \frac12}
   \+  \big(\fS^\e_{t_0}\big)^{p-1   } \Big( \int_{(0, \tau_n]}  \n  \int_\cX
         |U_s (x) |^2     N_\fp (  ds, dx) \Big)^{\n \frac12}   \bigg]  \nonumber \\
  && \hspace{-0.7cm}  \le \n  c_{p,l} e^{   C_A}   \n
    E   \bigg[  \big(\fS^\e_{t_0}\big)^{p    } \+ \Big( \int_0^{\tau_n}  \n
       |Z_s|^2    ds \Big)^{\n \frac{p}{2}}    \+  \Big( \int_{(0, \tau_n]}  \n  \int_\cX
         |U_s (x) |^2     N_\fp (  ds, dx) \Big)^{\frac{p}{2}} \bigg]
 \ls  c_{p,l}    e^{  C_A}
 \big(  \e^{\frac{p}{2}}  \+ \|Y\|^p_{\hD^p}  \+ n^{\frac{p}{2}}     \+  n
  \big)    \<  \infty .  \qq   \label{eq:a111}
 \eea
 So both $ M $ and    $ \cM  $ are   uniformly integrable    martingales.

 \no {\bf  2)} {\it Next, we   use Taylor's    expansion and some new analytic techniques to estimate
 the jump series $\underset{s \in (\tau_n \land t,\tau_n] }{\sum}  e^{A_s}   \Big( \vf^p_\e( Y_s)
      \- \vf^p_\e( Y_{s-})
      \-  \big\lan D  \vf^p_\e( Y_{s-}),   \D Y_s \big\ran \n \Big)$ and thus the equation \eqref{eqn-b190}. }

  \ss  Given   $s  \ins  [0,T]$,  \eqref{eqn-b150} implies that
 \bea
       \mbox{trace} \big( Z_s   Z^T_s  D^2  \vf^p_\e  ( Y_s )  \big)
         & \tn =  & \tn    p\,  \vf^{p-2}_\e  (  Y_s ) |Z_s|^2
           \+ p (p \- 2)     \vf^{p-4}_\e  (  Y_s ) \cd \sum^d_{j=1}
  \n \bigg(   \sum^l_{i=1}  Y^i_s Z^{ij}_s  \bigg)^{\n 2}   \nonumber  \\
       & \tn  \ge   & \tn    p \, \vf^{p-2}_\e  (  Y_s ) |Z_s|^2 \+ p (p \- 2)     \vf^{p-4}_\e (Y_s)   |Y_s|^2 |Z_s|^2
       \gs  p  (p \- 1) \vf^{p-2}_\e (Y_s) |Z_s|^2.  \q  \label{eqn-b198}
 \eea
  Setting  $   Y^\a_s  \df  Y_{s-} \+ \a \D Y_s  $, $\a \ins [0,1]$,
  we can deduce from   Taylor's Expansion Theorem and \eqref{eqn-b150}  that
  \bea
    && \hspace{-0.8cm}  \vf^p_\e \big( Y_s  \big) \- \vf^p_\e \big( Y_{s-}  \big) \- \big\lan D
 \vf^p_\e \big( Y_{s-} \big), \D Y_s   \big\ran
   \=   \int_0^1 \n  (1 \- \a) \big\lan  \D Y_s  ,   D^2  \vf^{p}_\e \big( Y^\a_s \big)
   \D Y_s   \big\ran  d \a   \nonumber \\
    &&      =  \n       p \n  \int_0^1 \n  (1 \- \a)
 \Big[ \vf^{p-2}_\e \big( Y^\a_s \big) |\D Y_s   |^2
  \+ (p \- 2) \vf^{p-4}_\e \big( Y^\a_s \big) \big\lan \D Y_s  ,  Y^\a_s   \big\ran^2 \Big] d \a
 \gs       p  (p \- 1)
 |\D Y_s   |^2 \n \int_0^1 \n (1 \- \a) \vf^{p-2}_\e  ( Y^\a_s  )  d \a . \q \qq  \label{eq:a114}
 \eea
   When $|Y_{s-}| \ls |\D Y_s|$, one has $ \vf^{p-2}_\e(  Y^\a_s)
 \gs \big( ( |Y_{s-}| \+   \a  |\D Y_s| )^2  \+  \e   \big)^{\frac{p}{2}-1}
 \gs \big( 4  |\D Y_s|^2  \+  \e   \big)^{\frac{p}{2}-1}
  \gs  2^{p-2} \big(    |\D Y_s|^2  \+  \e   \big)^{\frac{p}{2}-1}  $, $\fa \a \ins [0,1]$.
  So it follows from \eqref{eq:a114} and \eqref{eq:b617} that   for \pas ~ $\o \ins \O$
    \bea
       &      &   \hspace{-1.5cm}   \underset{s \in (\tau_n (\o) \land t,\tau_n (\o)] }{\sum}
  e^{A_s(\o)} \Big( \vf^p_\e \big( Y_s(\o)\big) \- \vf^p_\e\big(  Y_{s-}(\o)\big) \- \big\lan D
 \vf^p_\e\big(  Y_{s-}(\o)\big), \D Y_s (\o) \big\ran\Big)   \nonumber \\
 & & \ge   2^{p-3}  p  (p \- 1) \underset{s \in   (\tau_n (\o) \land t ,\tau_n (\o)] }{\sum}
 \b1_{\{|Y_{s-} (\o)|   \le   |\D Y_s (\o)| \}}
 e^{A_s(\o)}   \big|\D Y_s (\o)  \big|^2 \big(    \big|\D Y_s (\o) \big|^2  \+  \e   \big)^{\frac{p}{2}-1}  \nonumber \\
 & & =   2^{p-3}  p  (p \- 1) \underset{s \in D_{\fp(\o)} \cap   (\tau_n (\o) \land t ,\tau_n (\o)] }{\sum}
 \b1_{\{|Y_{s-} (\o)|   \le   |U (s,\o, \fp_s (\o) )| \}}
 e^{A_s(\o)}   \big|U (s,\o, \fp_s (\o) ) \big|^2 \big(    \big|U (s,\o, \fp_s (\o) ) \big|^2  \+  \e   \big)^{\frac{p}{2}-1}     \nonumber \\
  & & =   2^{p-3}  p  (p \- 1) \bigg( \n \int_{(\tau_n   \land t , \tau_n ]}
 \n \int_\cX \n   \b1_{\{ |Y_{s-}|  \le |U_s(x)|        \}}
  e^{A_s}   \big|U_s (x)   \big|^2
  \big(   |U_s(x)|^2  \+  \e   \big)^{\frac{p}{2}-1} N_\fp (ds,dx)  \n  \bigg)  (\o)  .    \label{eq:a451}
       \eea

   Multiplying $ \Big( \frac{|Y_s|}{\vf_\e ( Y_s)} \Big)^{2-p} \ls 1  $ to \eqref{nonlip-cond2}
   and applying Young's inequality  yield  that \pas
  \beas
   && \hspace{-0.8 cm} \vf^{p-2}_\e ( Y_s) \big\lan Y_s ,  f_1  (s,Y^1_s, Z^1_s, U^1_s  )
    \-  f_2  (s,Y^2_s,  Z^2_s, U^2_s  ) \big\ran  \nonumber \\
   && \le   \vf^{p-2}_\e ( Y_s) |Y_s| \big( g_s \+   \Phi_s  |Y_s  |   \big)
    \+  \L_s \vf^{p-2}_\e ( Y_s) |Y_s|    |Z_s   |
    \+ \G_s \vf^{p-2}_\e ( Y_s)  |Y_s|  \|U_s\|_{L^p_\nu}
     \+   \U_s   \nonumber \\
 &&  \le  \n  g_s \vf^{p-1}_\e ( Y_s) +   \Phi_s  \vf^p_\e ( Y_s)
 \+ \frac{\L^2_s}{p \- 1}     \vf^{p-2}_\e ( Y_s) |Y_s|^2
 \+  \frac{p \- 1}{4} \vf^{p-2}_\e ( Y_s) |Z_s|^2
    \+  \G_s \vf^{p-1}_\e ( Y_s)  \|U_s\|_{L^p_\nu}   \+  \U_s  \nonumber \\
   && \le \n   g_s \vf^{p-1}_\e ( Y_s)  \+     \Big(\Phi_s \+  \frac{\L^2_s}{p \- 1}
   \+ \frac{p \- 1}{p} \wp^{-q} \G^q_s \Big)
   \vf^p_\e ( Y_s) \+  \frac{p \- 1}{4} \vf^{p-2}_\e ( Y_s) |Z_s|^2
 \+ \frac{1}{p} \wp^p  \|U_s   \|^p_{L^p_\nu}
  \+  \U_s  ~ \hb{ for a.e. } s  \ins  [0,T]  .      \label{eq:b131}
  \eeas
 Since
 \beas
\|U_s\|^p_{L^p_\nu} \tn & = \tn &  \int_\cX |U_s (x)|^p \nu (dx)
 \ls  \int_\cX  \n \b1_{\{ |Y_{s-} | > |U_s (x)  |      \}}   \vf^p_\e ( Y_{s-} ) \nu (dx)
 \n + \dn  \int_\cX \n  \b1_{\{ |Y_{s-} | \le |U_s (x)  |    \}} |U_s (x)|^p \nu (dx) \nonumber \\
  \tn &  \le  \tn &      \vf^p_\e ( Y_{s-} )    \nu (\cX)
 \n + \dn  \int_\cX \n  \b1_{\{ |Y_{s-} | \le |U_s (x)  |    \}} |U_s (x)|^p \nu (dx) , \q \fa s \in [0,T] ,
\eeas
 it then follows from \eqref{Y_jump} that \pas
\beas
   && \hspace{-1 cm} \vf^{p-2}_\e( Y_s) \big\lan Y_s ,  f_1  (s,Y^1_s, Z^1_s, U^1_s  )  \-  f_2  (s,Y^2_s,  Z^2_s, U^2_s  ) \big\ran  \nonumber \\
   &&  \hspace{-4mm}  \le  g_s \vf^{p-1}_\e( Y_s)  \+    a_s \vf^p_\e( Y_s) \+  \frac{p \- 1}{4} \vf^{p-2}_\e( Y_s) |Z_s|^2
      \+  \frac{1}{p} \wp^p \n \int_\cX \b1_{\{ |Y_{s-} |   \le   |U_s (x)    | \}} |U_s (x)|^p \nu (dx)
  \+  \U_s    ~ \hb{ for a.e. } s  \ins  [0,T] . \q \qq
\eeas
 Plugging this inequality  together with  \eqref{eqn-b198}, \eqref{eq:a451}   into \eqref{eqn-b190}  leads to that
 for any $t \ins  [t_0,T]$
  \bea
  \qq &&  \hspace{-2 cm} e^{A_{\tau_n   \land t }} \vf^p_\e ( Y_{\tau_n   \land t})
   \+     \frac{ p }{4}  (p \- 1)   \dn  \int_{\tau_n \land  t}^{\tau_n }  e^{A_s}\,  \vf^{p-2}_\e ( Y_s) |Z_s|^2  ds   \+  2 \wp^p      \n \int_{(\tau_n   \land t , \tau_n ]}
 \n \int_\cX \n   \b1_{\{ |Y_{s-}|  \le |U_s(x)|         \}}
  e^{A_s}   \big|U_s (x)   \big|^2
  \big(   |U_s(x)|^2  \+  \e   \big)^{\frac{p}{2}-1} N_\fp (ds,dx)
     \nonumber \\
 && \le   \eta^\e_t \+  \wp^p  \n  \int_{\tau_n   \land t}^{  \tau_n }   \int_\cX  \n
 \b1_{\{ |Y_{s-} | \le |U_s (x)    | \}} e^{A_s} |U_s (x)|^p \nu (dx) ds
   \-  p \, ( M_T \- M_t  \+   \cM_T  \-   \cM_t  ), \q  \pas ,     \q     \label{eqn-b317}
 \eea
  where $\eta^\e_t \= \eta^{n, \e}_t \df  e^{C_A} \Big( \vf^p_\e (  Y_{\tau_n})
      \+    p  \n \int_{\tau_n \land t }^{\tau_n}  \n   g_s \vf^{p-1}_\e ( Y_s)   ds
      \+  p   \n  \int_{t_0  }^T  \n   \U_s ds  \Big)    $.
 Young's inequality and \eqref{eq:a143} show that 
     \bea \label{eq:a145}
     E[\eta^\e_t] \ls  e^{C_A} \n E \bigg[\big(\fS^\e_{t_0}\big)^p
      \dn + \n  p \big(\fS^\e_{t_0}\big)^{p-1}  \dn \int_{\tau_n \land t_0}^{\tau_n}  \n   g_s    ds
       \+  p   \n  \int_{t_0  }^T  \n   \U_s ds  \bigg]
       \dn \le \n  e^{C_A} \n E \bigg[p \big(\fS^\e_{t_0}\big)^p  \+     \bigg(\int_{t_0}^T \n  g_s    ds \bigg)^p
       \dn + \n  p   \n  \int_{t_0  }^T  \n   \U_s ds  \bigg]  \<  \infty   .
     \eea


      As $M$ and $\cM$ are uniformly integrable martingales,
   taking expectation in \eqref{eqn-b317}  give that
  \bea
  \qq &&  \hspace{-2cm}  \frac{ p }{4}  (p \- 1)
  E \n \int_{\tau_n \land  t }^{\tau_n }  e^{A_s}\,  \vf^{p-2}_\e ( Y_s) |Z_s|^2  ds
  +    2 \wp^p    E  \n  \int_{\tau_n   \land t }^{  \tau_n }
 \n \int_\cX \n   \b1_{\{|Y_{s-}| \le |U_s(x)|      \}}
  e^{A_s}   \big|U_s (x)   \big|^2
  \big(   |U_s(x)|^2  \+  \e   \big)^{\frac{p}{2}-1} \nu (dx) ds        \nonumber \\
 && \le   E [\eta^\e_t] + \wp^p E \n  \int_{\tau_n   \land t }^{  \tau_n } \n \int_\cX
 \b1_{\{ |Y_{s-} | \le |U_s (x)    | \}} e^{A_s} |U_s (x)|^p \nu (dx) ds  .         \label{eq:a141}
 \eea

 \no   {\bf  3)} {\it We continue our deduction, in which the analysis of
 $\hL^p-$norm of  random field $U$ is  quite technically involved. }

 \ss   Clearly, $ \lmtu{\e \to 0} |U (s,\o,x)|^2 \big(   |U (s,\o,x)|^2   +   \e   \big)^{\frac{p}{2}-1}
 \= |U (s,\o,x)|^p $, $\fa (s,\o,x)  \ins  [0,T]  \ti  \O  \ti  \cX$,
 so the monotone convergence theorem  implies  that
 \beas
  \lmtu{\e \to 0}  E \n \int_{\tau_n   \land t }^{  \tau_n }
 \n \int_\cX \n   \b1_{\{|Y_{s-}| \le |U_s(x)|      \}}
  e^{A_s}    |U_s (x)   |^2
  \big(   |U_s(x)|^2    +    \e   \big)^{\frac{p}{2}-1} \nu (dx) ds
  \=  E \n \int_{\tau_n   \land t }^{  \tau_n }
 \n \int_\cX \n   \b1_{\{|Y_{s-}| \le |U_s(x)|      \}}
  e^{A_s}    |U_s (x)   |^p     \nu (dx) ds  .
 \eeas
 On the other hand,
  since $\eta^\e_t  \ls  \eta^1_t$, $\fa \e  \ins  (0,1]$ and
  since  $E[ \eta^1_t ] \< \infty $ by \eqref{eq:a145}, the dominated convergence theorem shows that
  $ \lmt{\e \to 0} E[ \eta^\e_t ] \=   E    \big[ \wt{\eta}_t    \big] $,
  where $\wt{\eta}_t \df e^{C_A} \big(
|Y_{\tau_n}|^p     \+     p \int_{\tau_n \land t }^{\tau_n}   g_s | Y_s |^{p-1}    ds
       \+  p   \int_{t_0  }^T   \U_s ds \big)             $.

 Then letting   $\e  \n \to \n  0$ in \eqref{eq:a141} yields that
\beas
  2 \wp^p    E \n \int_{\tau_n   \land t }^{  \tau_n }
 \n \int_\cX \n   \b1_{\{|Y_{s-}| \le |U_s(x)|      \}}
  e^{A_s}   \big|U_s (x)   \big|^p \nu (dx) ds
    \ls   E \n  \big[ \wt{\eta}_t    \big]   \+  \wp^p    E  \int_{\tau_n   \land t }^{  \tau_n }
 \n \int_\cX \n   \b1_{\{|Y_{s-}| \le |U_s(x)|      \}}
  e^{A_s}   \big|U_s (x)   \big|^p \nu (dx) ds    .
 \eeas
 As $  E  \int_{\tau_n   \land t }^{  \tau_n }
 \n \int_\cX \n   \b1_{\{|Y_{s-}| \le |U_s(x)|      \}}
  e^{A_s}   \big|U_s (x)   \big|^p \nu (dx) ds
  \le e^{C_A} E  \int_0^{  \tau_n }
 \n \int_\cX \n      \big|U_s (x)   \big|^p \nu (dx) ds \le e^{C_A} n < \infty $, we obtain that
\bea
 \wp^p    E  \n  \int_{\tau_n   \land t }^{  \tau_n }
 \n \int_\cX \n   \b1_{\{|Y_{s-}| \le |U_s(x)|      \}}
     \big|U_s (x)   \big|^p \nu (dx) ds
 \ls  \wp^p    E  \n  \int_{\tau_n   \land t }^{  \tau_n }
 \n \int_\cX \n   \b1_{\{|Y_{s-}| \le |U_s(x)|      \}}
  e^{A_s}   \big|U_s (x)   \big|^p \nu (dx) ds
   \ls     E \big[ \wt{\eta}_t    \big]       .     \q     \label{eq:a151}
 \eea

 Now, fix $\e \ins  (0,1]$ again.
 As   $ \wt{\eta}_t  \ls  \eta^\e_t$, \eqref{eq:a141} and \eqref{eq:a151} show that
 \bea \label{eq:a153}
 \frac{ p }{4}  (p \- 1)  E  \n  \int_{\tau_n \land  t}^{\tau_n }   \vf^{p-2}_\e ( Y_s) |Z_s|^2  ds
 \le \frac{ p }{4}  (p \- 1)  E  \n  \int_{\tau_n \land  t}^{\tau_n }  e^{A_s}\,  \vf^{p-2}_\e ( Y_s) |Z_s|^2  ds
 \le   E \big[ \wt{\eta}_t  \+  \eta^\e_t \big] \le 2 E \big[ \eta^\e_t \big] .
 \eea
 Also, \eqref{eqn-b317} and \eqref{eq:a151} imply  that
 \bea \hspace{-5mm}
    E \big[  (\fS^\e_t )^p   \big]  & \tn \dn \le  & \tn  \dn   E \bigg[ \underset{s \in [ \tau_n \land t , \tau_n   ]}{\sup}
   e^{ A_s  }\vf^p_\e(  Y_s ) \bigg]
  \ls    E[\eta^\e_t ]  \+    \wp^p E  \n \int_{\tau_n   \land t }^{  \tau_n   } \n \int_\cX  \n
 \b1_{\{ |Y_{s-} | \le |U_s (x)    | \}} e^{A_s} |U_s (x)|^p \nu (dx) ds
       +  \n   2p \, E\bigg[  \underset{s \in [t , T]}{\sup}| M_s|
   \+ \underset{s \in [t , T]}{\sup}|\cM_s |  \bigg] \nonumber \\
    & \tn  \dn \le & \tn \dn   2 E \big[\eta^\e_t   \big]
     +  \n   2p \, E\bigg[  \underset{s \in [t , T]}{\sup}| M_s|
   \+ \underset{s \in [t , T]}{\sup}|\cM_s |  \bigg]  .    \label{eq:a161}
 \eea

 Similar to \eqref{eq:a111}, one can deduce from   the Burkholder-Davis-Gundy inequality, Young's inequality,
 \eqref{Y_jump}, \eqref{eq:a433} and \eqref{eq:a153} that
 \bea
    && \hspace{-1.2 cm}  2p E\bigg[  \underset{s \in [t , T]}{\sup}| M_s|
   \+ \underset{s \in [t , T]}{\sup}|\cM_s |  \bigg]
 \ls   c_l p  e^{C_A}   E \Bigg[ \big(\fS^\e_t\big)^{\frac{p}{2}}
     \bigg(  \n \int_{\tau_n \land t }^{\tau_n  }
   \vf^{p-2}_\e ( Y_{s-})  |Z_s|^2 ds  \bigg)^{\n \frac12} \+
    \big(\fS^\e_t\big)^{p-1}   \bigg(  \n  \int_{(\tau_n \land t , \tau_n  ]} \int_\cX
    |U_s(x)|^2 N_\fp(ds, dx) \bigg)^{\n \frac12}  \Bigg] \nonumber \\
&\tn \le & \tn \dn  \frac12 E   \big[ (\fS^\e_t )^p  \big]
  \+  c_l \,  p^2 e^{2 C_A}    E \n \int_{\tau_n \land t }^{\tau_n  }   \n   \vf^{p-2}_\e ( Y_s )  |Z_s|^2 ds
     \+    c_{p,l} \,   e^{p C_A}
      E \n  \left[ \bigg(  \n  \int_{(\tau_n \land t , \tau_n  ]} \n  \int_\cX \n
    |U_s(x)|^2 N_\fp(ds, dx) \bigg)^{\n \frac{p}{2}} \right] \nonumber \\
  &\tn \le & \tn \dn  \frac12 E   \big[ (\fS^\e_t )^p  \big]  \+ \fC  E[  \eta^\e_t]
     \+  \fC        E \n  \int_{\tau_n \land t }^{ \tau_n  } \n  \int_\cX \n
    |U_s(x)|^p \nu(dx) ds      .
    \label{eqn-b325}
 \eea
 By \eqref{eq:a151} again,
 \beas
 E \n  \int_{\tau_n \land t }^{ \tau_n  } \n  \int_\cX \n |U_s(x)|^p \nu(dx) ds
 & \tn \dn \le  & \tn  \dn    E \n  \int_{\tau_n \land t  }^{ \tau_n  } \n  \int_\cX \n
   \b1_{\{ |Y_{s-}| \le |U_s(x)|   \}} |U_s(x)|^p \nu(dx) ds
   \+ E \n  \int_{\tau_n \land t }^{ \tau_n  }    \int_\cX \n
  \b1_{\{  |U_s(x)| < |Y_{s-}|    \}}  |Y_{s-} |^p \nu(dx) ds \\
   & \tn  \dn  \le  & \tn  \dn  \wp^{-p} E \big[   \wt{\eta}_t \big]
  \+    \nu(\cX)    E \n  \int_{\tau_n \land t }^{ \tau_n  } \n  |Y_{s-} |^p   ds .
\eeas
 Since    \eqref{Y_jump}  and Fubini's Theorem imply that
   $   E \n  \int_{\tau_n \land t }^{ \tau_n  } \n     |Y_{s-} |^p   ds
      \=   E \n  \int_{\tau_n \land t }^{ \tau_n  } \n     |Y_s |^p   ds
      \ls   E \n  \int_{\tau_n \land t }^{ \tau_n  } \n     \big(\fS^\e_s\big)^p   ds
        \ls    E \n  \int_t^T \n     \big(\fS^\e_s\big)^p   ds
         \=     \int_t^T \n    E\big[ (\fS^\e_s )^p \big]   ds $,
   \bea \label{eq:a159}
   E \n  \int_{\tau_n \land t }^{ \tau_n  } \n  \int_\cX \n |U_s(x)|^p \nu(dx) ds
   \le  \wp^{-p} E \big[    \eta^\e_t \big] + \nu(\cX)  \int_t^T \n    E\big[  (\fS^\e_s )^p \big]   ds  .
   \eea

\no   {\bf  4)} {\it The remaining argument   is relatively routine \(c.f. Proposition 3.2 of \cite{BH_Lp_2003}\). }

 \ss    As $E \big[ (\fS^\e_t )^p   \big]  \ls  E \big[ (\fS^\e_{t_0} )^p   \big]  \<  \infty$ by \eqref{eq:a143},
 plugging \eqref{eq:a159}  back into \eqref{eqn-b325} and \eqref{eq:a161},
 we can deduce from Lemma \ref{lem_lp_esti} and Young's inequality that
 \bea
  E \big[  (\fS^\e_t )^p   \big] & \tn \le  & \tn   \fC E \big[   \eta^\e_t \big]
  \+  \fC  \int_t^T \n    E \big[ (\fS^\e_s )^p \big]   ds
   \le   \fC  E \bigg[  \big(|Y_{\tau_n}|^2  \+  \e  \big)^{\frac{p}{2}}
   \+     \big(\fS^\e_t\big)^{p-1} \n \int_{\tau_n \land t }^{\tau_n} \n   g_s  ds
   \+       \n  \int_{t_0}^T  \n  \U_s ds \bigg]
    \+  \fC  \int_t^T \n    E \big[ (\fS^\e_s )^p \big]   ds \nonumber \\
    & \tn   \le  & \tn  \frac12 E \big[ (\fS^\e_t )^p\big] \dn + \n \fC \sJ_\e
     \+  \fC  \int_t^T \n    E \big[  (\fS^\e_s )^p \big]   ds ,
    \label{eq:a163}
 \eea
 where  $   \sJ_\e \= \sJ^n_\e \df      \e^{\frac{p}{2}}  \+  E \Big[  |Y_{ \tau_n }|^p
   \+  \big( \int_{t_0}^T \n  g_s ds \big)^p  \+   \int_{t_0}^T    \U_s ds \Big]  \< \infty  $.
 So an application of  Gronwall's inequality  shows that
    \beas
      E [ \big(\fS^\e_t\big)^p   ] \le \fC \sJ_\e e^{\fC T}        = \fC \sJ_\e ,   \q \fa t \in [t_0,T] .
 \eeas

   Then we see from    \eqref{eq:a163} and \eqref{eq:a159} that
 \beas
 \qq  && \hspace{-1cm}   E\bigg[   \underset{s \in [ \tau_n \land t_0,  \tau_n]}{\sup}   | Y_s|^p  \bigg]  \ls
  E \big[  (\fS^\e_{t_0} )^p \big]  \ls  \fC \sJ_\e   , \qq
   E \big[ \eta^\e_{t_0} \big]   \ls
    \fC E \big[ (\fS^\e_{t_0} )^p \big] \dn + \n \fC \sJ_\e
     \+  \fC  \int_{t_0}^T \n    E \big[ (\fS^\e_s )^p \big]   ds   \ls    \fC \sJ_\e  \q   \hb{and} \qq \\
 && \qq \qq  E   \int_{\tau_n \land t_0}^{ \tau_n} \n  \int_\cX \n
    |U_s(x)|^p \nu(dx) ds   \ls   \fC \sJ_\e .
 \eeas
 These inequalities together with  Young's inequality   and \eqref{eq:a153} imply that
  \bea
    E\bigg[\Big(\int_{\tau_n \land t_0}^{\tau_n}   |Z_s|^2 ds \Big)^{\frac{p}{2}}    \bigg]
 & \tn \le  & \tn     E\Bigg[ \big(\fS^\e_{t_0}\big)^{\frac{p(2-p)}{2}}
  \bigg(\int_{\tau_n \land t_0}^{\tau_n}   \vf^{p-2}_\e ( Y_s) |Z_s|^2  ds \bigg)^{\n \frac{p}{2}}\Bigg] \nonumber \\
   & \tn  \le  & \tn  \frac{2-p}{2}  E \Big[ \big(\fS^\e_{t_0}\big)^p  \Big]
   \+  \frac{ p}{2} E \int_{\tau_n \land t_0}^{\tau_n}   \vf^{p-2}_\e ( Y_s) |Z_s|^2  ds
     \ls    \fC \sJ_\e  ,       \label{eqn-b338}
 \eea
 Letting $\e \n \to \n 0$, we obtain that
 \bea
 E\bigg[ \, \underset{s \in [ \tau_n \land  t_0,  \tau_n  ]}{\sup}   | Y_s|^p
  \+  \Big( \n \int_{ \tau_n \land  t_0}^{ \tau_n  } \n  |Z_s|^2 ds \Big)^{\n \frac{p}{2}}
   \dn + \n   \int_{ \tau_n \land t_0 }^{\tau_n  }   \n  \int_\cX \n
    |U_s(x)|^p \nu(dx) ds    \bigg] \ls
      \fC E \bigg[  |Y_{ \tau_n }|^p
   \+  \Big( \int_{t_0}^T \n  g_s ds \Big)^p  \+   \int_{t_0}^T  \n  \U_s ds \bigg]  . \q \label{eq:a165}
   \eea

  As  $(Z,U) \ins  \hZ^2_{\tiny \rm loc}   \ti \hU^p_{\tiny \rm loc}$,
  it holds for all   $\o \ins  \O$ except on a $P-$null set $\cN$ that
  $ \tau_\fn (\o) \= T$ for some $ \fn \= \fn (\o)  \ins  \hN$. Then
   \beas 
   \lmt{n \to \infty} Y (\tau_n   (\o),\o ) \=  Y (T,\o) \= \xi_1 (\o) \- \xi_2 (\o) ,
   \q \fa\o \ins \cN^c.
   \eeas
   (One can alternatively show this statement as follows:   Since the compensator
   $\nu(dx) dt$ of the counting measure $N_\fp(dt,dx)$ is absolutely continuous
with respect to $dt$, $P-$almost surely process $Y$ does not have a jump at time $T$. Thus
$ \lmt{n \to \infty} Y_{\tau_n} \= Y_{T-} \= Y_T $, \pas)
 \if{0}
 \beas
 \lmtu{n \to \infty} \underset{s \in [\tau_n \land t_0, \tau_n]}{\sup}   | Y_s|^p =
 \underset{s \in [  t_0, T]}{\sup}   | Y_s|^p
 \eeas
 \beas
 \lmtu{n \to \infty} \int_{(\tau_n \land t_0, \tau_n]} \n  \int_\cX \n     |U_s(x)|^p \nu(dx) ds
 = \int_{t_0}^\infty \n  \int_\cX \n     |U_s(x)|^p \nu(dx) ds
 \eeas
 \fi
  Eventually, letting  $n  \n \to \n  \infty$
  in \eqref{eq:a165},  we can derive \eqref{eq:b143}   from the monotone convergence
  theorem and  the dominated convergence  theorem.    \qed

 \no {\bf Proof of Proposition \ref{prop-a-priori}:}
 By    (H3)$-$(H5), it holds \dsp ~ that
  \beas
  &&  \hspace{-1cm} | Y_s  |^{p-1}   \big\lan \sD (  Y_s   ) ,    f  (s,Y_s, Z_s, U_s  )   \big\ran \\
  &&   \= | Y_s  |^{p-1} \Big(  \big\lan \sD (  Y_s   ) ,    f  (s,0, 0, 0  )   \big\ran
     \+    \big\lan \sD (  Y_s   ) ,    f  (s,Y_s, 0, 0  ) \- f  (s,0, 0, 0  )   \big\ran
     \+    \big\lan \sD (  Y_s   ) ,    f  (s,Y_s, Z_s, U_s  ) \- f  (s,Y_s, 0, 0  )   \big\ran    \Big)   \\
  &&   \ls | Y_s  |^{p-1}  \big(   | f  (s,0, 0, 0  ) |
     \+ \beta_s | Y_s  |   \+ |   f  (s,Y_s, Z_s, U_s  ) \- f  (s,Y_s, 0, 0  )    |  \big)  \\
  &&    \ls     |Y_s |^{p-1} \big(   |f(s,0,0,0)| \+   \beta_s  |Y_s |
   \+  \beta_s \+  c_1(s)      |Z_s  |   \+ c_2(s) \|U_s\|_{L^p_\nu}    \big)   .
 \eeas
 Clearly, $(0,0,0)$ is the solution to the BSDEJ\,$(0,0)$,
 applying Lemma \ref{lem-a-priori} with
  $(\xi_1,f_1,Y^1,Z^1,U^1)  \= (\xi ,f ,Y   ,   Z  ,U ) $, $(\xi_2,f_2,Y^2 \n ,Z^2 \n ,U^2)
   \= (0,0,0,0,0) $
  and $(g_s, \Phi_s, \L_s, \G_s, \U_s)  \=  \big( \beta_s \+ |f(s,0,0,0)|,\beta_s, c_1(s), c_2 (s), 0  \big)$,
  $s  \ins  [0,T]$ 
  yields the inequality \eqref{eqn-apropri}. \qed

  \no {\bf Proof of Proposition \ref{prop_stab}:}
  Given $  m, n \ins  \hN$ with $m \> n$,
  we set
  \beas
   \Xi^{m,n}_t  \df    \underset{s \in [t, T]}{\sup}  |Y^{m,n}_s|^p
    \+   \Big(\int_t^T   |Z^{m,n}_s|^2 ds \Big)^{\frac{p}{2}}
      \+     \int_t^T   \n  \int_\cX \n
    |U^{m,n}_s(x)|^p \nu(dx) ds    , \q t \in [0, T]  .
    \eeas
  Applying Lemma \ref{lem-a-priori} with
  $(\xi_1,f_1,Y^1,Z^1,U^1)  \= (\xi_m,f_m,Y^m  ,  Z^m \n ,U^m) $, $(\xi_2,f_2,Y^2 \n ,Z^2 \n ,U^2)
   \= (\xi_n,f_n,Y^n \n ,Z^n \n ,U^n) $
  and $(g_s,   \G_s, \U_s)  \=  (0,    c(s), \l(s) \, \th  (|Y^{m,n}_s|^p   \+  \eta_n  )
   \+  \U^{m,n}_s )$,
  $s  \ins  [0,T]$,  we can deduce from Fubini Theorem, the concavity of $\th$ and Jensen's inequality that
  for some constant $\fC$ depending on $T$, $\nu(\cX)$, $p$,
 $C_\Phi$, $ C_\L  $ and  $ \int_0^T \n  \big(c(t)\big)^q dt $
   \beas
    E \big[\Xi^{m,n}_t\big] 
  & \tn \le& \tn   \fC   \bigg( E \big[ |\xi_m \- \xi_n  |^p  \big]
  +     \int_t^T   \l(s) E \big[ \th  ( \Xi^{m,n}_s   \+  \eta_n   ) \big] ds
  +   E  \n   \int_t^T \U^{m,n}_s ds        \bigg)  \nonumber \\
  & \tn \le& \tn   \fC
  \bigg( E \big[ |\xi_m \- \xi_n  |^p  \big]
  +    \int_t^T   \l(s) \, \th \big( E \big[ \Xi^{m,n}_s \big]   \+  E[\eta_n]  \big) ds
  +   E   \n  \int_0^T \U^{m,n}_s ds        \bigg) , \q     t \in     [0,T] .
  \eeas
 Hence,   it holds for any $n \in \hN$ and $t \in [0, T]$ that
 \bea
 \underset{m>n}{\sup} E \big[\Xi^{m,n}_t\big] \ls      \fC
  \bigg( \underset{m>n}{\sup} E \big[ |\xi_m \- \xi_n  |^p \big]
   \+   \int_t^T \n  \l(s) \, \th \Big(\underset{m>n}{\sup} E \big[ \Xi^{m,n}_s \big]
  \+  E[\eta_n]  \Big) ds   \+    \underset{m>n}{\sup} \, E \n   \int_0^T \U^{m,n}_s ds     \n  \bigg).
  \label{eqn-a300}
 \eea

       Since $\{\xi_n\}_{n \in \hN}  $ is a Cauchy sequence in
 $L^p(\cF_T)$, one has
  \bea   \label{eqn-b350}
  \lmt{ n \to  \infty} \, \underset{m>n}{\sup} E \big[ |\xi_m-\xi_n  |^p \, \big]=0 .
 \eea
 If   $\int_0^T \l(t) dt  \=  0$, then
 $ \int_0^T \n  \l(s) \, \th \Big(\underset{m>n}{\sup} E \big[ \Xi^{m,n}_s \big]
  \+   E[\eta_n]  \Big) ds \= 0$.
 Taking $t \= 0$ and letting $n  \n \to \n  \infty$ in \eqref{eqn-a300},
 we see from \eqref{eqn-b350} and \eqref{eqn-b215}  that
  \bea    \label{eqn-b270}
   \lmt{ n \to    \infty}  \; \underset{m>n}{\sup} E \left[\Xi^{m,n}_0 \right] =0.
  \eea
 On the other hand, suppose that  $\int_0^T \l(t) dt > 0$.
 Lemma \ref{lem_lp_esti} implies that
  \bea
   \underset{m>n}{\sup}  E \big[\Xi^{m,n}_s \big] & \tn \le  & \tn
    \underset{m>n}{\sup}  E \big[\Xi^{m,n}_0 \big]
  \le  \underset{m>n}{\sup} \Big( \| Y^{m,n}\|_{\hD^p}
 \n +\n    \|Z^{m,n}\|_{\hZ^{2,p}} \n +\n  \| U^{m,n}\|_{\hU^p} \Big)^p   \nonumber \\
  & \tn \le & \tn  \left\{ 2\,\underset{i \in \hN}{\sup} \, \Big( \| Y^i\|_{\hD^p}
 \n +\n    \|Z^i\|_{\hZ^{2,p}} \n +\n  \| U^i\|_{\hU^p} \Big) \right\}^p
 \n< \n \infty , \q \fa (s,n)  \in [0,T] \ti  \hN .   \label{eqn-b380}
 \eea
 Since $  \l \ins  L^1_+[0, T]$
 and since $ \underset{n \in \hN}{\sup} \, E [\eta_n ]  \<  \infty $ by \eqref{eqn-b215},
 Fatou's Lemma, the monotonicity and the continuity of $\th$
  (real$-$valued concave functions are continuous) imply that for any $t \in [0, T]$
 \beas
  \lsup{ n \to \infty}  \int_t^T  \n  \l(s) \, \th \Big(\underset{m>n}{\sup} E \big[ \Xi^{m,n}_s \big]
    \+  E[\eta_n]  \Big) ds
  \ls         \int_t^T  \n   \l(s)   \lsup{ n \to \infty}
           \th \Big(\underset{m>n}{\sup} E \big[ \Xi^{m,n}_s \big]   \+  E[\eta_n]  \Big)   ds
           \ls           \int_t^T  \n   \l(s) \, \th \Big(\lsup{ n \to \infty} \, \underset{m>n}{\sup}
            E \big[ \Xi^{m,n}_s \big]    \Big) ds  .      
  \eeas
    Letting $n \to \infty$ in \eqref{eqn-a300}, we can deduce from \eqref{eqn-b350} and \eqref{eqn-b215}   that
   \beas
  \lsup{ n \to \infty}  \; \underset{m>n}{\sup} E \big[\Xi^{m,n}_t\big]
           \le    \fC     \int_t^T   \l(s) \, \th \Big(\lsup{ n \to \infty} \, \underset{m>n}{\sup} E \big[ \Xi^{m,n}_s \big]    \Big) ds ,
           \q  t \in [0, T].
   \eeas
  As $\th \n : [0, \infty)  \n \to \n  [0, \infty) $ is an increasing concave function,
  it is easy to see that either $\th  \n \equiv \n  0$ or $\th(t) \> 0$ for any $t \> 0$.
  Moreover,  one can deduce from \eqref{eqn-b380} that the function
  $\chi(t) \df  \lsup{ n \to \infty}  \; \underset{m>n}{\sup} E \big[\Xi^{m,n}_t\big]$,
  $t  \ins  [0,T]$ is  bounded.   Then  Bihari's inequality
   (see Lemma \ref{lem:theta-fcn}) and \eqref{eqn-b380} imply that $\lmt{ n \to    \infty}  \; \underset{m>n}{\sup} E \big[\Xi^{m,n}_t\big]  \= 0 $, $\fa t \in [0,T]$.
   Therefore,  \eqref{eqn-b270} always holds, which shows that $\big\{(Y^n, Z^n, U^n) \big\}_{n \in \hN}$ is a Cauchy sequence in $ \hS^p$.   \qed

 \no {\bf Proof of Proposition \ref{prop_exist_bdd}:}
 Let us  make the following settings first:

  \no $\bullet$ Set $C_f \df \big\|\int_0^T \n |f(t,0,0,0)| dt \big\|_{L^\infty_+(\cF_T)}$,
  $C_{p,\cX} \df \big(\nu(\cX)\big)^{\frac{2-p}{2p}}$  and
   \bea \label{eq:a335}
   R  \df  2 \+ \exp \Big\{ T \+  C_f  \+ 4 C_\beta \+ 2 \int_0^T \n \big(c_1(t)\big)^2 dt
    \n + \n   C^2_{p,\cX} \n \int_0^T \n  \big(c_2(t)\big)^2 dt   \Big\}
   \ti \sqrt{\|\xi\|^2_{L^\infty(\cF_T)} \+ 5T    \n+\n   C_f / 2 \+  7 C_\beta / 2   }  . \q
   \eea
 Let $\p \n : \hR^l  \n \to \n  [0,1]$ be a smooth function such  that
   $\p (x) \= 1$ (resp.\;\,$\p (x) \= 0$) if   $|x| \ls  R \- 1 $  \big(resp.\;\,$|x|  \gs  R$\big).

   \no $\bullet$ Let $\rho \n  : \hR^{l+l\times d} \to \hR^+$ be a smooth
  function that vanishes outside the unit open ball $\cB_1(0)$ of $\hR^{l+l\times d}$ and satisfies $\int_{\hR^{l+l\times d} }
\rho(x)dx=1$.  For   any $r \ins  (0,\infty) $, we set   $\rho_r(x) \df  r^{l(1+d)}\rho(rx)$, $\fa x \in \hR^{l+l\times d} $.

  \no $\bullet$  We say that $\big\{O_i\big\}^m_{i=1}$ is   a partition of the unit closed ball $\ol{\cB}_1(0)$ of $\hR^{l+l\times d}$  if $O_i$, $i  \= 1, \cds,  m$ are simple-connected, open subsets of $\cB_1(0)$ that are pairwisely disjoint, and if  $\dis \cup^m_{i=1} \ol{O}_i  \=  \ol{\cB}_1(0)$.
  Let $    \big\{O^k_i\big\}_{i=1}^{2^k}$,    $ k  \ins  \hN$  be   partitions  of  $\ol{\cB}_1(0)$ such that $\ol{O}^k_i \= \ol{O}^{k+1}_{2i-1} \cup  \ol{O}^{k+1}_{2i}$ holds  for any $k \ins  \hN$ and $i  \=  1,\cds, 2^k$.
  For each $O^k_i$, we pick up a $(y^k_i,z^k_i)  \ins  O^k_i$ with $y^k_i  \ins  \hR^l$, and let $\big\| O^k_i \big\|$ denote the volume of $O^k_i $.

 \no {\bf 1)} {\it To apply the existing wellposedness result on $\hL^p-$solutions of BSDEJs with Lipschitz generator,
 we first  approximate the monotonic generator $f$ by  a sequence of Lipschitz generators $\{f_n\}_{n \in \hN}$
 via convolution with mollifiers $ \{ \rho_n \}_{n \in \hN} $.}

 \ss   Fix  $n \ins  \hN$ with $n \> \k_0$. For any $u \ins L^2_\nu$,
 since H\"older's inequality shows that
 $ u $ also belongs to $L^p_\nu$ with $ \| u \|_{L^p_\nu} \ls \big(\nu(\cX)\big)^{\frac{2-p}{2p}} \| u \|_{L^2_\nu}
 \= C_{p,\cX} \| u \|_{L^2_\nu} $,  we define
\beas 
   \z_n(u)   \df  \Big( \frac{n}{n \ve \| u \|_{L^p_\nu}} \Big) \, u \ins L^p_\nu .
\eeas
    Applying Lemma \ref{lem_pi2} with $(\hE,\|\cd\|,r,x,y) \= \big( L^p_\nu,\|\cd\|_{L^p_\nu},n,u_1,u_2 \big)$  yields that
 \bea
 \|\z_n (u_1) \- \z_n (u_2) \|_{L^p_\nu}   \ls 2  \| u_1  \- u_2  \|_{L^p_\nu}
  \ls 2 C_{p,\cX} \| u_1 - u_2  \|_{L^2_\nu} , \q \fa u_1,u_2 \in L^2_\nu ,
   \label{eq:b415}
 \eea
 which shows that  the mapping $\z_n \n : L^2_\nu  \n \to \n   L^p_\nu $ is  $\sB(L^2_\nu) / \sB(L^p_\nu)-$measurable.
 (Note: As the space $L^p_\nu$ may not have an inner product, one may not apply Lemma \ref{lem_pi}.)

  Since   $\dis \beta^n_t  \df  \frac{n  }{n    \vee \beta_t   \ve |f(t,0,0,0)| }
  \ins  (0,1]$, $t  \ins  [0,T]$ is an $\bF-$progressively measurable
process,   we can deduce from \eqref{eqn-d051}, \eqref{eq:b415} and the
$\sP  \oti  \sB \big( \hR^l  \big)    \otimes \n  \sB \big( \hR^{l \times d} \big)
  \oti  \sB\big(L^p_\nu   \big)/\sB(\hR^l ) -$measurability of $f$ that  the mapping
  \beas
f^0_n(t,\o,y,z,u)   :=    \beta^n  (t,\o) \p(y)f\big(t,\o,y, \pi_n(z) ,\z_n(u)\big) , \q  \fa (t,\o,y,z,u) \in [0,T] \times \O \times   \hR^l \times \hR^{l \times d} \times  L^2_\nu
\eeas
 is $\sP \oti  \sB \big( \hR^l  \big)  \oti  \sB \big( \hR^{l\times d} \big)
  \oti  \sB\big(L^2_\nu   \big)/\sB(\hR^l ) -$measurable.
  Given  $(t,\o,y,z,u)  \ins  [0,T]  \ti  \O  \ti    \hR^l
     \ti  \hR^{l \times d}  \ti   L^2_\nu$, we further define
 \beas
     f_n(t,\o,y,z,u)  :=    \big(f^0_n(t,\o,\cdot,\cdot,u)\ast \rho_{n}\big)(y,z) .
\eeas
By (H1), the continuity of   mapping  $     f(t,\o,\cd,\cd,u)$ implies
that of mapping $      f^0_n(t,\o,\cd,\cd,u)$. Hence,   $f_n(t,\o,y,z,u)$ is indeed a Riemann integral:
  \bea
     f_n(t,\o,y,z,u) &=&  \int_{|(\tilde{y},\tilde{z})| \le 1}
     f^0_n\Big(t, \o, y-\frac{1}{n}\tilde{y},
 z-\frac{1}{n}\tilde{z}  , u \Big)    \rho(\tilde{y},\tilde{z})d\tilde{y}d\tilde{z}  \label{eqn-b410}  \\
 &=& \lmt{k   \to \infty}     \sum_{i=1}^{2^k}
    f^0_n\Big(t, \o, y-\frac{1}{n}y^k_i,
 z-\frac{1}{n}z^k_i  , u \Big)    \rho(y^k_i,z^k_i) \big\| O^k_i\big\|,   \nonumber
  \eea
  from which   one can deduce that $f_n$ is also
  $\sP \oti  \sB \big( \hR^l  \big)  \oti  \sB \big( \hR^{l\times d} \big)
   \oti  \sB\big(L^2_\nu   \big)/\sB(\hR^l ) -$measurable.

   Now, set $c_n(t)  \df  n \big( 3 \+ R   \+   c_1(t)  \+  c_2 (t) \big) $, $t  \ins  [0,T]$,
  which is clearly of   $  L^2_+ [0,T]$.
  As $ \beta^n_t   \big(     \beta_t \vee |     f (t,0, 0 ,0 )  | \big) \ls n $,
    (H2'), (H4) and (H5) shows that     \dtp
 \bea \label{eq:a457}
  |f^0_n(t,y,z,u)| 
   & \tn \ls  & \tn   \beta^n_t   \p(y) \big(    |     f (t,0, 0 ,0 )  |
 \+    \k_0 (1   \+  |y| )  \+ \beta_t    \+ c_1(t)\big|\pi_n(z)\big|
 \+  c_2(t)\big\|\z_n(u)\big\|_{L^p_\nu} \big)  \nonumber \\
  & \tn   \ls  & \tn  c_n(t)   ,
  \q  \fa (y,z,u)  \ins  \hR^l  \ti  \hR^{l \times d}  \ti   L^2_\nu .
 \eea
  This   implies that  for \dtp ~ $(t,\o) \ins [0,T] \ti \O$ and   any $u \ins L^2_\nu$,
  $f_n(t,\o,\cd,\cd,u)$ is a smooth function on $\hR^l \ti \hR^{l \times d}$ via convolution.

  Let $(y_1,z_1), (y_2,z_2) \ins \hR^l \ti \hR^{l \times d}$ and set $\fy_\a  \df  \a y_1 \+ (1 \- \a) y_2$,
  $\fz_\a  \df  \a z_1 \+ (1 \- \a) z_2$, $\fa \a  \ins  (0,1)$.
  Since
  \beas
   \rho_n (y_1 \- \wt{y}, z_1 \- \wt{z} ) \- \rho_n (y_2 \- \wt{y}, z_2 \- \wt{z} )
  & \tn \=  & \tn  \rho_n (\fy_1 \- \wt{y}, \fz_1 \- \wt{z} ) \- \rho_n (\fy_0 \- \wt{y}, \fz_0 \- \wt{z} )
  \= \int_0^1 d \rho_n (\fy_\a \- \wt{y}, \fz_\a \- \wt{z} )  \\
  &\=& \int_0^1 \Big\lan \big(  y_1 \-  y_2,  z_1 \-   z_2 \big),
  \nabla \rho_n \big(  \fy_\a\- \tilde{y},  \fz_\a \- \tilde{z}  \big) \Big\ran d \a ,  \q
  \fa (\wt{y}, \wt{z}) \ins \hR^l \ti \hR^{l \times d}   ,
  \eeas
  \eqref{eq:a457} also yields  that  \dtp~
 \bea
     |f_n(t, y_1,z_1,u)-f_n(t,y_2,z_2,u)|
     &  \tn=&  \tn   \bigg|\int_{\hR^{l+l\times d}} \big( \rho_n (y_1 \- \wt{y}, z_1 \- \wt{z} ) \- \rho_n (y_2 \- \wt{y}, z_2 \- \wt{z} ) \big)  f^0_n(t,\tilde{y},\tilde{z},u) d\tilde{y}d\tilde{z}\bigg|  \nonumber  \\
  &  \tn=&  \tn  \bigg|\int_{\hR^{l+l\times d}} \Big(\int_0^1 \Big\lan \big(  y_1 \-  y_2,  z_1 \-   z_2 \big),
  \nabla \rho_n \big(  \fy_\a\- \tilde{y},  \fz_\a \- \tilde{z}
\big) \Big\ran d \a\Big)  f^0_n(t,\tilde{y},\tilde{z},u) d\tilde{y}d\tilde{z}\bigg| \nonumber  \\
 & \tn \le &  \tn    c_n(t)    \int_0^1  \int_{\hR^{l+l\times d}} \big|(y_1 \- y_2,z_1 \- z_2)\big| \n \cd \n
 \big| \nabla \rho_n \big(  \fy_\a  \- \tilde{y},   \fz_\a \- \tilde{z}
\big)\big|     d\tilde{y}d\tilde{z}d \a   \nonumber  \\
  & \tn \le& \tn      \k^n_\rho \,  c_n(t) \big(|y_1 \- y_2| \+ |z_1 \- z_2|\big) , ~
   \fa (y_1,z_1 ), (y_2,z_2 )  \ins  \hR^l  \ti  \hR^{l \times d}, ~
     \fa u   \ins      L^2_\nu   , \qq   \q  \label{eqn-b414}
 \eea
  where    $\k^n_\rho  \df  \int_{\hR^{l+l\times d}} |\nabla \rho_n (x)|d x  \<  \infty $
  is a constant determined by $\rho$ and $n$.

    On the other hand,  \eqref{eqn-b410}, (H5) and \eqref{eq:b415} imply that   \dtp~
 \beas
  \q &&  \hspace{-1.2 cm}   |f_n(t, y,z,u_1)-f_n(t,y,z,u_2)| \\
   &=& \tn \bigg| \int_{|(\tilde{y},\tilde{z})| \le 1}
    \beta^n_t \p\Big( y-\frac{1}{n}\tilde{y}\Big)  \bigg( \n  f \Big(t,  y \-  \frac{1}{n}\tilde{y},
 \pi_n\Big(z \- \frac{1}{n}\tilde{z} \Big) , \z_n(u_1) \Big) \-   f \Big(t,  y  \-  \frac{1}{n}\tilde{y},
 \pi_n\Big( z  \-  \frac{1}{n}\tilde{z} \Big) , \z_n(u_2) \Big) \n \bigg)  \rho(\tilde{y},\tilde{z})d\tilde{y}d\tilde{z} \bigg| \\
 &\le & \tn   \int_{|(\tilde{y},\tilde{z})| \le 1}
     \bigg| f \Big(t,  y \-  \frac{1}{n}\tilde{y},
 \pi_n\Big(z \- \frac{1}{n}\tilde{z} \Big) , \z_n(u_1) \Big) \-   f \Big(t,  y  \-  \frac{1}{n}\tilde{y},
 \pi_n\Big( z  \-  \frac{1}{n}\tilde{z} \Big) , \z_n(u_2) \Big)  \bigg|  \rho(\tilde{y},\tilde{z})d\tilde{y}d\tilde{z}  \\
  &\le& \tn c_2 (t) \big\| \z_n(u_1)  \- \z_n(u_2)  \|_{L^p_\nu}
   \le 2 c_2 (t) C_{p,\cX} \big\|u_1 \- u_2\|_{L^2_\nu}    ,    \q
   \fa (y,z, u_1, u_2 ) \in \hR^l \times \hR^{l \times d} \times  L^2_\nu  \times  L^2_\nu   ,
 \eeas
  which together with \eqref{eqn-b414} shows that  $f_n$ is Lipschitz continuous in $(y,z,u) \ins
  \hR^l  \ti  \hR^{l \times d}  \ti   L^2_\nu $ with $ L^2_+ [0,T] -$coefficients.

     Moreover,   \eqref{eqn-b410}, (H2') and (H4) imply that   \dtp
       \beas
       \hspace{-4mm}
         | f_n(t,0,0,0) | & \tn   \le& \tn       \int_{|(\tilde{y},\tilde{z})| \le 1} \beta^n_t \bigg|
   f\Big(t,-\frac{1}{n}\tilde{y}, \pi_n\Big( \dn - \n \frac{1}{n}\tilde{z}\Big) , 0 \Big)     \bigg|    \rho(\tilde{y},\tilde{z})d\tilde{y}d\tilde{z} \\
   & \tn    \le  & \tn       \int_{|(\tilde{y},\tilde{z})| \le 1} \dn
 \beta^n_t \bigg(  |f(t,0,0,0)| \+   \k_0   \+  \frac{\k_0}{n} \big|\tilde{y}\big| \+ \beta_t
 \+ c_1(t)\Big|\pi_n\Big( \dn - \n \frac{1}{n}\tilde{z}\Big)\Big| \bigg) \rho(\tilde{y},\tilde{z})d\tilde{y}d\tilde{z}  \\
   & \tn   \le& \tn       \int_{|(\tilde{y},\tilde{z})| \le 1}
 \Big(    n \+ \k_0 \+ 1  \+ n  \+ \frac{1}{n}c_1(t)  \Big) \rho(\tilde{y},\tilde{z})d\tilde{y}d\tilde{z}
    \=   2 n \+ \k_0 \+ 1    \+ \frac{1}{n}c_1(t) ,
 \eeas
 so $          E   \Big[ \big( \int_0^T \n  |f_n (t,0,0,0)| dt \big)^2 \Big]  \ls
           \big( (  2 n \+ \k_0 \+ 1 )T \+  \frac{1}{n} \int_0^T  \n   c_1(t)  dt \big)^2    \<  \infty$.
 Then we know from the classical wellposedness  result of BSDEJs  in $\hL^2-$case (see e.g. Lemma 2.2 of \cite{Yin_Situ_2003}) that   the  BSDEJ\,$(\xi, f_n)$ has   a unique solution $(Y^n,Z^n,U^n)
 \ins  \hD^2   \ti   \hZ^2      \ti  \hU^2$.

   \no {\bf 2)} {\it In this part, we will use  regular argument to show that
    the $\hL^2-$norms of $\{(Y^n,Z^n,U^n)\}_{n \in \hN}$ are bounded.}

  \ss    Next, we define $a_t \df 1 \+ |f(t,0,0,0)| \+
  4 \beta_t    \+  2 \big(c_1(t)\big)^2
  \+ C^2_{p,\cX} \big(c_2(t)\big)^2  $
  and $A_t \df   2  \n \int_0^t \n  a_s ds$, $t  \ins  [0, T]$.
  Clearly,  $A_T  \ins  L^\infty_+(\cF_T)$ with
  $   C_A  \df  \big\| A_T \big\|_{L^\infty_+(\cF_T)}
   \ls 2 T  \+ 2 C_f \n+ \n 8 C_\beta \n+ \n 4 \int_0^T \n \big(c_1(t)\big)^2 dt
    \n+ \n 2 C^2_{p,\cX} \n \int_0^T \n  \big(c_2(t)\big)^2 dt  \<  \infty $.

     Fix $n \ins  \hN$ with $n \> \k_0$  and fix $t  \ins  [0,T]$.
   Applying It\^o's formula to  process  $e^{A_s}|Y^n_s|^2$ over
interval $[t, T]$ and using \eqref{Y_jump} yield  that
 \bea
 && \hspace{-1.5cm} e^{A_{t}}|Y^n_{t}|^2+  \int_{t}^T   e^{A_s} |Z^n_s|^2 ds    +    \int_{(t, T]} \int_\cX   e^{A_s}      |U^n_s(x) |^2    N_\fp(ds,dx)  \nonumber \\
  & \tn = & \tn  e^{A_T}|\xi|^2
   \+  2 \int_{t}^T    e^{A_s} \big[ \big\lan  Y^n_s , f_n (s, Y^n_s, Z^n_s, U^n_s) \big\ran
  \-    a_s    |Y^n_s|^2 \big]  ds
  \-  2     \big( M_T   \- M_t   \+  \cM_T  \-  \cM_t \big)  , ~ \pas,    \qq   \label{eqn-b441}
\eea
where $M_s \df   \int_0^s  \n  e^{A_r}  \lan  Y^n_{r-}, Z^n_r \, dB_r \ran$ and $\cM_s  \df
\int_{(0, s]}  \n \int_\cX e^{A_r}  \n   \lan Y^n_{r-}, U^n_r (x)\ran   \tnp(dr,dx)  $, $\fa s  \ins  [0, T]$.

 Since (H2') and (H3)   show  that  \dtp
 \beas
   \Big\lan Y^n_s  ,  f \Big(s, Y^n_s \- \frac{1}{n}y, 0,0 \Big) \- f(s,0,0,0) \Big\ran
 & \tn  \ls  & \tn  \beta_s \Big| Y^n_s \- \frac{1}{n}y \Big|^2
 \+ \frac{1}{n} |y| \Big| f \Big(s, Y^n_s \- \frac{1}{n}y, 0,0 \Big) \- f(s,0,0,0) \Big| \\
  & \tn  \ls  & \tn  \beta_s \Big| Y^n_s \- \frac{1}{n}y \Big|^2
 \+   |y|   \Big(1 \+ \Big|Y^n_s \- \frac{1}{n}y \Big|\Big)   ,
 \q \fa  (y,z) \ins \hR^l \times \hR^{l \times d}
 \eeas
  and since $\big\|\z_n(U^n_s)\big\|_{L^p_\nu} \ls   \|U^n_s\|_{L^p_\nu} \ls C_{p,\cX} \|U^n_s\|_{L^2_\nu} $,
  we can deduce from \eqref{eqn-b410},  (H4) and (H5)    that \pas
  \bea
 && \hspace{-1cm}    \lan Y^n_s,f_n(s,Y^n_s,Z^n_s, U^n_s)\ran
    \=      \int_{|(y,z)| \le 1}
  \beta^n_s\p\Big(Y^n_s \- \frac{1}{n}y\Big) \Big\lan Y^n_s,  f \Big(s, Y^n_s \- \frac{1}{n}y,
\pi_n \Big(Z^n_s \- \frac{1}{n}z\Big),\z_n \big(U^n_s\big) \Big) \Big\ran \rho(y,z)dydz  \nonumber  \\
  & &   \le       \int_{|(y,z)| \le 1} \bigg\{  | Y^n_s|
 \left[ |f(s,0,0,0)|
  \+ \beta_s \+ c_1(s)\Big|\pi_n\Big(Z^n_s \- \frac{1}{n}z\Big)\Big|
 \+ c_2(s) \big\|\z_n(U^n_s)\big\|_{L^p_\nu}\right] \+  \nonumber \\
 &&  \q    \beta_s \Big| Y^n_s \- \frac{1}{n}y \Big|^2
 \+   |y|   \Big(1 \+ \Big|Y^n_s \- \frac{1}{n}y \Big|\Big) \bigg\} \rho(y,z)dydz     \nonumber  \\
   & &  \le   2  \+ \beta_s  \+     | Y^n_s|    \Big( 1 \+ |f(s,0,0,0)| \+
   3  \beta_s \+   c_1(s) \big(1 \+ |Z^n_s| \big)
  \+   c_2 (s) C_{p,\cX} \|U^n_s\|_{L^2_\nu} \Big)
  \+ \beta_s    | Y^n_s |^2     \nonumber  \\
   & &  \le    \frac52  \+  \frac14  |f(s,0,0,0)| \+  \frac74  \beta_s
     \+  a_s |Y^n_s|^2 \+ \frac14 |Z^n_s|^2
     \+ \frac14  \|U^n_s\|^2_{L^2_\nu}  \; \; \hb{ for a.e. } s  \ins  [0,T] \, ,  \qq \q \label{eqn-b445}
  \eea
 where we used the inequality $\a \ls \frac14 \+ \a^2 $, $\fa \a \ins [0,\infty) $.

  Moreover, Burkholder-Davis-Gundy inequality and H\"older's inequality imply  that
    \beas
        E\bigg[ \underset{s \in [0, T]  }{\sup}|M_s|   \+
        \underset{s \in [0, T]}{\sup}|\cM_s|   \bigg]
   & \tn \le   & \tn     c_l e^{ C_A}
   E \n \left[  Y^n_*  \n  \left( \int_0^T \n    |Z^n_s|^2   ds \right)^{\frac12}
   \+  Y^n_*  \n
   \left( \int_0^T  \dn  \int_\cX \n  \big|U^n_s (x)\big|^2 N_\fp (ds,dx) \right)^{\frac12} \right] \nonumber \\
  & \tn  \le  & \tn   c_l e^{ C_A} \|Y^n\|_{\hD^2} \big(\|Z^n\|_{\hZ^2}
   \dn + \n  \|U^n\|_{\hU^2}\big)  \<  \infty ,
 \eeas
    which shows that both $M  $ and $\cM$ are uniformly integrable martingales. Since
 \beas
     E \n \left[   \int_{(t, T]}  \int_\cX   e^{A_s}    |U^n_s(x) |^2    N_\fp(ds,dx)  \bigg|\cF_t\right]
   \n = E \n \left[  \int_{t}^T  \int_\cX e^{A_s} |U^n_s(x)|^2 \nu(dx) \, ds  \bigg|\cF_t\right]
   \n = E \n \left[  \int_{t}^T \n e^{A_s}  \|U^n_s\|^2_{L^2_\nu}   ds  \bigg|\cF_t\right] , ~ \pas ,
 \eeas
 taking conditional expectation $E[\cdot|\cF_t]$ in \eqref{eqn-b441}, one can deduce from \eqref{eqn-b445} that \pas
 \beas
&& \hspace{-1.5cm}
|Y^n_{ t}|^2 \+ \frac{1}{2}E \bigg[ \int_{t}^T
 \Big( |Z^n_s|^2 \+ \|U^n_s\|^2_{L^2_\nu} \Big) ds \Big|\cF_t \bigg]
 \ls   e^{A_{ t}}|Y^n_{ t}|^2
\+ \frac{1}{2}E \bigg[ \int_{t}^T e^{A_s}
 \Big( |Z^n_s|^2 \+ \|U^n_s\|^2_{L^2_\nu} \Big) ds \Big|\cF_t \bigg] \nonumber  \\
 &&      \le        e^{C_A } \Big(  \big\|\xi \big\|^2_{L^\infty(\cF_T)}  \+ 5 T    \+
   C_f / 2 \+  7 C_\beta / 2   \Big)  \ls  (R \- 2)^2 .   
 \eeas

 This together with the right-continuity of $Y^n$ implies that
 \bea \label{norm-esti}
 \| Y^n \|_{\hD^\infty} \le R-2 \q \mbox{and} \q
 \big\|Z^n\big\|^2_{\hZ^2} + \|U^n\|^2_{\hU^2} \le 2(R-2)^2, \q \fa  n  \in \hN.
 \eea

 \no {\bf 3)} {\it Next, we carefully verify conditions \eqref{nonlip-cond3} and \eqref{eqn-b215}
  for $(Y^n,Z^n,U^n)$'s, so the sequence has a limit $(Y,Z,U)$ according to   Proposition \ref{prop_stab}.}

 \ss  For any $(t,\o) \ins [0,T]  \ti  \O$ except on
 a $dt \ti  dP-$null  set $\fN $,
    we may assume that (H2'),  (H4)$-$(H6) hold, that   $|Y^n_t (\o)|  \ls  R \- 2 $,   $ \fa  n   \ins  \hN$,
 and that $U^n_t (\o) \in L^2_\nu \sb L^p_\nu $,   $ \fa  n   \ins  \hN$.

  Fix $(t, \o) \in  \fN^c$.
 By    (H5) and (H6),   it holds    for any $(y_1,z_1,u_1 ) , (y_2,z_2,u_2 ) \ins  \hR^l
   \ti  \hR^{l \times d}   \ti   L^p_\nu$ that
 \bea
 && \hspace{-1.4cm} |y_1 \- y_2|^{p-1} \big\lan \sD ( y_1 \- y_2 ) ,  f (t,\o,y_1,z_1,u_1   )
 \- f (t,\o,y_2,z_2,u_2   )   \big\ran \nonumber \\
 &&  \hspace{-0.7cm}   \le \n   |y_1 \- y_2|^{p-1} \Big( \big\lan \sD ( y_1 \- y_2 ) ,  f (t,\o,y_1,z_1,u_1   )
 \- f (t,\o,y_2,z_2,u_1   )   \big\ran
  \+    \big|  f (t,\o,y_2,z_2,u_1   ) \- f (t,\o,y_2,z_2,u_2   )\big| \Big)  \nonumber  \\
 &&  \hspace{-0.7cm}  \le \n  \l(t) \, \th \big( |y_1 \- y_2|^p \big) \+  \Phi (t,\o) |y_1 \- y_2|^p
 \+   |y_1 \- y_2|^{p-1} \big(  \L (t,\o)     |z_1 \- z_2|
 \+  c_2(t) \| u_1   \-  u_2 \|_{L^p_\nu}   \big) . \qq  \label{eq:b411}
 \eea
  Let us also fix $m, n \ins  \hN$ with $m \> n$.
  Since $(Y^m,Z^m,U^m) 
  $ is the  unique solution of   BSDEJ\,$(\xi, f_m)$ and since $ \psi (x) \n \equiv \n 1$ for all $|x| \ls R \- 1$,
       \eqref{eqn-b410} and \eqref{norm-esti} show    that
       $(Y^{m,n},Z^{m,n}, U^{m,n})  \df  (Y^m \- Y^n, Z^m \- Z^n,  U^m \- U^n)$ satisfies
    \bea
 && \hspace{-1cm}   |Y^{m,n}_t (\o)|^{p-1} \big\lan \sD  (Y^{m,n}_t (\o) ),  f_m(t,\o,Y^m_t (\o), Z^m_t (\o), U^m_t (\o))-f_n(t,\o,Y^n_t (\o),  Z^n_t (\o), U^n_t (\o)) \big\ran \nonumber \\
 &&=  \int_{|(\tilde{y},\tilde{z})| \le 1}  |Y^{m,n}_t (\o)|^{p-1} \big\lan \sD  (Y^{m,n}_t (\o) ),   \beta^m_t (\o)
 h^m_{t,\o}(\tilde{y},\tilde{z}) -\beta^n_t (\o)  h^n_{t,\o}(\tilde{y},\tilde{z}) \big\ran \rho(\tilde{y},\tilde{z})d\tilde{y}d\tilde{z} , \qq   \label{eqn-b471}
  \eea
where  $h^n_{t,\o}(\tilde{y},\tilde{z}) \df
f\big(t,\o,Y^n_t (\o) \- \frac{1}{n}\tilde{y},\pi_n(Z^n_t (\o) \- \frac{1}{n}\tilde{z}),\z_n(U^n_t (\o))\big)  $.
      Next, we fix  $( \tilde{y}, \tilde{z})  \ins  \hR^l  \ti  \hR^{l \times d}$
      with $|(\tilde{y},\tilde{z})| \< 1$ and set
   $(\wt{y}_{m,n}, \wt{z}_{m,n}) \df
    \Big(\big(\frac{1}{m} \- \frac{1}{n}\big)  \tilde{y},
    \big(\frac{1}{m} \- \frac{1}{n}\big)  \tilde{z} \Big)$.
   Consider the following  decomposition:
   \beas
   && \hspace{-0.8cm}   |Y^{m,n}_t (\o)|^{p-1}  \big\lan  \sD( Y^{m,n}_t (\o) ) ,  \beta^m_t (\o) h^m_{t,\o}(\tilde{y},\tilde{z})
 \- \beta^n_t (\o) h^n_{t,\o}(\tilde{y},\tilde{z}) \big\ran \\
 &&   =    \beta^m_t (\o) \big|Y^{m,n}_t (\o) \- \wt{y}_{m,n}
 \big|^{p-1} \big\lan \sD \left( Y^{m,n}_t (\o) \- \wt{y}_{m,n} \right),
 h^m_{t,\o}(\tilde{y},\tilde{z}) \- h^n_{t,\o}(\tilde{y},\tilde{z}) \big\ran   \\
 && \q + \beta^m_t (\o) \big\lan |Y^{m,n}_t (\o)|^{p-1} \sD\left(Y^{m,n}_t (\o)\right)
 \-  \big|Y^{m,n}_t (\o) \- \wt{y}_{m,n}
 \big|^{p-1}   \sD \left( Y^{m,n}_t (\o) \- \wt{y}_{m,n} \right) ,
h^m_{t,\o}(\tilde{y},\tilde{z})  \- h^n_{t,\o}(\tilde{y},\tilde{z}) \big\ran \\
&& \q + |Y^{m,n}_t (\o)|^{p-1}   \big\lan \sD\left(Y^{m,n}_t (\o)\right), ( \beta^m_t (\o)
  \- \beta^n_t (\o)) h^n_{t,\o}(\tilde{y},\tilde{z}) \big\ran  \df   I^1_{t,\o}(\tilde{y},\tilde{z}) \+ I^2_{t,\o}(\tilde{y},\tilde{z}) \+ I^3_{t,\o}(\tilde{y},\tilde{z})   .
 \eeas

  \no {\bf 3a)}  We see from \eqref{eq:b411} that
   \bea
&& \hspace{-1.2cm}  I^1_{t,\o}(\tilde{y},\tilde{z}) \ls  \l(t)\th(
  |Y^{m,n}_t (\o) \- \wt{y}_{m,n}|^p) \+  \Phi_t (\o)|Y^{m,n}_t (\o)  \-  \wt{y}_{m,n}|^p
   \n  +  \,    |Y^{m,n}_t (\o)  \-  \wt{y}_{m,n}|^{p-1} \ti \nonumber \\
&&  \hspace{-0.6cm}  \Big( \L_t (\o) \big|\pi_m(Z^m_t (\o)  \-  \hb{$\frac{1}{m }$} \tilde{z})
  \-  \pi_n(Z^n_t (\o)  \-  \hb{$\frac{1}{n}$} \tilde{z})\big|
   \+ c_2 (t)    \big\| \z_m \big( U^m_t (\o)\big)
     \-   \z_n \big( U^n_t (\o)\big)      \big\|_{L^p_\nu}   \Big)     .  \qq \q  \label{eqn-b501}
 \eea
 Applying Lemma \ref{lem_a02}    with
 $(b,c) \= \big(|Y^{m,n}_t (\o) \- \wt{y}_{m,n}|,|Y^{m,n}_t (\o)  | \big) $ and
 $p\=p \- 1$ (then $p\=p$)  yields that
  \bea     \label{eqn-b505}
    |Y^{m,n}_t (\o) \- \wt{y}_{m,n}|^{p-1}
    & \tn \le& \tn  \big|Y^{m,n}_t (\o)  \big|^{p-1}  \+   |\wt{y}_{m,n}|^{p-1}
     \ls   \big|Y^{m,n}_t (\o)  \big|^{p-1}  \+  n^{1-p}\\
 \hb{and} \qq  |Y^{m,n}_t (\o) \- \wt{y}_{m,n}|^p
 & \tn \le& \tn  \big|Y^{m,n}_t (\o)  \big|^p  \+ p\Big( |Y^{m,n}_t (\o) | \+  |\wt{y}_{m,n}|  \Big)^{p-1}
  \big|   \wt{y}_{m,n}  \big|   \ls    \big|Y^{m,n}_t (\o)  \big|^p  \+ \eta_n \qq \qq
  \eea
  with $\eta_n \df  \frac{p}{n}(2R \- 3)^{p-1}  $.  Also,   \eqref{eqn-d051} implies   that
 \bea
  && \hspace{-1.2cm} \big|\pi_m(Z^m_t (\o) \-  \hb{$ \frac{1}{m }$} \tilde{z})  \-  \pi_n(Z^n_t (\o)
  \- \hb{$\frac{1}{n}$}\tilde{z})\big|
 \ls  \big|\pi_m(Z^m_t (\o) \- \hb{$\frac{1}{m }$}\tilde{z}) \- \pi_m(Z^n_t (\o)
   \- \hb{$\frac{1}{n}$}\tilde{z})\big|
    \+ \big|\pi_m(Z^n_t (\o) \- \hb{$\frac{1}{n}$}\tilde{z}) \- \pi_n(Z^n_t (\o) \- \hb{$\frac{1}{n}$}\tilde{z})\big|  \nonumber \\
 & & \le   |Z^{m,n}_t (\o) \- \wt{z}_{m,n}|
 \+ \b1_{\{|Z^n_t (\o)   -   \hb{$\frac{1}{n}$}  \tilde{z}|>n\}}|Z^n_t (\o) \- \hb{$\frac{1}{n}$}\tilde{z}|
  \ls  |Z^{m,n}_t (\o)| \+  \hb{$\frac{2}{n }$}  \+ \b1_{\{|Z^n_t (\o)|>n-1\}}|Z^n_t (\o)|    .
    \qq \label{eqn-b509}
 \eea

 For any $u \ins L^2_\nu $, since $ \frac{k}{k \vee \| u \|_{L^p_\nu}} \= \frac{1}{1 \vee (\| u \|_{L^p_\nu}/k)}
 \nearrow 1$ as $k \n \to \n \infty$,  one can deduce that
 \beas
 \big\| \z_m (u)     \-   \z_n (u)     \big\|_{L^p_\nu}
 \= \b1_{\{\| u \|_{L^p_\nu} > n\}} \Big( \frac{m}{m \ve \| u \|_{L^p_\nu}}
 \- \frac{n}{n \ve \| u \|_{L^p_\nu}}  \Big) \| u \|_{L^p_\nu}
 \ls \b1_{\{\| u \|_{L^p_\nu} > n\}} \| u \|_{L^p_\nu} ,
 \eeas
 which together with the first inequality of \eqref{eq:b415} implies    that
 \bea
   \big\| \z_m \big( U^m_t (\o)\big)
     \-   \z_n \big( U^n_t (\o)\big)      \big\|_{L^p_\nu}
     & \tn  \ls & \tn  \big\| \z_m \big( U^m_t (\o)\big)
     \-   \z_m \big( U^n_t (\o)\big)      \big\|_{L^p_\nu}
     \+  \big\| \z_m \big( U^n_t (\o)\big)
     \-   \z_n \big( U^n_t (\o)\big)      \big\|_{L^p_\nu} \nonumber \\
       & \tn   \ls  & \tn   2 \big\| U^{m,n}_t (\o)  \big\|_{L^p_\nu}
       \+ \b1_{\{\| U^n_t (\o) \|_{L^p_\nu} > n\}} \big\| U^n_t (\o) \big\|_{L^p_\nu}
     \label{eq:a471} \\
      & \tn  \ls  & \tn  2 \big\| U^{m,n}_t (\o)  \big\|_{L^p_\nu} \+ n^{\frac{p-2}{p}} \big\| U^n_t (\o) \big\|^{\frac{2}{p}}_{L^p_\nu}
       \ls 2 \big\| U^{m,n}_t (\o)  \big\|_{L^p_\nu}
      \+ n^{\frac{p-2}{p}} C^{\frac{2}{p}}_{p,\cX} \big\| U^n_t (\o) \big\|^{\frac{2}{p}}_{L^2_\nu}     . \nonumber
 \eea
 Since $\big\| U^{m,n}_t (\o)  \big\|_{L^p_\nu} \ls C_{p,\cX} \big\| U^{m,n}_t (\o)  \big\|_{L^2_\nu}$,
 plugging  this inequality and \eqref{eqn-b505}$-$\eqref{eqn-b509}     into \eqref{eqn-b501},
  we can deduce from the monotonicity of function $\th $     that
 \bea
 && \hspace{-1.3cm} I^1_{t,\o}(\tilde{y},\tilde{z}) \ls
 \l(t)\th \big( |Y^{m,n}_t (\o)|^p \+ \eta_n \big)  \+  \Phi_t (\o)\big( |Y^{m,n}_t (\o)|^p \+ \eta_n \big)
    \+    \big(|Y^{m,n}_t (\o) |^{p-1}  \+   n^{1-p} \big) \ti \nonumber  \\
 &   & \hspace{-0.8cm} \q  \Big[ \Psi^n_t (\o) \+  \L_t (\o) |Z^{m,n}_t (\o)|
 \+  2 c_2 (t) \| U^{m,n}_t (\o) \|_{L^p_\nu} \Big]   \nonumber  \\
 & & \hspace{-0.8cm} \le    \l(t)\th \big( |Y^{m,n}_t (\o)|^p  \+  \eta_n \big)  \+  \Phi_t (\o)
  |Y^{m,n}_t (\o)|^p    \+ \eta_n \Phi_t (\o)  \+ \big[1 \+  (2R \- 4)^{p-1} \big]  \Psi^n_t (\o)
  \+    |Y^{m,n}_t (\o) |^{p-1}  \Big[ \L_t (\o) |Z^{m,n}_t (\o)|    \nonumber  \\
  &   & \hspace{-0.8cm} \q      +   2  c_2 (t) \| U^{m,n}_t (\o) \|_{L^p_\nu}  \Big]
       \+ \frac12 n^{1-p} \Big(    \L^2_t (\o)      \+     |Z^{m,n}_t (\o)|^2
    \+  4 C^2_{p,\cX} (c_2(t))^2 \+   \big\|U^{m,n}_t (\o)\big\|^2_{L^2_\nu} \Big)  ,
       \qq   \label{eqn-b477a}
 \eea
 where $\Psi^n_t (\o) \df  \L_t (\o) \big( \frac{2}{n }  \+ \b1_{\{|Z^n_t (\o)|>n-1\}}|Z^n_t (\o)| \big)
  \+   n^{\frac{p-2}{p}} c_2 (t) C^{\frac{2}{p}}_{p,\cX} \big\| U^n_t (\o) \big\|^{\frac{2}{p}}_{L^2_\nu}   \ls
     \frac{2}{n } \L_t (\o) \+ (n \- 1)^{\frac{-\fe}{2+\fe}}  \L_t (\o)
     |Z^n_t (\o)|^{\frac{2+2\fe}{2+\fe}}
   \dn + \n   n^{\frac{p-2}{p}} c_2 (t) C^{\frac{2}{p}}_{p,\cX} \big\| U^n_t (\o) \big\|^{\frac{2}{p}}_{L^2_\nu}  $.

   \no {\bf 3b)} As $\big\|\z_n(U^n_t (\o)) \big\|_{L^p_\nu} \ls \big\| U^n_t (\o)  \big\|_{L^p_\nu}
   \ls C_{p,\cX} \|U^n_t (\o)\|_{L^2_\nu} $,        (H2'), (H4) and (H5) show  that
   \bea
   | h^n_{t,\o}(\tilde{y},\tilde{z})| & \tn \le & \tn
    |f(t,\o,0,0,0)| \+   \k_0 \big( 1   \+   |Y^n_t (\o) \- \hb{$\frac{1}{n}$}\tilde{y} | \big)
     \+ \beta_t(\o) \+ c_1(t) \big|\pi_n(Z^n_t (\o) \- \hb{$\frac{1}{n}$}\tilde{z})\big| \+
   c_2 (t)  \big\|\z_n(U^n_t (\o)) \big\|_{L^p_\nu} \nonumber  \\
     & \tn \le & \tn     |f(t,\o,0,0,0)| \+ \k_0 R   \+ \beta_t(\o)  \+  c_1(t)
     \big(1 \+ |Z^n_t (\o)| \big)  \+   c_2 (t) C_{p,\cX}
       \|U^n_t (\o)\|_{L^2_\nu}   ,  \q     \label{eqn-b481}
   \eea
      which together with Lemma \ref{lem_a04} yields that
 \bea
   I^2_{t,\o}(\tilde{y},\tilde{z}) & \tn \le  & \tn
  \left| |Y^{m,n}_t (\o)|^{p-1}\sD\big(Y^{m,n}_t (\o)\big) \-  \big|\,Y^{m,n}_t (\o)  \- \wt{y}_{m,n}
 \big|^{p-1}  \sD \big( Y^{m,n}_t (\o) \- \wt{y}_{m,n} \big) \right| \Big(
\big|h^m_{t,\o}(\tilde{y},\tilde{z})\big| \+ \big|h^n_{t,\o}(\tilde{y},\tilde{z})\big|\Big) \nonumber \\
  & \tn   \le  & \tn   (1 \dn + \n 2^{p-1})n^{1-p} \Big[
 2  |f(t,\o,0,0,0)| \+ 2 \k_0 R  \+ 2 \beta_t(\o)
 \dn + \n  c_1 (t)\big(2 \+ |Z^m_t (\o)| \+ |Z^n_t (\o)| \big)
   \nonumber \\
  & \tn  & \tn    +    c_2 (t) C_{p,\cX}  \big( \|U^m_t (\o)\|_{L^2_\nu}
   \+   \|U^n_t (\o)\|_{L^2_\nu} \big)\Big] \df  \wt{I}^{\,2}_t (\o) .  \q \qq
 \label{eqn-b477b}
  \eea

     Since $0 \<  \beta^n_t (\o)  \ls  \beta^m_t (\o) \ls  1$, $\fa t  \ins    [0, T]$,
       \eqref{eqn-b481} also implies that
  \bea   \label{eqn-b477c}
  I^3_{t,\o}(\tilde{y},\tilde{z})  
  & \tn \ls  & \tn  (2R \- 4)^{p-1} ( 1  \- \beta^n_t (\o))  \Big[  |f(t,\o,0,0,0)| \+
 \k_0 R  \+ \beta_t(\o)  \+   c_1(t)
     \big(1 \+ |Z^n_t (\o)| \big) \nonumber \\
      & \tn  & \tn     +   c_2 (t) C_{p,\cX}   \|U^n_t (\o)\|_{L^2_\nu}   \Big]  \df  \wt{I}^{\,3}_t (\o) .
 \eea
  Putting \eqref{eqn-b477a}, \eqref{eqn-b477b} and \eqref{eqn-b477c} back  into \eqref{eqn-b471} shows that  \eqref{nonlip-cond3} is satisfied with $c(\cd) = 2 c_2 (\cd)$ and
     \beas
 \U^{m,n}_t  \=   \eta_n \Phi_t   \+   \big[1 \+  (2R \- 4)^{p-1} \big]    \Psi^n_t
   \+  \frac12 n^{1-p} \Big(    \L^2_t (\o)      \+     |Z^{m,n}_t (\o)|^2
    \+  4  C^2_{p,\cX} (c_2(t))^2 \+   \big\|U^{m,n}_t (\o)\big\|^2_{L^2_\nu} \Big)
     \+  \wt{I}^{\,2}_t  \+ \wt{I}^{\,3}_t  , ~     t    \ins   [0, T]   .
 \eeas

   H\"older's inequality, Young's inequality and \eqref{norm-esti} give rise to the following four estimates:
   \bea
 {\it a)}  && \hspace{-5mm} \underset{m>n}{\sup} \, E \n   \int_0^T \U^{m,n}_t dt
     \ls  \frac{p}{n}(2R-3)^{p-1}  C_\Phi
 + \big[1+  (2R-4)^{p-1} \big]   E\int_0^T   \Psi^n_t dt \nonumber \\
   &    &  \hspace{2.6cm}   +  \frac12 n^{1-p} \bigg(    C_\L
   \+  C^2_{p,\cX} \int_0^T (c_2 (t))^2 dt    \+   8 (R \- 2)^2  \bigg)
    \+  E \n \int_0^T  \n   \big( \wt{I}^{\,2}_t  \+  \wt{I}^{\,3}_t \big) dt  .
     \qq \hspace{1.7cm}   \label{eqn-a500} \\
  {\it b)} && \hspace{-5mm}   E\int_0^T  \Psi^n_t dt  \ls  \frac{2}{n} C^{(1)}_{\L} +
 (n-1)^{\frac{-\fe}{2+\fe}}  C^{(2)}_{\L}   \big\|Z^n\big\|_{\hZ^2}^{\frac{2+2\fe}{2+\fe}}
  +   n^{\frac{p-2}{p}} C^{\frac{2}{p}}_{p,\cX} \bigg(  \int_0^T \big( c_2 (t) \big)^q dt \bigg)^{\frac{1}{q}}
   \|U^n\|^{\frac{2}{p}}_{\hU^2}    \hspace{4.7 cm}  \nonumber \\
  & &  \hspace{1.3cm}  \le       \frac{2}{n} C^{(1)}_{\L} + 2
 (n \- 1)^{\frac{-\fe}{2+\fe}} C^{(2)}_{\L} (R \- 2)^{\frac{2+2\fe}{2+\fe}}
   \+   2^{\frac{1}{p}} n^{\frac{p-2}{p}}  C^{\frac{2}{p}}_{p,\cX}
    \ol{C}^{\frac{1}{q}}   (R \- 2)^{\frac{2}{p}}   , \nonumber
 \eea
 where $C^{(1)}_{\L} \df   E \n \int_0^T  \n  \L_t  dt $ and
 $C^{(2)}_{\L}  \df  \big( E \n \int_0^T  \n  \L_t^{2+\fe} dt \big)^{\frac{1}{2+\fe}}$.
 \beas
  {\it c)} \hspace{1mm}   E \n \int_0^T    \wt{I}^{\,2}_t  dt
 & \tn \dn \le& \tn  \dn   (1 \+ 2^{p-1})n^{1-p}
 \bigg\{ 2C_f \+ 2\k_0 R T  \+ 2 C_\beta    \+   \int_0^T \dn \Big( 2 c_1(t) \+  \frac12   (c_1(t))^2  \+
 \frac12 C^2_{p,\cX}   (c_2(t))^2 \Big) dt \\
 && \tn \+  \sum_{i=m,n}\Big(\big\|Z^i\big\|^2_{\hZ^2}  \+  \big\|U^i \big\|^2_{\hU^2}  \Big)\bigg\} \hspace{2cm}  \\
 &  \tn \dn \le&  \tn \dn    (1 \+ 2^{p-1})n^{1-p}
 \bigg\{ 2C_f \+ 2\k_0 R T   \+ 2 C_\beta     \+  \int_0^T \dn \Big( 2 c_1(t) \+  \frac12   (c_1(t))^2  \+
 \frac12 C^2_{p,\cX}   (c_2(t))^2 \Big) dt   \+  4  (R \- 2)^2\bigg\}  . \\
 {\it d)}  \hspace{1mm}  E \n \int_0^T    \wt{I}^{\,3}_t   dt
&  \tn   \dn \le&  \tn   \dn    (2R \- 4)^{p-1} \bigg\{  \n E  \dn \int_0^T  \n ( 1  \-  \beta^n_t) \big(
 |f(t,\o,0,0,0)| \+\k_0 R \+ \beta_t  \+ c_1 (t)\big) dt   \+  \|Z^n\|_{\hZ^2}
    \Big(    E \dn \int_0^T  \n  \big(c_1 (t)\big)^2( 1  \-  \beta^n_t)^2 dt \Big)^{\frac12} \nonumber  \\
& \tn   \dn & \tn   \dn  +   C_{p,\cX}  \|U^n\|_{\hU^2} \Big(    E \dn \int_0^T  \n  \big(c_2 (t)\big)^2( 1  \-  \beta^n_t)^2 dt \Big)^{\frac12}   \bigg\} \qq \qq \nonumber \\
&  \tn   \dn \le &  \tn   \dn    (2R \- 4)^{p-1} \bigg\{ \n E  \dn \int_0^T \n  ( 1  \-  \beta^n_t) \big(
 |f(t,\o,0,0,0)| \+ \k_0 R  \+ \beta_t \+ c_1 (t)\big) dt  \+ \sqrt{2} (R \- 2)\Big(
E  \dn  \int_0^T  \n   \big(c_1(t)\big)^2( 1  \- \beta^n_t)^2 dt \Big)^{\frac12}    \nonumber \\
& \tn   \dn & \tn   \dn  +   \sqrt{2} C_{p,\cX}  (R \- 2) \Big(    E \dn \int_0^T  \n  \big(c_2 (t)\big)^2( 1  \-  \beta^n_t)^2 dt \Big)^{\frac12}   \bigg\}  := J_n . \q
 \eeas

Because $  \beta^n_t \= \frac{ 1}{  1  \vee (\beta_t/n) \vee (|f(t,0,0,0)|/n) }  \n \nearrow \n  1$
as $n  \n \to \n  \infty$, $\fa t  \ins  [0,T]$, the dominated convergence theorem shows that
 $     \lmt{n \to \infty} J_n  \= 0$.
     Thus, letting $n  \n \to \n  \infty$ in \eqref{eqn-a500}
     yields that $ \lmt{n \to \infty} \, \underset{m>n}{\sup} \, E \n   \int_0^T \U^{m,n}_t dt \= 0  $.
 Moreover, since   $\|\cd\|_{\hD^p} \le \|\cd\|_{\hD^\infty} $, $\|\cd\|_{\hZ^{2,p}} \le \|\cd\|_{\hZ^2}  $ and  $\|\cd\|_{\hU^p} \le \big(\nu(\cX)T \big)^{\frac{2-p}{2p}} \|\cd\|_{\hU^2}  $ by H\"older's inequality, we see from \eqref{norm-esti} that
  \eqref{eqn-b210} also holds.      Then Proposition \ref{prop_stab} shows that $\big\{(Y^n, Z^n, U^n)\big\}_{n \in \hN} $
 is a Cauchy sequence in $ \hS^p$. Let
 $(Y, Z, U)$ be its limit.

  \no  {\bf 4)} {\it In this part, we will extract  an almost-surely convergent   and summable subsequence
   $ \big\{(Y^{m_i},Z^{m_i},U^{m_i}) \big\}_{i \in \hN} $   from $\{(Y^n, Z^n, U^n)\}_{n \in \hN}$. }

    Since
 \beas 
 \lmt{n \to \infty} E    \bigg[ \, \underset{t \in [0,T]}{\sup} \big|Y^n_t   \-  Y_t  \big|^p
+ \bigg( \n \int_0^T \n  \big|Z^n_t  \-  Z_t  \big|^2    \,  dt \bigg)^{\frac{p}{2}}
+  \n \int_0^T  \n
   \int_\cX   |U^n_t (x)  \-  U_t (x)    |^p \nu(dx) dt \bigg]=0 ,
    \eeas
    we can extract a subsequence $\{m_i\}_{i \in \hN}$ from $\hN$    such that
        \bea
     (i)&&   \hspace{-3mm}   \lmt{i \to \infty} \,   \underset{t \in [0,T]}{\sup} \big|Y^{m_i}_t  \-  Y_t  \big|
           \=  \lmt{i \to \infty}   \int_0^T \n  \big|Z^{m_i}_t  \-  Z_t  \big|^2    \,  dt
          \= 0 , \q \pas , \hspace{7.5cm}  \label{eqn-b356} \\
      (ii)&&  \hspace{-3mm}   \lmt{i \to \infty}   \big\|U^{m_i}_t  \-  U_t  \big\|_{L^p_\nu}  = 0 , \q \dtp , \label{eq:a411} \\
   (iii)&&  \hspace{-3mm} \big\|Y^{m_{i+1}} \- Y^{m_i} \big\|_{\hD^p}
      \ve  \big\|Z^{m_{i+1}} \- Z^{m_i} \big\|_{\hZ^{2,p}}
      \ve   \big\|U^{m_{i+1}} \- U^{m_i} \big\|_{\hU^p} \le 2^{-i}   , \q  \fa i \in \hN  . \label{eq:b515}
      \eea

     By \eqref{norm-esti}, it holds \pas~ that
 $  \underset{t \in [0, T]}{\sup}| Y_t|
 \ls  \underset{t \in [0,   T]}{\sup}\big|Y_t \- Y^{m_i}_t\big| \+  \underset{t \in [0, T]}{\sup}\big|Y^{m_i}_t\big|
    \ls  \underset{t \in [0,  T]}{\sup}\big|Y^{m_i}_t \- Y_t\big| \+ R \- 2 $,
    $ \fa i  \ins  \hN$.
 Letting $ i  \n \to \n  \infty$, we see from \eqref{eqn-b356} that
 \bea \label{eq:a317}
  \underset{t \in [0, T]}{\sup}| Y_t|  \ls R \- 2 , \q  \pas, \q \hb{thus} \q \|Y \|_{\hD^\infty}    \ls  R \- 2 .
 \eea

      For any $i \in \hN$,    we define two $[0,\infty)-$valued, $\bF -$predictable processes
    \beas
        \cZ^i_t :=   |Z_t|    +     \sum_{j=1}^i   \big|Z^{m_j}_t - Z^{m_{j-1}}_t  \big|
           \q \hb{and}   \q   \cU^i_t :=    \|  U_t  \|_{L^p_\nu}+  \sum_{j=1}^i  \big\|U^{m_j}_t - U^{m_{j-1}}_t \big\|_{L^p_\nu}   \,, \q    t \in [0,T]
        \eeas
            with $ Z^{m_0} :=  Z$ and $ U^{m_0} :=  U$.
        Minkowski's inequality and \eqref{eq:b515} imply  that
       \bea   \label{eqn-b360a}
       \bigg\{ E \bigg[ \Big( \n \int_0^T \n  ( \cZ^i_t )^2  \, dt \Big)^{\frac{p}{2}} \bigg]\bigg\}^{\frac{1}{p}}
         \ls    \big\|  Z \big\|_{\hZ^{2,p}}
           \+       \sum_{j=1}^i     \big\|Z^{m_{j}}  \-  Z^{m_{j-1}} \big\|_{\hZ^{2,p}}
         \ls    1+\big\|  Z \big\|_{\hZ^{2,p}}
          \+    \big\| Z^{m_1}  \-  Z \big\|_{\hZ^{2,p}}   , \q \fa i \in \hN         .
        \eea
   As $\big\{\cZ^i\big\}_{i \in \hN}$ is an increasing sequence,
   \bea
  \cZ_t   \df  \lmtu{i \to \infty} \cZ^i_t
   \=   |Z_t|    \+  \sum_{j=1}^\infty  \big|Z^{m_{j}}_t  \-  Z^{m_{j-1}}_t \big|    \,, \q
    t \in [0,T]     \label{eq:a425}
 \eea
  defines an $[0,\infty]-$valued, $\bF-$predictable process.
        The monotone convergence theorem shows that
    \beas
    \q   \int_0^T \n \big( \cZ_t (\o)\big)^2 dt
     \=  \lmtu{i \to \infty} \int_0^T  \n  \big( \cZ^i_t (\o)\big)^2 dt
        \q       \hb{and thus} \q
    \bigg(  \int_0^T  \n  \big( \cZ_t (\o)\big)^2 dt \bigg)^{\frac{p}{2}} \=  \lmtu{i \to \infty}
    \bigg(  \int_0^T  \n  \big( \cZ^i_t (\o)\big)^2 dt \bigg)^{\frac{p}{2}}   ,
    ~ \; \fa \o  \ins  \O .
    \eeas
     Applying the monotone convergence theorem once again,
    we can deduce from \eqref{eqn-b360a} and Lemma \ref{lem_lp_esti}  that
    \bea
            E \bigg[ \n  \Big(  \int_0^T  \cZ_t^2 dt \Big)^{\frac{p}{2}}    \bigg]
        \n =  \lmtu{i \to \infty} E \bigg[ \n \Big(  \int_0^T  \big( \cZ^i_t \big)^2 dt \Big)^{\frac{p}{2}}   \bigg]
        \le 3^{p-1}  \Big(1\n +\n \big\|   Z \big\|^p_{\hZ^{2,p}}
        \n +\n \big\| Z^{m_1} \n -\n  Z \big\|^p_{\hZ^{2,p}}   \Big) < \infty .       \label{eqn-b513}
    \eea

     Minkowski's inequality and \eqref{eq:b515} also imply  that
       \bea   \label{eqn-b360d}
        \bigg\{  E   \n \int_0^T \n   ( \cU^i_t )^p   dt \bigg\}^{\frac{1}{p}}
         \ls    \big\|  U \big\|_{\hU^p}
           \+       \sum_{j=1}^i     \big\|U^{m_{j}}  \-  U^{m_{j-1}} \big\|_{\hU^p}
         \ls    1+\big\|  U \big\|_{\hU^p}
          \+    \big\| U^{m_1}  \-  U \big\|_{\hU^p}   , \q \fa i \in \hN         .
        \eea
   As $\big\{\cU^i\big\}_{i \in \hN}$ is an increasing sequence,
 \bea
 \cU_t \df  \lmtu{i \to \infty} \cU^i_t
     \=   \|U_t\|_{L^p_\nu}     \+  \sum_{j=1}^\infty
    \big\|U^{m_j}_t  \-  U^{m_{j-1}}_t \big\|_{L^p_\nu}    \,, \q t  \ins  [0,T] \label{eq:a426}
 \eea
    defines an $[0,\infty]-$valued, $\bF-$predictable process.
    Applying the monotone convergence theorem again, we can deduce from \eqref{eqn-b360d}
    and Lemma \ref{lem_lp_esti}  that
    \bea  \label{eqn-b513d}
     E   \n \int_0^T \n  \cU^p_t  dt   = \lmtu{i \to \infty} E   \n \int_0^T \n   ( \cU^i_t )^p   dt
     \le 3^{p-1} \Big(  1+\big\|  U \big\|^p_{\hU^p}
          \+    \big\| U^{m_1}  \-  U \big\|^p_{\hU^p}  \Big) < \infty .
    \eea

 \no {\bf 5)}  {\it Finally, we will send $i \n \to \n \infty$ in BSDEJ\,$(\xi,  f_{m_i})$
 to demonstrate that the processes  $(Y,Z,U)$ solve  BSDEJ\,$(\xi,  f)$.}

 \ss    Fix $ k \ins  \hN$.  We  define an $\bF -$stopping time
  \bea \label{def_tau_k}
  \t_k  \df  \inf  \bigg\{ t  \ins  [0,T] \n :  \n
  \int_0^t  \n \cZ^2_s    ds   \>  k  \bigg\}    \land   T   .
 \eea
      Since    $   \int_0^{\t_k} \n     | Z^{m_i}_t \- Z_t  |^2 dt
         \ls  \int_0^{\t_k}  \n       (\cZ^i_t  )^2  dt
        \ls   \int_0^{\t_k}  \n   \cZ_t^2 dt            \ls  k     $,
        $  \fa \o  \ins  \O$,
        the dominated convergence theorem and \eqref{eqn-b356} show that
    \bea    \label{eqn-b362}
         \lmt{i \to \infty}     E \n  \int_0^{\t_k} \n  \big|Z^{m_i}_t \-  Z_t  \big|^2    \,  dt = 0 .
         \eea
       Hence,  there exists a subsequence $\big\{m^k_i \big\}_{i \in \hN}$ of $\{m_i\}_{i \in \hN}$
   such that for \dtp~ $(t,\o) \in [0,T] \times \O $
    \bea
         \lmt{i \to \infty}    \b1_{\{t \le  \t_k\}}    \big|Z^{m^k_i}_t  \-  Z_t  \big|
          = 0  .  \label{eqn-b521}
         \eea

 We shall  show that
  \bea
   \lmt{i \to \infty}  E \n \int_0^{  \t_k} \left|   f_{m^k_i}   \Big(t, Y^{m^k_i}_t, Z^{m^k_i}_t, U^{m^k_i}_t\Big) -f \big(t, Y_t, Z_t,
U_t\big) \right| dt   =0 .   \label{eqn-b519}
  \eea
 Since $ \psi (x) \n \equiv \n 1$ for all $|x| \ls R \- 1$,    \eqref{norm-esti} implies that for any $i \ins  \hN$
  \bea
   \q  && \hspace{-1.5cm}  E \n \int_0^{  \t_k} \left|   f_{m^k_i}   \Big(t, Y^{m^k_i}_t, Z^{m^k_i}_t, U^{m^k_i}_t\Big) -f \big(t, Y_t, Z_t,
U_t\big) \right| dt     \nonumber   \\
& \= & \tn \n       E \n  \int_0^{ \t_k} \dn \int_{|(\tilde{y},\tilde{z})|<1} \n  \left|\beta^{m^k_i}_t f\bigg(t,Y^{m^k_i}_t\n -\n \hb{$\frac{1}{{m^k_i}}$}\tilde{y},
\pi_{m^k_i}\Big(Z^{m^k_i}_t\n -\n\hb{$\frac{1}{{m^k_i}}$}\tilde{z}\Big),\z_{m^k_i}\Big(U^{m^k_i}_t\Big) \n \bigg)
\n-\n f \big(t, Y_t, Z_t, U_t\big) \right| \rho(\tilde{y},\tilde{z})d\tilde{y}d\tilde{z} dt       .   \qq  \label{eqn-b517}
   \eea

  For any $(t,\o) \ins [0,T]  \ti  \O$
  except on a $dt  \oti  dP-$null  set $\fN_k \n \supset \n \fN $,  we may assume further that
  \eqref{eq:a411}, \eqref{eqn-b521} hold, that
  $\lmt{i \to \infty}   \big| Y^{m^k_i}_t (\o)     \-  Y_t (\o) \big| \= 0$ (by  \eqref{eqn-b356}),
    that   $|Y_t (\o)|  \ls  R \- 2$ (by \eqref{eq:a317}),
    and that $U_t (\o) \ins  L^p_\nu $.

   Let $(t, \o)  \ins  \fN^c_k \n \cap \n \[0,  \t_k \]$
and let $(\tilde{y},\tilde{z})  \ins  \hR^l  \ti  \hR^{l \times d} $
with $|(\tilde{y},\tilde{z})| \< 1$.
Since
  \bea \label{eq:a347}
  \lmt{i \to \infty}         \big|Z^{m^k_i}_t (\o) \-  Z_t (\o)  \big|  \=  0
  \eea
   by \eqref{eqn-b521},  Lemma \ref{lem_pi} and the first inequality of \eqref{eq:b415} imply that

  \no (e1)   $\Big|Y^{m^k_i}_t (\o) \- \frac{1}{m^k_i}\tilde{y} \- Y_t(\o)\Big|
 \ls  \frac{1}{m^k_i} \+ \Big|Y^{m^k_i}_t (\o) \- Y_t(\o)\Big|   \n \to \n  0$, as $i  \n \to \n  \infty$;

 \no (e2)  $\Big|\pi_{m^k_i}\Big(Z^{m^k_i}_t (\o)  \- \frac{1}{m^k_i}\tilde{z}\Big) \- Z_t (\o) \Big| \ls
 \Big|\pi_{m^k_i}\Big(Z^{m^k_i}_t  (\o) \- \frac{1}{m^k_i}\tilde{z}\Big) \- \pi_{m^k_i}\big(Z_t (\o) \big)\Big|
    \+   \big|\pi_{m^k_i}\big(Z_t (\o) \big) \- Z_t (\o) \big|
  \ls  \Big|Z^{m^k_i}_t(\o) \- \frac{1}{m^k_i}\tilde{z} \- Z_t(\o)\Big|
     \+   \big|\pi_{m^k_i}\big(Z_t(\o)\big) \- Z_t(\o)\big|
      \ls  \frac{1}{m^k_i}+\Big|Z^{m^k_i}_t(\o) \- Z_t(\o)\Big|
  \+ \big|\pi_{m^k_i}\big(Z_t(\o)\big) \- Z_t(\o)\big| \n \to \n  0$, as $i  \n \to \n  \infty$;

 \no (e3) $\Big\|\z_{m^k_i}\Big(U^{m^k_i}_t(\o)\Big) \- U_t(\o)\Big\|_{L^p_\nu}  \ls
\big\|\z_{m^k_i}\Big(U^{m^k_i}_t(\o)\Big) \- \z_{m^k_i}(U_t(\o))\big\|_{L^p_\nu}
 \+ \big\|\z_{m^k_i}(U_t(\o))    \-    U_t(\o)\big\|_{L^p_\nu}
  \ls 2 \big\|U^{m^k_i}_t(\o)  \-  U_t(\o)\big\|_{L^p_\nu}
  \+ \Big( 1 \- \frac{m^k_i}{m^k_i \vee \|U_t(\o)\|_{L^p_\nu} } \Big) \big\|  U_t(\o) \big\|_{L^p_\nu}    \to 0$, as $i \to
 \infty$.

    Since the mapping $     f \big(t, \o, \cd,\cd,  U_t (\o) \big)$ is  continuous by (H1) and since
$\lmtu{i \to \infty}  \beta^{m^k_i}_t (\o) \=  1$, we can deduce  from (e1) and (e2)  that
 \bea   \label{eqn-b525}
 \lmt{i \to \infty} \beta^{m^k_i}_t (\o) f\Big(t, \o, Y^{m^k_i}_t (\o) \n -\n \hb{$\frac{1}{{m^k_i}}$}\tilde{y},
\pi_{m^k_i}\Big(Z^{m^k_i}_t (\o) \n -\n\hb{$\frac{1}{{m^k_i}}$}\tilde{z}\Big),U_t (\o) \Big)
 =  f \big(t, \o, Y_t (\o), Z_t(\o), U_t(\o)\big) .
 \eea
 Moreover, (H5) shows that
  \beas
 \q \qq && \hspace{-2cm}  \bigg| 
 f\Big(\n t,\o, Y^{m^k_i}_t (\o) \- \hb{$\frac{1}{{m^k_i}}$}\tilde{y},
    \pi_{m^k_i}\Big(\n Z^{m^k_i}_t (\o)\n-\n\hb{$\frac{1}{{m^k_i}}$}\tilde{z}\Big),\z_{m^k_i}\Big(\n U^{m^k_i}_t(\o)\Big) \n \Big)   \- 
f\Big(t,Y^{m^k_i}_t (\o)\n-\n \hb{$\frac{1}{{m^k_i}}$}\tilde{y},
\pi_{m^k_i}\Big(\n Z^{m^k_i}_t (\o)\n-\n \hb{$\frac{1}{{m^k_i}}$}\tilde{z}\Big),U_t (\o) \Big) \bigg| \\
&&   \le   c_2(t)    \Big\|\z_{m^k_i}\Big(U^{m^k_i}_t (\o) \Big) \- U_t (\o) \Big\|_{L^p_\nu} \, ,
 \eeas
  which together with \eqref{eqn-b525} and (e3) implies that
  \bea  \label{eqn-b523}
 \lmt{i \to \infty}  \Big| \beta^{m^k_i}_t (\o)  f \Big(t, \o, Y^{m^k_i}_t (\o) \n -\n \hb{$\frac{1}{{m^k_i}}$}\tilde{y},
 \pi_{m^k_i}\Big(Z^{m^k_i}_t (\o)
    \n -\n\hb{$\frac{1}{{m^k_i}}$}\tilde{z}\Big),\z_{m^k_i}\Big(U^{m^k_i}_t(\o)\Big)   \Big)
   \-   f \big(t, \o,  Y_t(\o), Z_t(\o), U_t(\o)\big) \Big| \= 0 . \q
 \eea

 Given $ i \ins  \hN $, there exists an $\wh{j} = \wh{j} (k,i) \ins \hN$ such that
 $m^k_i \= m_{\wh{j}}$.  Since
 \bea \label{eq:a477}
   \Big|Z^{m^k_i}_t (\o)\Big|  \ls  \cZ^{\wh{j}}_t (\o)  \ls  \cZ_t (\o)  \q \hb{and} \q
   \Big\|  \z_{m^k_i}\Big(U^{m^k_i}_t (\o) \Big)  \Big\|_{L^p_\nu}
  \ls \Big\|U^{m^k_i}_t (\o) \Big\|_{L^p_\nu}  \ls  \cU^{\wh{j}}_t (\o)  \ls  \cU_t (\o)  ,
 \eea
  one can deduce from (H2'), (H4) and (H5)  that
 \beas
 \qq   && \hspace{-1.5cm} \Big| \beta^{m^k_i}_t (\o)  f \Big(t, \o, Y^{m^k_i}_t (\o) \n -\n \hb{$\frac{1}{{m^k_i}}$}\tilde{y},
 \pi_{m^k_i}\Big(Z^{m^k_i}_t (\o)
    \n -\n\hb{$\frac{1}{{m^k_i}}$}\tilde{z}\Big),\z_{m^k_i}\Big(U^{m^k_i}_t(\o)\Big)   \Big)
   \-   f \big(t, \o,  Y_t(\o), Z_t(\o), U_t(\o)\big) \Big|  \\
&  \dn   \dn \le&  \dn  \dn  \bigg|  f \Big(t, \o, Y^{m^k_i}_t (\o) \n -\n \hb{$\frac{1}{{m^k_i}}$}\tilde{y},
 \pi_{m^k_i}\Big(Z^{m^k_i}_t (\o)
    \n -\n\hb{$\frac{1}{{m^k_i}}$}\tilde{z}\Big),\z_{m^k_i}\Big(U^{m^k_i}_t(\o)\Big)   \Big) \bigg|
 + \Big| f \big(t, \o,  Y_t(\o), Z_t(\o), U_t(\o)\big) \Big|  \\
 & \dn  \dn \le& \dn  \dn 2 |f(t,\o,0,0,0)| \+
 \k_0 \Big( 2 \+ \big|Y^{m^k_i}_t (\o) \n -\n \hb{$\frac{1}{{m^k_i}}$}\tilde{y}\big|
 \+ |Y_t (\o) | \Big)
  \+ 2 \beta_t(\o) \+  c_1 (t) \Big(    \big| Z^{m^k_i}_t  (\o)
    \n -\n\hb{$\frac{1}{{m^k_i}}$}\tilde{z}  \big| \+ |Z_t (\o) | \Big) \\
 & \dn  \dn &  \dn  \dn   +   c_2 (t) \Big(   \Big\|  \z_{m^k_i}\big(U^{m^k_i}_t  (\o)\big)  \Big\|_{L^p_\nu}
 \+ \|U_t (\o)\|_{L^p_\nu} \Big) \\
 & \dn  \dn \le& \dn  \dn  2 |f(t,\o,0,0,0)| \+ (2R \- 1) \k_0 \+ 2 \beta_t(\o)
\+ c_1(t)     \big(   1 \+  \cZ_t  (\o) \+ |Z_t (\o)|  \big)
 \+ c_2(t)    \big( \cU_t (\o) \+ \|U_t (\o)\|_{L^p_\nu}    \big)  := H_t (\o).
 \eeas
   Applying Holder's inequality, we see from \eqref{eqn-b513} and \eqref{eqn-b513d}  that
   \bea
 \qq  && \hspace{-1.2cm}   E \int_0^{\t_k} \n \int_{|(\tilde{y},\tilde{z})|<1} H_t \,
\rho(\tilde{y},\tilde{z})d\tilde{y}d\tilde{z}    dt  = E \int_0^{\t_k}   H_t dt \le  E \int_0^T   H_t dt   \nonumber  \\
 & \tn \le & \tn  C \+ E  \n \left[ \Big(\int_0^T \n \big(c_1(t)\big)^2 dt \Big)^{\frac12}
   \bigg\{ \Big(\int_0^T \n \cZ^2_t dt
 \Big)^{\frac12}    \+  \Big(\int_0^T  \n  |Z_t|^2 dt \Big)^{\frac12}    \bigg\} \right]
    \+  \bigg(\int_0^T \n  (c_2(t))^q dt\bigg)^{\frac{1}{q}}   \bigg\{
 \Big( E \int_0^T  \n  \cU^p_t dt       \Big)^{\frac{1}{p}}
 \+  \|U\|_{\hU^p} \bigg\}   \nonumber   \\
 & \tn \le & \tn   C   \+ \ol{C}^{\frac12}  \Bigg\{ \bigg(E \bigg[ \Big(\int_0^T \n \cZ^2_t dt
 \Big)^{\frac{p}{2}} \bigg]\bigg)^{\frac{1}{p}} + \|Z\|_{\hZ^{2,p}} \Bigg\}
 \+ \ol{C}^{\frac{1}{q}}   \bigg\{
 \Big( E \int_0^T  \n  \cU^p_t dt       \Big)^{\frac{1}{p}} \+  \|U\|_{\hU^p} \bigg\}  < \infty  \label{eqn-b539}
\eea
with $C := 2 C_f \+ (2R \-  1) \k_0  T  \+ 2 C_\beta \+  \int_0^T \n c_1(s) ds \< \infty $.
Hence,  the dominated convergence theorem and \eqref{eqn-b523} show that
 \beas
  \lmt{i \to \infty}  E \n  \int_0^{ \t_k} \dn \int_{|(\tilde{y},\tilde{z})|<1} \n  \left|\beta^{m^k_i}_t f\bigg(t,Y^{m^k_i}_t\n -\n \hb{$\frac{1}{{m^k_i}}$}\tilde{y},
\pi_{m^k_i}\Big(Z^{m^k_i}_t\n -\n\hb{$\frac{1}{{m^k_i}}$}\tilde{z}\Big),\z_{m^k_i}\Big(U^{m^k_i}_t\Big) \n \bigg)
\n-\n f \big(t, Y_t, Z_t, U_t\big) \right| \rho(\tilde{y},\tilde{z})d\tilde{y}d\tilde{z} dt  =0,
 \eeas
 which together with \eqref{eqn-b517} leads to \eqref{eqn-b519}.

   Since  $\int_{( \tau_k \land t, \t_k]} =\int_{(0,  \t_k]} -\int_{(0,\tau_k \land t]}  $, $\fa   t \in [0,T]$, the Burkholder-Davis-Gundy inequality, H\"older's inequality
  and an analogy to \eqref{eq:a017} imply that
    \bea
  && \hspace{-0.8cm}   E\left[   \underset{t \in [0,T]}{\sup} \bigg|  \int_{( \tau_k \land t, \t_k]} \int_\cX \big(U^{m^k_i}_s(x) \- U_s(x)\big)\tnp(ds,dx) \bigg|  \right]
  \ls  2    E\left[   \underset{t \in [0,T]}{\sup} \bigg|  \int_{(0,\tau_k \land t]} \int_\cX \big(U^{m^k_i}_s(x)-U_s(x)\big)\tnp(ds,dx) \bigg|  \right] \nonumber  \\
    && \le   c_l  \,   E   \bigg[      \Big(  \int_{(0,\t_k]} \int_\cX \big| (U^{m^k_i}_s(x) \-  U_s(x)) \big|^2  N_\fp(ds,dx) \Big)^{\frac{1}{2}} \bigg]
   \ls   c_l  \bigg\{ E     \bigg[   \Big(  \int_{(0,T]} \int_\cX \big| (U^{m^k_i}_s(x) \- U_s(x)) \big|^2  N_\fp(ds,dx) \Big)^{\frac{p}{2}} \bigg]   \bigg\}^{\frac{1}{p}} \nonumber  \\
   &&  \le     c_l  \bigg\{ E \n  \int_{(0,T]} \n \int_\cX \big| (U^{m^k_i}_s(x) \- U_s(x)) \big|^p  N_\fp(ds,dx)    \bigg\}^{\frac{1}{p}}
    \=  c_l  \bigg\{ E  \n    \int_0^T \n \int_\cX \big| (U^{m^k_i}_s(x) \- U_s(x)) \big|^p \nu (dx) ds   \bigg\}^{\frac{1}{p}} \nonumber \\
   && = c_l \|U^{m^k_i} \- U \|_{\hU^p}
   \to 0, \hb{\, as }  i \to \infty ,    \label{eqn-b364}
   \eea
  and   that
   \bea
  E\left[   \underset{t \in [0,T]}{\sup} \bigg|  \int_{\tau_k \land t}^{\t_k}
    \big(Z^{m^k_i}_s \- Z_s\big) dB_s \bigg|       \right]
        \ls  
          c_l   E\left[   \bigg(   \int_0^{ \t_k}
          \big| Z^{m^k_i}_s \- Z_s \big|^2  ds  \bigg)^{\frac{1}{2}}       \right]
       \ls         c_l  \,   \big\| Z^{m^k_i} \- Z \big\|_{\hZ^{2,p}} \to 0, \hb{\, as } i \to \infty  .       \label{eqn-b366}
   \eea

     In light of \eqref{eqn-b356}, \eqref{eqn-b519},  \eqref{eqn-b364} and  \eqref{eqn-b366}, there exists a subsequence $\big\{\wt{m}^k_i \big\}_{i \in \hN}$ of $\big\{m^k_i \big\}_{i \in \hN}$   such that
    except on a $P-$null set $\cN^k_1$
    \beas
   \q && \hspace{-2cm}  \lmt{i \to \infty}   \bigg\{ \underset{t \in [0,T]}{\sup} \Big|Y^{\wt{m}^k_i}_t - Y_t  \Big|
   +     \int_0^{  \t_k} \left|   f_{\wt{m}^k_i}   \Big(t, Y^{\wt{m}^k_i}_t, Z^{\wt{m}^k_i}_t, U^{\wt{m}^k_i}_t\Big)
   - f \big(t, Y_t, Z_t,  U_t\big) \right| dt   \nonumber  \\
  && +    \underset{t \in [0,T]}{\sup} \bigg|  \int_{\tau_k \land t}^{\t_k} \dn
  \Big(Z^{\wt{m}^k_i}_s \- Z_s\Big) dB_s \bigg|
    +   \underset{t \in [0,T]}{\sup}  \bigg|  \int_{( \tau_k \land t, \t_k]}  \n \int_\cX \n \Big(U^{\wt{m}^k_i}_s \n (x)
   \- U_s(x)\Big)\tnp(ds,dx) \bigg|   \bigg\}   \n  =0.   
    \eeas
     Since  $ \big(Y^{\wt{m}^k_i}, Z^{\wt{m}^k_i}, U^{\wt{m}^k_i}\big) $
     solves   BSDEJ\,$\big( \xi, f_{\wt{m}^k_i} \big)$ for any $ i \ins \hN$,
     it holds  except on a $P-$null set $\cN^k_2$ that
     \beas
 Y^{\wt{m}^k_i}_{\tau_k \land t} & \tn =  & \tn    \b1_{\{\tau_k <T\}} Y^{\wt{m}^k_i}_{  \t_k}  \n +\n \b1_{\{\tau_k = T\}} \xi
 \n + \dn \int_{\tau_k \land t}^{  \t_k} \n f_{\wt{m}^k_i}\Big(s, Y^{\wt{m}^k_i}_s, Z^{\wt{m}^k_i}_s, U^{\wt{m}^k_i}_s\Big) \, ds \n - \dn \int_{\tau_k \land t}^{  \t_k} \n   Z^{\wt{m}^k_i}_s \, dB_s \nonumber \\
 & \tn & \tn  -   \int_{(\tau_k \land t, \t_k]} \n  \int_\cX \n  U^{\wt{m}^k_i}_s(x)\tnp(ds,dx),
  \q \fa t  \ins  [0,T] , \q \fa i \ins \hN . 
   \eeas
 Letting $i  \n \to \n  \infty$, we obtain that  over
 $\O_k  \df   \big( \cN^k_1 \big)^c   \n \cap \n   \big( \cN^k_2 \big)^c  $
\bea
 \hspace{-3mm}
Y_{\tau_k \land t }   \=   \b1_{\{\tau_k <T\}} Y_{  \t_k}  \n +\n \b1_{\{\tau_k = T\}} \xi
  \n  +   \dn  \int_{\tau_k \land t }^{  \t_k } f \big(s, Y_s, Z_s, U_s\big)  ds
  \n - \dn    \int_{\tau_k \land t}^{  \t_k}  Z_s \, dB_s
  \n  -  \dn  \int_{(\tau_k \land t, \t_k]} \int_\cX U_s(x)\tnp(ds,dx) , ~  t \ins  [0, T]. \qq \label{eqn-b376}
   \eea

  By \eqref{eqn-b513}, it holds   for all $\o \ins \O$ except on $P-$null set $\cN_Z$ that
  $    \int_0^T  \cZ_t^2 (\o) dt      \<   \infty $, and thus $ \t_{\fk}(\o) \=  T $
  for some $ \fk \= \fk(\o) \ins  \hN$.
  Then letting $k  \n \to \n \infty$ in \eqref{eqn-b376} shows that \eqref{BSDEJ} holds  over
  $  \Big(  \underset{k \in \hN}{\cap}  \O_k  \Big) \n \cap   \cN^c_Z $,
  which together with Remark \ref{rem_drift} shows that   $(Y, Z, U) $ is a solution of BSDEJ\,$(\xi,  f)$.   \qed

      \no {\bf Proof of Theorem \ref{thm_BSDEJ1}:}  {\bf (Uniqueness)}
     Suppose that $(Y,Z,U) , (Y',Z',U') \ins   \hS^p$ are two solutions of the BSDEJ\,$(\xi,f)$.
     For any $n  \ins  \hN$, we set
  \beas
 (\xi_n,f_n) := (\xi,f) \q  \hb{ and } \q  (Y^n,Z^n,U^n) := \left\{\ba{ll}  (Y,Z,U)  & \hb{if $n$ is odd,}   \\  (Y',Z',U') \q & \hb{if $n$ is even.} \ea\right.
 \eeas
   By an analogy to \eqref{eq:b411},
   the inequality \eqref{nonlip-cond3} holds for  $\eta_n \= 0$, $c(\cd) \= c_2 (\cd)$
    and $\U^{m,n} \n \equiv \n 0$.
      Proposition \ref{prop_stab} then shows that
   $\big\{(Y^n, Z^n, U^n)\big\}_{n \in \hN} $ is a Cauchy sequence in $  \hS^p   $,    which implies that
   $ \| Y \- Y'\|_{\hD^p}   \=    \|Z \- Z'\|_{\hZ^{2,p}}  \=  \| U \- U'\|_{\hU^p} \= 0 $.
   Hence, one has  that   $P\{ Y_t \=  Y'_t , \, \fa t  \ins  [0,T] \} \= 1$,
  that   $Z_t (\o) \= Z'_t  (\o)$ for \dtp ~ $(t,\o )  \ins  [0,T]  \ti  \O $,
  and that $U(t,\o,x) \= U'(t,\o,x)$ for $dt \ti dP  \ti  \nu (dx)-$a.s.
  $(t,\o,x)  \ins  [0,T]  \ti  \O  \ti  \cX$.

   \no {\bf (Existence)}

   \ss \no {\bf 1)}
 Let us first assume that $ \xi \ins L^\infty (\cF_T) $ and $  \int_0^T \n  |f(t,0,0,0)| dt   \ins  L^\infty_+(\cF_T) $.
 We   set
    \bea \label{eq:a335b}
   R  \df  2 \+ \exp \Big\{ T \+  C_f  \+ 4 C_\beta \+ 2 \int_0^T \n \big(c_1(t)\big)^2 dt
    \n + \n  4 \big(\nu(\cX)\big)^{\frac{2-p}{p}}   \n \int_0^T \n  \big(c_2(t)\big)^2 dt   \Big\}
   \ti \sqrt{\|\xi\|^2_{L^\infty(\cF_T)} \+ 5T    \n+\n   C_f / 2 \+  7 C_\beta / 2   }   , \q
   \eea
   and let
   $\p \n : \hR^l  \n \to \n  [0,1]$ be a smooth function such  that
   $\p (x) \= 1$ (resp.\;\,$\p (x) \= 0$) if   $|x| \ls  R \- 1 $  \big(resp.\;\,$|x|  \gs  R$\big).

    Let $n \ins  \hN$. For any $u \ins L^p_\nu$, we define
   $   \pi_n(u)   \df  \big( \frac{n}{n \vee \| u \|_{L^p_\nu}} \big) \, u  \ins L^p_\nu  $.
   An application of Lemma \ref{lem_pi2} with $(\hE,\|\cd\|) \= \big( L^p_\nu,\|\cd\|_{L^p_\nu}\big)$ shows that
   $          \big\|\pi_n (u_1)   \-   \pi_n(u_2)\big\|_{L^p_\nu} \ls  2 \|u_1    \-    u_2   \|_{L^p_\nu}   $,
        $  \fa u_1 , u_2 \ins   L^p_\nu  $,
   which together with \eqref{eqn-d051} and the
$\sP  \oti  \sB \big( \hR^l  \big)    \otimes \n  \sB \big( \hR^{l \times d} \big)
  \oti  \sB\big(L^p_\nu   \big)/\sB(\hR^l ) -$measurability of $f$ shows that
 \beas
   f_n (t,\o,y,z,u)  \df \frac{n}{n \ve  \phi^R  ( t , \o ) }\psi (y)
 \big( f (t,\o,y,\pi_n(z),\pi_n (u))  \-   f (t,\o,0,0,0) \big)  \+    f (t,\o,0,0,0)   ,
 \eeas
 $(t,\o,y,z,u)  \ins   [0, T]  \ti  \O  \ti  \hR^l  \ti  \hR^{l \times d}   \ti   L^p_\nu$
 defines a $   \hR^l-$valued,   $\sP  \oti  \sB  ( \hR^l   )   \oti  \sB \big( \hR^{l\times d} \big)
  \oti  \sB\big(L^p_\nu   \big)/\sB(\hR^l ) -$measurable mapping   satisfying (H1),   (H3)$-$(H5)
  with the same   coefficients as $f$ except for $c^n_2(\cd) \= 2 c_2 (\cd)$.
  By (H2),   it holds for \dtp ~ $(t,\o) \in [0,T]  \ti  \O$ that
    \beas
     |f_n(t, \o, y,0,0) \- f_n(t, \o, 0,0,0)|
     \= \frac{n}{n \ve  \phi^R_t }\psi (y) \big| f (t, \o, y,0,0) \- f (t, \o, 0,0,0) \big|
     \ls \frac{n}{n \ve  \phi^R_t }\psi (y) \phi^R_t    \ls  n  ,  \q \fa  y  \ins   \hR^l  ,
     \eeas
 so $f_n$ satisfies (H2') with  $\k_0 \= n $.

 Also, let  (H2) and (H4)$-$(H6) hold for $f$ except on a $dt \ti dP-$null set $ \fN $ and let $(t,\o) \ins \fN^c$.
 Given $(y_1,z_1 ) , (y_2,z_2 ) \ins  \hR^l  \ti  \hR^{l \times d} $ and $   u   \ins   L^p_\nu $,
 if $|y_1| \ld |y_2| \gs R $, then we automatically have
 $f_n (t,\o,y_1,z_1,u   ) \- f_n (t,\o,y_2,z_2,u   )      \=  0$
 and thus $ |y_1 \- y_2|^{p-1} \big\lan \sD ( y_1 \- y_2 ) ,  f_n (t,\o,y_1,z_1,u   )
 \- f_n (t,\o,y_2,z_2,u   )   \big\ran \= 0 $; on the other hand, let us assume without loss of generality that
 $|y_1| \< R $,  then (H2), (H4)$-$(H6) and \eqref{eqn-d051}   imply that
 \beas
 && \hspace{-1.5cm} |y_1 \- y_2|^{p-1} \big\lan \sD ( y_1 \- y_2 ) ,  f_n (t,\o,y_1,z_1,u   )
 \- f_n (t,\o,y_2,z_2,u   )   \big\ran \\
 && \= \frac{n}{n \ve  \phi^R_t (\o) } \big( \psi (y_1) \- \psi (y_2) \big)
  |y_1 \- y_2|^{p-1} \big\lan \sD ( y_1 \- y_2 ) ,  f  (t,\o,y_1, \pi_n(z_1), \pi_n(u)   )
 \- f  (t,\o,0,0,0   )   \big\ran \\
&& \q   +   \frac{n}{n \ve  \phi^R_t (\o) }\psi (y_2) |y_1 \- y_2|^{p-1} \big\lan \sD ( y_1 \- y_2 ) ,
 f  \big( t,\o,y_1,\pi_n(z_1),\pi_n(u)  )  \- f  (t,\o,y_2,\pi_n(z_2),\pi_n(u)  \big)   \big\ran \\
&& \ls \frac{n}{n \ve  \phi^R_t (\o) } C_\psi |y_1 \- y_2|^p \big( \phi^R_t (\o) \+ \beta_t (\o) \+ c_1(t) |\pi_n(z_1)|
\+ c_2(t) \|\pi_n(u)\|_{L^p_\nu} \big)  \\
&& \q    + \frac{n}{n \ve  \phi^R_t (\o) }\psi (y_2) \Big[ \l(t) \, \th \big( |y_1 \- y_2|^p \big)
\+  \Phi_t (\o) |y_1 \- y_2|^p
  \+    \L_t (\o)  |y_1 \- y_2|^{p-1}   \big|\pi_n(z_1) \- \pi_n(z_2)\big| \Big]         \\
    && \ls  \l(t) \, \th \big( |y_1 \- y_2|^p \big) \+
    \big[ \Phi_t (\o) \+ C_\psi \beta_t   (\o) \+ n C_\psi (1\+c_1(t)\+c_2(t) ) \big]    |y_1 \- y_2|^p
      +    \L_t (\o)  |y_1 \- y_2|^{p-1}   |z_1 \- z_2|   ,
 \eeas
 where $C_\psi$ denotes the Lipschitz coefficient of the smooth function $\psi$.
 Hence, $f_n$ satisfies (H6)   with  the same   coefficients as $f$ except for
 with $ \Phi^n_t \= \Phi_t   \+ C_\psi \beta_t     \+ n C_\psi (1\+c_1(t)\+c_2(t) )  $, $t \ins [0,T]$.
 Clearly, $\int_0^T \n \Phi^n_t dt \ins L^\infty_+ (\cF_T) $.

 Since $f_n$ satisfies    (H3)$-$(H5)    with the same   coefficients as $f$
 except for $c^n_2(\cd) \= 2 c_2 (\cd)$ and
 since $\int_0^T \n \big|f_n (t,0,0,0)\big|dt  \=  \int_0^T \n \big|  f (t, 0,0,0)  )\big|dt  \ins  L^\infty_+(\cF_T)  $,
   the constant $R$ in \eqref{eq:a335b}  is exactly that   for $f_n$   in \eqref{eq:a335}.
 According to the proof of  Proposition \ref{prop_exist_bdd} \big(in particular, \eqref{eq:a317}\big),
   the BSDEJ\,$(\xi , f_n)$ has a  solution $(Y^n,Z^n,U^n)
 \ins   \hD^\infty   \ti   \hZ^{2,p}      \ti  \hU^p$ such that
 \bea \label{eq:a341}
  \|Y^n\|_{\hD^\infty} \ls R \- 2  .
  \eea
  We also see from    Proposition \ref{prop-a-priori}    that
  \bea \label{eq:b511}
   \|Y^n\|^p_{\hD^p} \+  \|Z^n\|^p_{\hZ^{2,p}} \+ \|U^n\|^p_{\hU^p}
 \ls    \cC  \bigg( 1 \+  \|\xi\|^p_{L^\infty (\cF_T)}  \+
  \Big\| \int_0^T \n   |  f (t, 0,0,0) |   dt \Big\|^p_{L^\infty_+ (\cF_T)}   \bigg)
  \df \wh{\cC} ,
 \eea
 where $ \cC $ is a constant depending on $T$, $\nu(\cX)$, $p$,     $ \ol{C}   $  and $C_\beta$.

 Set $\varpi \df p (1 \- \frac{1}{q'}) \> p (1 \- \frac{1}{q}) \= 1 $
 and let $m,\,n  \ins  \hN$ with $m \> n$.
 Since $ \psi (x) \n \equiv \n 1$ for all $|x| \ls R \- 1$
 and since an analogy to \eqref{eq:a471} shows that
 $  \big\| \pi_m \big( U^m_t  \big)     \-   \pi_n \big( U^n_t  \big)      \big\|_{L^p_\nu} \ls
 2 \big\| U^m_t \- U^n_t  \big\|_{L^p_\nu} \+ \b1_{\{ \| U^n_t \|_{L^p_\nu}>n\}} \big\| U^n_t   \big\|_{L^p_\nu} \ls
 2 \big\| U^m_t \- U^n_t  \big\|_{L^p_\nu} \+ n^{1-\varpi} \big\| U^n_t   \big\|^\varpi_{L^p_\nu} $,
  $t \ins [0,T]$, we can deduce from \eqref{eq:a341}, (H2) and (H4)$-$(H6)      that   \dtp
  \beas
    \q   && \hspace{-1.2cm}  \big|Y^m_t \- Y^n_t\big|^{p-1}
     \big\lan  \sD \big( Y^m_t \- Y^n_t \big)  ,  f_m \big(t,Y^m_t, Z^m_t, U^m_t \big) - f_n\big(t,Y^n_t,
  Z^n_t, U^n_t \big) \big\ran    \nonumber    \\
  &&  \=   \frac{m}{m \ve  \phi^R_t }   \big| Y^m_t \- Y^n_t \big|^{p-1}
     \big\lan  \sD \big( Y^m_t \- Y^n_t \big)  ,
   f \big(t,Y^m_t, \pi_m( Z^m_t), \pi_m(U^m_t) \big)
    \-  f \big(t,Y^n_t, \pi_n( Z^n_t), \pi_n(U^n_t) \big)   \big\ran   \nonumber    \\
   && \q   +    \Big( \frac{m}{m \ve  \phi^R_t }  \- \frac{n}{n \ve  \phi^R_t }  \Big) \big| Y^m_t \- Y^n_t \big|^{p-1}
     \big\lan  \sD \big( Y^m_t \- Y^n_t \big)  ,
   f \big(t,Y^n_t, \pi_n( Z^n_t), \pi_n(U^n_t) \big)  \-   f (t, 0,0,0)    \big\ran  \nonumber    \\
  &&   \le   \l(t) \, \th \big(|Y^m_t \- Y^n_t|^p   \big) \+  \Phi_t      |Y^m_t \- Y^n_t|^p
   \+  |Y^m_t \- Y^n_t|^{p-1} \Big[ \, \L_t \big|\pi_m( Z^m_t)\- \pi_n( Z^n_t) \big|
    \+ c_2 (t) \big\| \pi_m(U^m_t) \- \pi_n(U^n_t) \big\|_{L^p_\nu}    \Big] \\
  && \q +  \Big( 1  \- \frac{n}{n \ve \phi^R_t }  \Big) (2R \- 4 )^{p-1}
     \Big( \phi^R_t \+ \beta_t  \+ c_1(t)  |\pi_n( Z^n_t)| \+ c_2 (t) \|\pi_n(U^n_t)\|_{L^p_\nu} \Big) \\
  && \le  \l(t) \, \th \big(|Y^m_t \- Y^n_t|^p   \big) \+  \Phi_t      |Y^m_t \- Y^n_t|^p
   \+  |Y^m_t \- Y^n_t|^{p-1} \Big[ \, \L_t \big|  Z^m_t \-   Z^n_t  \big|
    \+ 2 c_2 (t) \big\|  U^m_t  \-  U^n_t  \big\|_{L^p_\nu}    \Big] \+ \U^{m,n}_t  ,
 \eeas
 where $ \U^{m,n}_t \df \big( 1  \- \frac{n}{n \vee \phi^R_t }  \big)
  (2R \- 4 )^{p-1} \big( \phi^R_t \+ \beta_t  \+ c_1(t)  |Z^n_t| \+ c_2 (t) \|U^n_t\|_{L^p_\nu} \big)
  \+ (2R \- 4 )^{p-1}  c_2(t) n^{1-\varpi} \big\| U^n_t   \big\|^\varpi_{L^p_\nu} $.
 Thus,  \eqref{nonlip-cond3} holds for   $\eta_n \= 0$, $c(\cd) \= 2 c_2 (\cd)$ and the above process $\U^{m,n}$.
 By H\"older's inequality and \eqref{eq:b511},
 \beas
 \hspace{-0.2cm}    (2R \- 4 )^{1-p} \, E \n \int_0^T \U^{m,n}_t dt
 & \tn \ls & \tn
    E \n \int_0^T \Big( 1  \- \frac{n}{n \ve \phi^R_t }  \Big)  ( \phi^R_t \+ \beta_t )  dt
  \+ \bigg\{ E   \bigg[ \Big( \int_0^T \Big( 1  \- \frac{n}{n \ve \phi^R_t }  \Big)^2 c^2_1(t) dt \Big)^{\frac{q}{2}}
  \bigg] \bigg\}^{\frac{1}{q}} \|Z^n\|_{\hZ^{2,p}} \\
   & \tn  &  \dn + \bigg\{ E \n \int_0^T \Big( 1  \- \frac{n}{n \ve \phi^R_t }  \Big)^q c^q_2(t) dt \bigg\}^{\frac{1}{q}} \|U^n\|_{\hU^p}
  \+    n^{1-\varpi} \Big(\int_0^T \n \big(c_2 (t) \big)^{q'} dt  \Big)^{\frac{1}{q'}}
  \bigg(E \n \int_0^T \n \big\| U^n_t   \big\|^p_{L^p_\nu} dt  \bigg)^{1-\frac{1}{q'}} \\
    & \tn  \ls  & \tn   E \n \int_0^T \Big( 1  \- \frac{n}{n \ve \phi^R_t }  \Big)  ( \phi^R_t \+ \beta_t )  dt
  \+ \wh{\cC}^{\frac{1}{p}} \bigg\{ E   \bigg[ \Big( \int_0^T \Big( 1  \- \frac{n}{n \ve \phi^R_t }  \Big)^2 c^2_1(t) dt \Big)^{\frac{q}{2}}  \bigg] \bigg\}^{\frac{1}{q}} \\
    & \tn  &  \dn  +  \wh{\cC}^{\frac{1}{p}}  \bigg\{ E \n \int_0^T \Big( 1  \- \frac{n}{n \ve \phi^R_t }  \Big)^q c^q_2(t) dt \bigg\}^{\frac{1}{q}}
   \+   \wh{\cC}^{1 \- \frac{1}{q'}}  n^{1-\varpi} \Big(\int_0^T \n \big(c_2 (t) \big)^{q'} dt  \Big)^{\frac{1}{q'}}
  \df I_n .
 \eeas

 Because $  \frac{n}{n \vee \phi^R_t } \= \frac{1}{1 \vee (\phi^R_t / n) }  \n \nearrow \n  1$
 as $n  \n \to \n  \infty$, $\fa t  \ins  [0,T]$, the dominated convergence theorem shows that
 $     \lmt{n \to \infty} I_n  \= 0$.
    It follows  that $ \lmt{n \to \infty} \, \underset{m>n}{\sup} \, E \n   \int_0^T \U^{m,n}_t dt \= 0  $.
 Since $ \underset{n \in \hN}{\sup} \Big( \|Y^n\|^p_{\hD^p} \+  \|Z^n\|^p_{\hZ^{2,p}} \+ \|U^n\|^p_{\hU^p} \Big)
     \ls    \wh{\cC}  $  by \eqref{eq:b511}, we see from
   Proposition \ref{prop_stab}   that  $\big\{(Y^n, Z^n, U^n)\big\}_{n \in \hN} $
 is a Cauchy sequence in $ \hS^p$.
 Let $(Y, Z, U)$  be its limit.
 As in the proof of Proposition \ref{prop_exist_bdd},
 one can extract a subsequence $\{m_i\}_{i \in \hN}$ from $\hN$    such that
 \eqref{eqn-b356}$-$\eqref{eq:b515}   hold, then we still have     \eqref{eq:a317}.
 Also,   similar to \eqref{eq:a425} and \eqref{eq:a426},
 we can define  two $[0,\infty)-$valued, $\bF-$predictable processes $\cZ$ and $\cU$
 that satisfy \eqref{eqn-b513} and  \eqref{eqn-b513d} respectively.

  Fix $ k \ins  \hN$ and   define the $\bF-$stopping time $\tau_k$ as in \eqref{def_tau_k}.
 We can still derive \eqref{eqn-b362} from  the dominated convergence theorem and \eqref{eqn-b356}.
       Hence,  there exists a subsequence $\big\{m^k_i \big\}_{i \in \hN}$ of $\{m_i\}_{i \in \hN}$
   such that \eqref{eqn-b521} holds   \dtp ~
     For any $(t,\o) \ins [0,T]  \ti  \O$  except on a $dt  \ti  dP-$null  set $ \fN_k $,
  we may assume that    (H2), (H4), (H5), \eqref{eqn-b356}, \eqref{eq:a411}, \eqref{eqn-b521} hold,
  that $|Y_t (\o)|  \ls  R \- 2 $,  $|Y^{m^k_i}_t (\o)|  \ls  R \- 2 $,   $ \fa  i   \ins  \hN$
  \big(by \eqref{eq:a341}, \eqref{eq:a317}\big),
  and that $ U_t (\o)  \ins  L^p_\nu $, $ U^{m^k_i}_t(\o)  \ins  L^p_\nu $, $\fa i  \ins  \hN$.

    Let  $(t, \o)  \ins   \fN^c_k \n \cap \n \[0, \tau_k\]  $.
    Since $\lmtu{i \to \infty} \frac{m^k_i}{m^k_i \vee \phi^R_t } \= 1$
    and since $\psi \big( Y^{m^k_i}_t(\o) \big) \= 1$,
 \bea   \label{eqn-b525b}
 \lmt{i \to \infty}    f_{m^k_i} \Big(t,\o,Y^{m^k_i}_t(\o) , Z^{m^k_i}_t(\o)  ,U^{m^k_i}_t(\o) \Big)
 \=  \lmt{i \to \infty}    f \Big(t,\o,Y^{m^k_i}_t(\o) , \pi_{m^k_i} \Big( Z^{m^k_i}_t(\o) \Big) ,
 \pi_{m^k_i} \Big(U^{m^k_i}_t(\o) \Big) \Big)  .
 \eea
  Using (H5), \eqref{eq:a411} and
   an analogy to the inequality (e3) in Part 5 of Proposition \ref{prop_exist_bdd}'s proof, we obtain
 \bea
 && \hspace{-1.2cm}  \bigg|   f  \Big(  t,\o,Y^{m^k_i}_t(\o) , \pi_{m^k_i} \Big( Z^{m^k_i}_t(\o) \Big), \pi_{m^k_i} \Big(U^{m^k_i}_t(\o) \Big)  \Big)
    \-       f   \Big( t,\o,Y^{m^k_i}_t(\o) , \pi_{m^k_i} \Big( Z^{m^k_i}_t(\o) \Big),U_t(\o)  \Big) \Big| \nonumber \\
 && \hspace{-0.7cm}   \ls     c_2(t)    \Big\|  \pi_{m^k_i} \Big( U^{m^k_i}_t(\o) \Big)  \-   U_t (\o)   \Big\|_{L^p_\nu}
    \ls  c_2(t)  \Big( 2 \Big\|    U^{m^k_i}_t(\o)    \-   U_t (\o)   \Big\|_{L^p_\nu}
   \+     \big\|  \pi_{m^k_i} \big( U_t (\o) \big) \-  U_t (\o) \big\|_{L^p_\nu} \Big)
   \n \to \n 0 , ~ \hb{as } i \to \infty ,   \qq
    \eea
 Also, similar to the inequality (e3) in Part 5 of Proposition \ref{prop_exist_bdd}'s proof,
 one can deduce from \eqref{eqn-d051}  and \eqref{eq:a347}   that
 $ \Big|  \pi_{m^k_i} \Big( Z^{m^k_i}_t(\o) \Big) \- Z_t(\o) \Big|
 \ls \Big| Z^{m^k_i}_t(\o)   \-   Z_t(\o) \Big| \+ \Big|  \pi_{m^k_i} \big( Z_t(\o) \big) \- Z_t(\o) \Big|
 \to 0 $ as $i \to \infty$, which together with  \eqref{eqn-b356}
 and  the  continuity of the mapping $     f \big(t, \o, \cd,\cd,  U_t (\o) \big)$, shows  that
        \bea   \label{eq:a473}
  \lmt{i \to \infty}    f \Big(t,\o,Y^{m^k_i}_t(\o) , \pi_{m^k_i} \Big( Z^{m^k_i}_t(\o) \Big) ,  U_t(\o)   \Big)
 \=  f  \big(t,\o, Y_t(\o), Z_t(\o), U_t(\o) \big) .   \q
 \eea
 Combining \eqref{eqn-b525b}$-$\eqref{eq:a473} leads to that
  \bea  \label{eqn-b523b}
 \lmt{i \to \infty}  \Big|  f_{m^k_i} \Big(  t,\o,Y^{m^k_i}_t(\o) , Z^{m^k_i}_t(\o),  U^{m^k_i}_t(\o)   \Big)
  \-   f  \big(t,\o, Y_t(\o), Z_t(\o), U_t(\o) \big) \Big| \= 0 .
 \eea

 Given $ i \ins  \hN $, since $\psi \big( Y^{m^k_i}_t(\o) \big) \= 1$,
 one can deduce from  (H2), (H4), (H5) and an analogy to \eqref{eq:a477}  that
 \beas
 && \hspace{-1.2cm}  \Big|  f_{m^k_i} \Big(  t,\o,Y^{m^k_i}_t(\o) , Z^{m^k_i}_t(\o),  U^{m^k_i}_t(\o)   \Big)
   \-  f \big(t, \o,  Y_t(\o),  Z_t(\o), U_t(\o)\big) \Big| \\
   && \hspace{-0.7cm} \= \bigg| \frac{m^k_i}{m^k_i \ve \phi^R_t (\o) } 
   \Big(  f   \Big(  t,\o,Y^{m^k_i}_t(\o) , \pi_{m^k_i} \Big( Z^{m^k_i}_t(\o)\Big), \pi_{m^k_i} \Big( U^{m^k_i}_t(\o) \Big)  \Big) \- f (t,\o,0,0,0) \Big) \\
  && \hspace{-0.3cm}
   +  f (t,\o,0,0,0) \- f \big(t, \o,  Y_t(\o),  Z_t(\o), U_t(\o)\big)  \bigg| \\
   && \hspace{-0.7cm} \ls
   \Big|  f   \Big(  t,\o,Y^{m^k_i}_t(\o) , \pi_{m^k_i} \Big( Z^{m^k_i}_t(\o)\Big), \pi_{m^k_i} \Big( U^{m^k_i}_t(\o) \Big)  \Big) \- f (t,\o,0,0,0) \Big|
   \+  \big|  f \big(t, \o,  Y_t(\o),  Z_t(\o), U_t(\o)\big) \- f (t,\o,0,0,0) \big| \\
 & & \hspace{-0.7cm} \ls  2 \phi^R_t (\o)   \+ 2 \beta_t(\o)
  \+  c_1(t) \Big(    \Big| Z^{m^k_i}_t (\o)     \Big| \+ |Z_t(\o)|   \Big)
  \+  c_2(t) \Big(    \Big\|  \pi_{m^k_i} \Big( U^{m^k_i}_t(\o) \Big) \Big\|_{L^p_\nu}
 \+   \|U_t(\o)\|_{L^p_\nu} \Big) \\
 & & \hspace{-0.7cm} \ls   2 \phi^R_t (\o)  \+ 2 \beta_t(\o)
 \+ c_1(t)     \big(  \cZ_t (\o) \+ |Z_t(\o)| \big) \+ c_2(t)     \big(  \cU_t (\o)
 \+   \|U_t(\o)\|_{L^p_\nu}    \big) \df H_t (\o)  .   \qq
 \eeas
 Analogous  to  \eqref{eqn-b539},
 we can deduce from Holder's inequality,   \eqref{eqn-b513} and \eqref{eqn-b513d} that
    \beas
 E \n \int_0^T \n  H_t dt \ls 2 E \n \int_0^T \n  \phi^R_t dt \+ 2 C_\beta  \+ \ol{C}^{\frac12}
 \Bigg\{ \bigg(E \bigg[ \Big(\int_0^T \n \cZ^2_t dt
 \Big)^{\frac{p}{2}} \bigg]\bigg)^{\frac{1}{p}} + \|Z\|_{\hZ^{2,p}} \Bigg\}
 \+ \ol{C}^{\frac{1}{q}}   \bigg\{
 \Big( E \int_0^T  \n  \cU^p_t dt       \Big)^{\frac{1}{p}} \+  \|U\|_{\hU^p} \bigg\}  < \infty .
 \eeas
 The  dominated convergence theorem and \eqref{eqn-b523b} yield that
     \beas
   \lmt{i \to \infty}  E \n \int_0^{\tau_k} \Big|   f_{m^k_i}    \Big(t, Y^{m^k_i}_t, Z^{m^k_i}_t, U^{m^k_i}_t\Big)
   \- f  (t, Y_t, Z_t,  U_t ) \Big| dt   \= 0 .  
  \eeas
  \if{0}
 As the two limits \eqref{eqn-b364} and  \eqref{eqn-b366} still hold,
 one can extract a subsequence $\big\{\wt{m}^k_i \big\}_{i \in \hN}$ from $\big\{m^k_i \big\}_{i \in \hN}$   such that
 \eqref{eq:b531} holds   except on a $P-$null set $\cN^k_1$.

     Since  $ \Big(Y^{\wt{m}^k_i}, Z^{\wt{m}^k_i}, U^{\wt{m}^k_i} \Big) $
      solves   BSDEJ\,$\big( \xi , f_{\wt{m}^k_i} \big)$ for any $i \ins  \hN$, \eqref{eq:a345} holds
       except on a $P-$null set $\cN^k_2$.
  For any $   \o \ins  \O_k  \df  \big( \cN^k_1 \big)^c  \n  \cap  \n  \big( \cN^k_2 \big)^c  $,
  letting $i  \n \to \n  \infty$ in \eqref{eq:a345},
  we obtain \eqref{eqn-b376} again  over $\O_k $.
  By \eqref{eqn-b513}, it holds   for all $\o \ins \O$ except on $P-$null set $\cN_Z$ that
  $    \int_0^T  \cZ_t^2 (\o) dt      \<   \infty $, and thus $ \t_{\fk}(\o) \=  T $
  for some $ \fk \= \fk(\o) \ins  \hN$.
  Then letting $k  \n \to \n \infty$ in \eqref{eqn-b376} shows that \eqref{BSDEJ} holds  over
  $  \Big(  \underset{k \in \hN}{\cap}  \O_k  \Big) \n \cap   \cN^c_Z $,
  which together with Remark \ref{rem_drift} shows that   $(Y, Z, U) $ is a solution of BSDEJ\,$(\xi,  f)$.
  \fi
  Then following similar arguments to Part 5 in the proof of Proposition \ref{prop_exist_bdd}, one can show that
   $(Y, Z, U) $ is a solution of BSDEJ\,$(\xi,  f)$.

 \no {\bf 2)} Next, let consider the general case that
 $ \xi \ins L^p (\cF_T) $ and $  \int_0^T \n  |f(t,0,0,0)| dt   \ins  L^p_+(\cF_T) $.
    For any $n \ins  \hN$, we set $\xi_n  \df  \pi_n (\xi )$ and define
 \beas
 \q \wt{f}_n (t,\o,y,z,u)  \df   f (t,\o,y,z,u)  \-   f (t,\o,0,0,0)    \+  \pi_n \big( f (t,\o,0,0,0) \big) ,
 ~ (t,\o,y,z,u)  \ins   [0, T]  \ti  \O  \ti  \hR^l  \ti  \hR^{l \times d}
  \ti   L^p_\nu  .
 \eeas
 Clearly, $\wt{f}_n$ is  an $   \hR^l-$valued,
  $\sP  \oti  \sB  ( \hR^l   )   \oti  \sB \big( \hR^{l\times d} \big)
  \oti  \sB\big(L^p_\nu   \big)/\sB(\hR^l ) -$measurable mapping   satisfying   (H1)$-$(H6)
   with the same   coefficients as $f$.
 As $\int_0^T |\wt{f}_n (t,0,0,0)|dt  \=  \int_0^T \big|\pi_n ( f (t, 0,0,0)  )\big|dt  \ls  nT $,
  Part 1 shows that
   the BSDEJ\,$\big( \xi_n,\wt{f}_n \big)$ has a  solution $ \big(Y^n,Z^n,U^n \big)
 \ins   \hD^\infty   \ti   \hZ^{2,p}      \ti  \hU^p$ \big(For easy reference, we still denote this solution by
 $  (Y^n,Z^n,U^n  ) $. Note its difference from the triple   considered in Part 1\big).
 Also,  we know from Proposition \ref{prop-a-priori}   that
  \bea
   \|Y^n\|^p_{\hD^p} \+  \|Z^n\|^p_{\hZ^{2,p}} \+ \|U^n\|^p_{\hU^p}
& \tn  \ls  & \tn   \cC      E\bigg[ 1 \+ |\xi_n|^p
 \+ \Big( \int_0^T \n \big|\pi_n  ( f (t, 0,0,0)  ) \big| dt \Big)^p \bigg] \nonumber \\
& \tn  \ls  & \tn   \cC      E\bigg[ 1 \+  | \xi  |^p
 \+ \Big( \int_0^T \n \big|  f (t, 0,0,0)   \big| dt \Big)^p \bigg] \df \wt{\cC}  , \q \label{eq:b511d}
 \eea
 where $ \cC $ is a constant depending on $T$, $\nu(\cX)$, $p$,     $ \ol{C}   $
  and $C_\beta$.

       Given $m,\,n  \ins  \hN$ with $m \> n$, an analogy to  \eqref{eq:b411}  shows  that
 \eqref{nonlip-cond3} holds for $f_n \= \wt{f}_n$,  $\eta_n \= 0$, $c(\cd) \= c_2(\cd)$ and
 \beas
 \U^{m,n}_t  \= \big|Y^m_t \- Y^n_t\big|^{p-1}
  \big|\pi_m \big( f (t, 0,0,0) \big) \- \pi_n \big( f (t, 0,0,0) \big) \big| , \q \fa t \ins [0,T] .
  \eeas
  By H\"older's inequality and \eqref{eq:b511d},
  \beas
  E \int_0^T \U^{m,n}_t dt
  & \tn \ls  & \tn  E \bigg[ \underset{t \in [0,T]}{\sup} \big|Y^m_t \- Y^n_t\big|^{p-1} \int_0^T \big|  f (t, 0,0,0)   \- \pi_n \big( f (t, 0,0,0) \big) \big| dt \bigg] \\
   & \tn \ls  & \tn  \big\| Y^m  \- Y^n \big\|_{\hD^p}^{\frac{p}{q}}
  \bigg\{ E \bigg[ \Big(\int_0^T \big|  f (t, 0,0,0)   \- \pi_n \big( f (t, 0,0,0) \big) \big| dt \Big)^p \bigg] \bigg\}^{\frac{1}{p}} \\
   & \tn \ls & \tn    2^{\frac{p}{q}} \wt{\cC}^{\frac{1}{q}}
   \bigg\{ E \bigg[ \Big(\int_0^T \big|  f (t, 0,0,0)   \- \pi_n \big( f (t, 0,0,0) \big) \big| dt \Big)^p \bigg] \bigg\}^{\frac{1}{p}} \df \wt{I}_n .
  \eeas

  As $  E \Big[ \big(\int_0^T  |  f (t, 0,0,0)    | dt \big)^p \Big] \< \infty $,
  the dominated convergence theorem implies that  $\lmt{n \to \infty} \wt{I}_n \= 0$.
  It follows that   $ \lmt{n \to \infty} \, \underset{m>n}{\sup} \, E \n   \int_0^T \U^{m,n}_t dt \= 0  $.
  Since $ \underset{n \in \hN}{\sup} \Big( \big\|Y^n \big\|^p_{\hD^p}
 \+  \big\| Z^n \big\|^p_{\hZ^{2,p}} \+ \big\| U^n \big\|^p_{\hU^p} \Big)
     \ls  \wt{\cC} $
  by \eqref{eq:b511d}, we see from
   Proposition \ref{prop_stab}   that  $\big\{(Y^n, Z^n, U^n)\big\}_{n \in \hN} $
 is a Cauchy sequence in $ \hS^p$.
 Let $(Y, Z, U)$  be its limit.
 As in the proof of Proposition \ref{prop_exist_bdd},
 one can extract a subsequence $\{m_i\}_{i \in \hN}$ from $\hN$    such that
 \eqref{eqn-b356}$-$\eqref{eq:b515}   hold.

           For any $i \ins  \hN$, we define an $\cF_T-$measurable
     random variable  $ \dis \fy_i  \df Y_*   \+  \sum_{j=1}^i
    \big( Y^{m_j}   \-  Y^{m_{j-1}}   \big)_* $  with $ Y^{n_0}  \df   Y$.
    Minkowski's inequality and \eqref{eq:b515} show that
  \bea
   \big\| \fy_i \big\|_{L^p_+(\cF_T)}
    \ls  
 \big\|  Y \big\|_{\hD^p }  \+  \sum_{j=1}^i     \big\|Y^{m_j}  \dn - \n  Y^{m_{j-1}} \big\|_{\hD^p }
        \ls   1 \+ \big\|  Y \big\|_{\hD^p}  \+    \big\| Y^{m_1}  \dn - \n  Y \big\|_{\hD^p}  .
        \label{eqn-b360b}
        \eea
   Since $\big\{\fy_i\big\}_{i \in \hN}$ is an increasing sequence,
     $\dis \fy  \df  \lmtu{i \to \infty} \fy_i
        \=    Y_*    + \n  \sum_{j=1}^\infty    \big(Y^{m_j} \- Y^{m_{j-1}}  \big)_*   $
        defines a $[0,\infty]-$valued, $\cF_T-$measurable random variable.
    Then the monotone convergence theorem and   \eqref{eqn-b360b} imply  that
    \bea
    \big\| \fy \big\|_{L^p_+(\cF_T)}    \=   \lmtu{i \to \infty}  \big\| \fy_i \big\|_{L^p_+(\cF_T)}
     \ls   1\n +\n \big\|  Y \big\|_{\hD^p}
     \+    \big\| Y^{n_1}  \-  Y \big\|_{\hD^p}   \<  \infty.   \label{eqn-b513b}
    \eea
 Moreover, as in \eqref{eq:a425} and \eqref{eq:a426},
 we can define  two $[0,\infty)-$valued, $\bF-$predictable processes $\cZ$ and $\cU$
 that satisfy \eqref{eqn-b513} and  \eqref{eqn-b513d} respectively.

 Fix $ k \ins  \hN$ and   define the $\bF-$stopping time $\tau_k$ as in \eqref{def_tau_k}.
 One can again derive \eqref{eqn-b362} from  the dominated convergence theorem and \eqref{eqn-b356}.
       Hence,  there exists a subsequence $\big\{m^k_i \big\}_{i \in \hN}$ of $\{m_i\}_{i \in \hN}$
   such that \eqref{eqn-b521} holds   \dtp ~
    For any $(t,\o) \ins [0,T]  \ti  \O$  except on a $dt  \ti  dP-$null  set $ \wt{\fN}_k $,
  we may assume that    (H2), (H4), (H5), \eqref{eqn-b356}, \eqref{eq:a411}, \eqref{eqn-b521} hold
  and that $ U_t (\o)  \ins  L^p_\nu $, $ U^{m^k_i}_t(\o)  \ins  L^p_\nu $, $\fa i  \ins  \hN$.

   Let us also fix $ \ell \ins \hN$ and define $A_\ell \df \{ \fy \vee Y_* \ls \ell \} \ins \cF_T$.

   Let $(t, \o)  \ins    \wt{\fN}_k^c \n \cap \n \[0, \tau_k\]  $.
        The  continuity of the mapping $     f \big(t, \o, \cd,\cd,  U_t (\o) \big)$,
       \eqref{eqn-b356} and \eqref{eq:a347} yield    that
  \bea   \label{eqn-b525d}
  \hspace{-5mm} \lmt{i \to \infty}    \wt{f}_{m^k_i} \Big(t,\o,Y^{m^k_i}_t(\o) , Z^{m^k_i}_t(\o)  ,U_t(\o) \Big)
  \=   \lmt{i \to \infty}    f \Big(t,\o,Y^{m^k_i}_t(\o) , Z^{m^k_i}_t(\o)  ,U_t(\o) \Big)
  \=     f  \big(t,\o, Y_t(\o), Z_t(\o), U_t(\o) \big) . \q
  \eea
  By  (H5), it holds for any $i \ins  \hN$ that
 \beas
  \q && \hspace{-1.2cm}  \Big|   \wt{f}_{m^k_i} \Big(  t,\o,Y^{m^k_i}_t(\o) , Z^{m^k_i}_t(\o),  U^{m^k_i}_t(\o)   \Big)
    \-       \wt{f}_{m^k_i}  \Big( t,\o,Y^{m^k_i}_t(\o) , Z^{m^k_i}_t(\o),U_t(\o)  \Big) \Big| \\
   &&  \=  \Big|   f  \Big(  t,\o,Y^{m^k_i}_t(\o) , Z^{m^k_i}_t(\o),  U^{m^k_i}_t(\o)   \Big)
    \-       f   \Big( t,\o,Y^{m^k_i}_t(\o) , Z^{m^k_i}_t(\o),U_t(\o)  \Big) \Big|
    \ls     c_2(t)    \Big\| U^{m^k_i}_t (\o)  \-  U_t (\o) \Big\|_{L^p_\nu}  ,
    \eeas
    which together with \eqref{eq:a411} and \eqref{eqn-b525d} shows that
  \bea  \label{eqn-b523d}
 \lmt{i \to \infty}  \Big|  \wt{f}_{m^k_i} \Big(  t,\o,Y^{m^k_i}_t(\o) , Z^{m^k_i}_t(\o),  U^{m^k_i}_t(\o)   \Big)
  \-   f  \big(t,\o, Y_t(\o), Z_t(\o), U_t(\o) \big) \Big| \= 0 .
 \eea

  Given $ i \ins  \hN $, there exists an $j = j (k,i) \ins \hN$ such that
 $m^k_i \= m_{j}$.  Since
 $ Y^{m^k_i}_*(\o) \ls  Y_*(\o)  \+  \sum_{j=1}^i   \n  \big(Y^{m^k_j}
  \-  Y^{m^k_{j-1}}   \big)_*(\o)   \=  \fy_{j} (\o)  \ls    \fy(\o)     $,
    one can deduce from  (H2), (H4), (H5) and an analogy to \eqref{eq:a477} that
 \beas
 && \hspace{-1.2cm} \b1_{A_\ell} \Big|  \wt{f}_{m^k_i} \Big(  t,\o,Y^{m^k_i}_t(\o) , Z^{m^k_i}_t(\o),  U^{m^k_i}_t(\o)   \Big)
   \-  f \big(t, \o,  Y_t(\o),  Z_t(\o), U_t(\o)\big) \Big| \\
 && \hspace{-0.7cm}  \ls \big| f(t,\o,0,0,0)  \- \pi_n ( f(t,\o,0,0,0) )  \big| \+
  \b1_{A_\ell} \Big| f \Big(  t,\o,Y^{m^k_i}_t(\o) , Z^{m^k_i}_t(\o),  U^{m^k_i}_t(\o)   \Big)
 \-   f \big(t, Y_t(\o),  Z_t(\o), U_t(\o)\big)   \Big| \\
 && \hspace{-0.7cm}  \ls \big| f(t,\o,0,0,0)   \big|
 \+ \b1_{A_\ell} \bigg\{ \Big| f \Big(  t,\o,Y^{m^k_i}_t(\o) , Z^{m^k_i}_t(\o),  U^{m^k_i}_t(\o)   \Big)
 \- f(t,\o,0,0,0) \Big|
   \+   \big| f \big(t, Y_t(\o),  Z_t(\o), U_t(\o)\big)  \- f(t,\o,0,0,0)  \big| \bigg\}  \\
& & \hspace{-0.7cm} \ls   \big| f(t,\o,0,0,0)   \big| \+  2 \phi^\ell_t(\o)
  \+ 2 \beta_t(\o) \+  c_1(t) \Big(    \Big| Z^{m^k_i}_t (\o)     \Big| \+ |Z_t(\o)|   \Big)
  \+  c_2(t) \Big(    \Big\|   U^{m^k_i}_t (\o)  \Big\|_{L^p_\nu}
 \+   \|U_t(\o)\|_{L^p_\nu} \Big) \\
 & & \hspace{-0.7cm} \ls   \big| f(t,\o,0,0,0)   \big| \+      2 \phi^\ell_t(\o)
  \+ 2 \beta_t(\o) \+  c_1(t)     \big(  \cZ_t (\o) \+ |Z_t(\o)| \big) \+ c_2(t)     \big(  \cU_t (\o)
 \+   \|U_t(\o)\|_{L^p_\nu}    \big) \df H^\ell_t (\o)  .   \q
 \eeas
 Similar to  \eqref{eqn-b539},
 we can deduce from Holder's inequality,  \eqref{eqn-b513} and \eqref{eqn-b513d} that
 \beas
 \hspace{-7mm}  E \n \int_0^T \n  H^\ell_t dt \ls E \n \int_0^T \n \big( | f(t, 0,0,0)   |
  \+  2   \phi^\ell_t \big) dt  \+ 2 C_\beta  \+ \ol{C}^{\frac12}
 \Bigg\{ \bigg(E \bigg[ \Big(\int_0^T \n \cZ^2_t dt
 \Big)^{\frac{p}{2}} \bigg]\bigg)^{\frac{1}{p}} \+ \|Z\|_{\hZ^{2,p}} \Bigg\}
 \+ \ol{C}^{\frac{1}{q}}   \bigg\{
 \Big( E \int_0^T  \n  \cU^p_t dt       \Big)^{\frac{1}{p}} \+  \|U\|_{\hU^p} \bigg\}  \< \infty .
 \eeas
 Then the  dominated convergence theorem and \eqref{eqn-b523d} yield that
     \bea
   \lmt{i \to \infty}  E \n \int_0^{\tau_k} \n \b1_{A_\ell}
    \Big|   \wt{f}_{m^k_i}    \Big(t, Y^{m^k_i}_t, Z^{m^k_i}_t, U^{m^k_i}_t\Big)
   \- f  (t, Y_t, Z_t,  U_t ) \Big| dt   \= 0 .   \label{eqn-b519d}
  \eea

     In light of \eqref{eqn-b356}, \eqref{eqn-b519d},   \eqref{eqn-b364} and  \eqref{eqn-b366},
     there exists a subsequence $\big\{m^{k,\ell}_i \big\}_{i \in \hN}$ of $\big\{m^k_i \big\}_{i \in \hN}$   such that
    except on a $P-$null set $\cN^{k,\ell}_1$
    \beas
   \q && \hspace{-2cm}  \lmt{i \to \infty}   \bigg\{ \underset{t \in [0,T]}{\sup} \Big|Y^{m^{k,\ell}_i}_t - Y_t  \Big|
   +   \b1_{A_\ell} \n   \int_0^{  \t_k}  \left|   \wt{f}_{m^{k,\ell}_i}   \Big(t, Y^{m^{k,\ell}_i}_t, Z^{m^{k,\ell}_i}_t, U^{m^{k,\ell}_i}_t\Big)
   - f \big(t, Y_t, Z_t,  U_t\big) \right| dt   \nonumber  \\
  && +    \underset{t \in [0,T]}{\sup} \bigg|  \int_{\tau_k \land t}^{\t_k} \dn
  \Big(Z^{m^{k,\ell}_i}_s \- Z_s\Big) dB_s \bigg|
    +   \underset{t \in [0,T]}{\sup}  \bigg|  \int_{( \tau_k \land t, \t_k]}  \n \int_\cX \n \Big(U^{m^{k,\ell}_i}_s \n (x)
   \- U_s(x)\Big)\tnp(ds,dx) \bigg|   \bigg\}   \n  =0.
    \eeas

     Since  $ \Big(Y^{m^{k,\ell}_i}, Z^{m^{k,\ell}_i}, U^{m^{k,\ell}_i} \Big) $
      solves   BSDEJ\,$\Big( \xi_{m^{k,\ell}_i}, \wt{f}_{m^{k,\ell}_i} \Big)$ for any $i \ins  \hN$, it holds
       except on a $P-$null set $\cN^{k,\ell}_2$ that
     \bea
 Y^{m^{k,\ell}_i}_{\tau_k \land t}  & \tn = &  \tn  \b1_{\{\tau_k < T\}}Y^{m^{k,\ell}_i}_{  \tau_k}
 + \b1_{\{\tau_k = T\}} \pi_{m^{k,\ell}_i} (\xi)
 + \int_{\tau_k \land t}^{  \tau_k} \wt{f}_{m^{k,\ell}_i} \Big(s, Y^{m^{k,\ell}_i}_s, Z^{m^{k,\ell}_i}_s, U^{m^{k,\ell}_i}_s\Big) \, ds -\int_{\tau_k \land t}^{  \tau_k}  Z^{m^{k,\ell}_i}_s \, dB_s  \nonumber \\
    &&   -\int_{(\tau_k \land t, \tau_k]} \int_\cX U^{m^{k,\ell}_i}_s(x)\tnp(ds,dx),
    \qq    \fa    t \in  [0,T],   \q \;\;  \fa i \in \hN.  \label{eqn-b368c}
   \eea
   Set  $\wt{A}^k_\ell  \df  \big( \cN^{k,\ell}_1 \big)^c  \n  \cap  \n   \big( \cN^{k,\ell}_2 \big)^c \n \cap \n A_\ell $,
   which includes  the set   $ \Big( \underset{\ell \in \hN}{\cup}  \cN^{k,\ell}_1 \Big)^c
   \n  \cap  \n   \Big( \underset{\ell \in \hN}{\cup}  \cN^{k,\ell}_2 \Big)^c \n \cap \n A_\ell  $.
  For any $   \o \ins  \wt{A}^k_\ell   $,   letting $i  \n \to \n  \infty$ in \eqref{eqn-b368c},
  we obtain \eqref{eqn-b376}  over $ \wt{A}^k_\ell $.
  As $\ell$ varies over $\hN$, \eqref{eqn-b376} further holds over $ \O_k \df \Big( \underset{\ell \in \hN}{\cup}  \cN^{k,\ell}_1 \Big)^c
   \n  \cap  \n   \Big( \underset{\ell \in \hN}{\cup}  \cN^{k,\ell}_2 \Big)^c \n \cap
   \Big( \underset{\ell \in \hN}{\cup} \n A_\ell \Big) $.
   By \eqref{eqn-b513b}  and $Y \ins \hD^p$,
  one has $P \big(\O_k\big) \= P \Big( \underset{\ell \in \hN}{\cup} A_\ell \Big) \= 1$.

  We see from \eqref{eqn-b513}  that for all $\o \ins \O$ except on $P-$null set $\cN_Z$,
  $    \int_0^T  \cZ_t^2 (\o) dt      \<   \infty $  and thus $ \t_{\fk}(\o) \=  T $
  for some $ \fk \= \fk(\o) \ins  \hN$.
  Then letting $k  \n \to \n \infty$ in \eqref{eqn-b376} shows that \eqref{BSDEJ} holds  over
  $  \Big(  \underset{k \in \hN}{\cap}  \O_k  \Big) \n \cap   \cN^c_Z $,
  which together with Remark \ref{rem_drift} shows that   $(Y, Z, U) $ is a solution of BSDEJ\,$(\xi,  f)$.   \qed

 \if{0}

    \no {\bf Proof of Corollary \ref{cor-case1}:}
   As $f_\tau$ is a $\sP \otimes \sB \big( \hR^l  \big) \otimes \sB \big( \hR^{l\times d} \big)
\otimes \sB\big(L^p_\nu   \big)/\sB(\hR^l ) -$measurable function   satisfying
 (H1)-(H6), Theorem \ref{thm_BSDEJ1}   shows that the BSDEJ\,$\big(\xi,f_\tau \big)$ admits
 a unique solution $(Y, Z, U) \in  \hS^p  $.  So it holds  except on a $P-$null set $\cN_1$ that
   \bea   \label{eqn-b621}
 Y_t=\xi+\int_t^T f_\tau(s, Y_s, Z_s, U_s) ds -\int_t^T Z_s \, dB_s-\int_{(t,T]} \int_\cX U_s(x)\tnp(ds,dx), \q t\in  [0,T].
  \eea

   Fix $T \in (0,T]$ and $n \in \hN$.  We  define
 \bea
 \ga_n   := \inf\left\{t \in [0, T]:
   \int_0^t \left( |Z_s|^2 +  \|U_s \|^2_{L^p_\nu} \right) ds > n\right\}  \ins \cT  .      \label{sp_tau_n3}
 \eea
  By   \eqref{eqn-b621},  it holds  on $\cN^c_1$ that
 \beas
 Y_{\tau \land \ga_n \land T} &=& Y_{  \ga_n \land T}+ \int_{\tau \land \ga_n \land T}^{  \ga_n \land T} f_\tau(s, Y_s, Z_s, U_s) ds
   -\int_{\tau \land \ga_n \land T}^{  \ga_n \land T} Z_s \, dB_s-\int_{ (\tau \land \ga_n \land T,  \ga_n \land T]} \int_\cX U_s(x)\tnp(ds,dx)  \\
 &=& Y_{  \ga_n \land T}    -\int_{\tau \land \ga_n \land T}^{  \ga_n \land T} Z_s \, dB_s-\int_{ (\tau \land \ga_n \land T,  \ga_n \land T]} \int_\cX U_s(x)\tnp(ds,dx)   .
 \eeas
  Taking conditional expectation $E\big[\cd\big|\cF_{ \tau \land \ga_n \land T }\big]$ and  multiplying $\b1_{\{\tau \le  \ga_n \land T\}}$ on both sides yield  that
  \bea   \label{eqn-b625}
    \b1_{\{\tau \le  \ga_n \land T\}} Y_\tau     =   \b1_{\{\tau \le  \ga_n \land T\}} Y_{\tau \land \ga_n \land T}
 =\b1_{\{\tau \le  \ga_n \land T\}} E\big[Y_{  \ga_n \land T}\big|\cF_{ \tau \land \ga_n \land T }\big]
 =\b1_{\{\tau \le  \ga_n \land T\}} E\big[Y_{  \ga_n \land T}\big|\cF_\tau \big]     , \q \pas
  \eea
As $ (Z , U  ) \in \hZ^{2,p}     \times \hU^p $, we see that $\int_0^T \big( |Z_s|^2 +  \|U_s \|^2_{L^p_\nu} \big) ds < \infty$, \pas ~
    Thus  for \pas~$\o \in \O$,    $ \ga_n (\o)=\infty $ for some   $n = n(\o) \in \hN$, which implies  that
  $   \lmt{n \to \infty} Y_{\ga_n \land T} =   Y_T   $,    \pas ~although the process $Y $ may not be left-continuous.
  Since $ E\Big[\underset{t \in [0, T] }{\sup}|Y_t| \Big]
 \ls   \| Y \|_{\hD^p }  \<  \infty   $ by H\"older's inequality,
 the dominated convergence theorem   implies that
  \beas
   \lmt{n \to \infty} E\big[Y_{  \ga_n \land T}\big|\cF_\tau \big] = E\big[Y_T  \big|\cF_\tau \big], \q \hb{and} \q
 \lmt{T \to \infty}   E\big[Y_T  \big|\cF_\tau \big] =E\big[\xi  \big|\cF_\tau \big] =\xi   , \q \pas
  \eeas
 It is clear that  $ \lmtu{n \to \infty} \b1_{\{\tau \le  \ga_n \land T\}} = \b1_{\{\tau \le  T  \}}$ and that $ \lmtu{T \to \infty} \b1_{\{\tau \le  T  \}} = \b1_{\{\tau < \infty  \}}$. Thus, letting $n \to \infty $ and then letting  $T \to \infty $ in \eqref{eqn-b625} give  that $\b1_{\{\tau < \infty  \}} Y_\tau = \b1_{\{\tau < \infty  \}} \xi$, \pas, which
 implies that $Y_\tau = \xi$ holds except on a $P-$null set $\cN_2$. Let $\cN \df \cN_1 \cup \cN_2$.  It then  holds on $N^c$ that
  \beas
  \int_\tau^T  f_\tau(s, Y_s, Z_s, U_s) ds -  \int_\tau^T  Z_s \, dB_s-\int_{ (\tau ,  T]} \int_\cX U_s(x)\tnp(ds,dx) =Y_\tau -\xi =0.
 \eeas
Therefore, one can deduce from \eqref{eqn-b621}    that on $\cN^c$
 \beas
   Y_{\tau \land t}&= & \xi + \int_{\tau \land t}^T  f_\tau(s, Y_s, Z_s, U_s) ds -  \int_{\tau \land t}^T  Z_s \, dB_s-\int_{ (\tau \land t,  T]} \int_\cX U_s(x)\tnp(ds,dx) \\
 &=&    \xi + \int_{\tau \land t}^\tau  f(s, Y_s, Z_s, U_s) ds -  \int_{\tau \land t}^\tau  Z_s \, dB_s-\int_{ (\tau \land t,  \tau]} \int_\cX U_s(x)\tnp(ds,dx)   , \q t \in [0,T],
  \eeas
   which shows that $\Big\{\big(Y_t(\o), Z_t(\o), U_t(\o)\big)\Big\}_{(t,\o) \in \[0, \tau\]}   $ is a solution of    \eqref{BSDEJ_random}.
   Moreover, since $(Y, Z, U) \in  \hS^p  $, we easily see that
   $ \Big\{\Big(Y_{\tau \land t}, \b1_{\{t \le \tau\}} Z_t, \b1_{\{t \le \tau\}} U_t\Big)\Big\}_{t \in [0,T]}$  belongs to $  \hS^p  $ as well.

    On the other hand, if $\Big\{\big(\wt{Y}_t(\o), \wt{Z}_t(\o), \wt{U}_t(\o)\big)\Big\}_{(t,\o) \in \[0, \tau\]}   $ is another solution of   \eqref{BSDEJ_random}
   such that  $ \Big\{\Big(\wt{Y}_{\tau \land t}, \b1_{\{t \le \tau\}} \wt{Z}_t, \\ \b1_{\{t \le \tau\}} \wt{U}_t\Big)\Big\}_{t \in [0,T]} \in  \hS^p  $, then it holds
 \pas~ that
   \beas
\wt{Y}_{\tau \land t}& \tn \dn =&  \n  \tn \xi  \+ \int_{\tau \land t}^\tau f(s, \wt{Y}_s, \wt{Z}_s, \wt{U}_s) ds     - \dn \int_{\tau \land t}^\tau \wt{Z}_s \, dB_s \n - \dn \int_{( \tau \land t,\tau]} \int_\cX \wt{U}_s(x)\tnp(ds,dx) \\
  & \tn  \dn =& \tn  \n \xi+ \n  \int_t^T \b1_{\{s \le \tau\}}   f(s, \wt{Y}_s, \wt{Z}_s, \wt{U}_s) ds -  \dn  \int_t^T \b1_{\{s \le \tau\}}  \wt{Z}_s \, dB_s  \n  -  \dn  \int_{(t, T]} \b1_{\{s \le \tau\}}  \int_\cX \wt{U}_s(x)\tnp(ds,dx) \\
 & \tn  \dn =& \tn  \n  \xi  \+ \int_t^T  \n f_\tau(s, \wt{Y}_{\tau \land s}, \b1_{\{s \le \tau\}} \wt{Z}_s, \b1_{\{s \le \tau\}} \wt{U}_s) ds
   - \dn \int_t^T  \n \b1_{\{s \le \tau\}}  \wt{Z}_s \, dB_s  \n - \dn  \int_{(t, T]} \n \int_\cX  \n \b1_{\{s \le \tau\}} \wt{U}_s(x)\tnp(ds,dx) , ~     t\in  [0,T],
   \eeas
which shows that $ \Big\{\Big(\wt{Y}_{\tau \land t}, \b1_{\{t \le \tau\}} \wt{Z}_t,  \b1_{\{t \le \tau\}} \wt{U}_t\Big)\Big\}_{t \in [0,T]}$ also solves BSDEJ\,$\big(\xi,f_\tau \big)$.
Hence, the uniqueness of the solution of BSDEJ\,$\big(\xi,f_\tau \big)$ in $   \hS^p $ entails that
 \beas
 0 & \tn  \n = & \tn \n   E\Big[ \, \underset{t \in [0,T] }{\sup} \big|Y_t  \-  \wt{Y}_{\tau \land t} \big|^p\Big]
+ E \bigg[   \Big( \n \int_0^T \n  \big| Z_t  \-  \b1_{\{t \le \tau\}} \wt{Z}_t \big|^2   \, dt \Big)^{\frac{p}{2}} \bigg]
+ E \bigg[ \Big( \n \int_0^T \n  \big\| U_t  \-  \b1_{\{t \le \tau\}} \wt{U}_t \big\|^2_{L^p_\nu}  \, dt \Big)^{\frac{p}{2}} \bigg] \\
 & \tn \n \ge& \tn \n   E\Big[ \, \underset{t \in [0,\tau] }{\sup} \big|Y_t  \- \wt{Y}_t \big|^p\Big]
+ E \bigg[ \Big( \n \int_0^\tau \n  \big| Z_t  \-    \wt{Z}_t \big|^2   \, dt \Big)^{\frac{p}{2}} \bigg]
+ E \bigg[ \Big( \n \int_0^\tau \n  \big\| U_t  \-    \wt{U}_t \big\|^2_{L^p_\nu}  \, dt \Big)^{\frac{p}{2}} \bigg]
\ge 0 ,
 \eeas
which implies that    $\Big\{\big(Y_t(\o), Z_t(\o), U_t(\o)\big)\Big\}_{(t,\o) \in \[0, \tau\]}   $ is the unique solution of   \eqref{BSDEJ_random} such that  $ \Big\{\Big(Y_{\tau \land t}, \b1_{\{t \le \tau\}} Z_t, \\ \b1_{\{t \le \tau\}} U_t\Big)\Big\}_{t \in [0,T]} \in  \hS^p  $.   \qed

 \fi

     \no {\bf Proof of Corollary \ref{cor_martingale}:}
 Clearly, $  f (t,\o,y,z,u) \df  0 $,
 $   \fa (t,\o, y,z,u)  \ins  [0,T]  \ti  \O  \ti  \hR  \ti  \hR^d
  \ti  L^2_\nu $
 satisfies (H1)$-$(H6).
 In light of Theorem \ref{thm_BSDEJ1}, BSDEJ$(\xi,0)$ admits a unique solution
 $(Y , Z , U) \ins  \hS^p $.
  Since \eqref{eq:f347} and Lemma \ref{lem_stoch_integr_Lp} show that
  $\int_0^t \n  Z_s dB_s   \+  \int_{(0,t]}   \int_\cX   U_s (x) \tnp (ds,dx)$, $t \ins [0,T]$ is a
 uniformly integrable  martingale, it holds for any    $t  \ins  [0,T]$ that
   \beas
    Y_t  \=   E \bigg[   \xi  \-  \int_t^T Z_s dB_s   \-  \int_{(t, T]}   \int_\cX   U_s (x) \tnp (ds,dx) \Big|\cF_t \bigg]
      \= E  [   \xi |\cF_t  ]  , \q  \pas ~
 \eeas
   In particular, one has $Y_0 =  E[\xi]$. It follows that for any  $t  \ins  [0,T]$
  \beas
   E[ \xi |\cF_t ] \=  Y_t  \=  Y_0  \+  \int_0^t \n  Z_s dB_s
   \+     \int_{(0,t]}  \n   \int_\cX  \n   U_s (x) \tnp (ds,dx)
   \=   E[\xi]  \+  \int_0^t \n  Z_s dB_s
    \+     \int_{(0,t]}  \n   \int_\cX  \n   U_s (x) \tnp (ds,dx) , \q \pas ,
  \eeas
  which together with the right continuity of processes $  E[ \xi |\cF_t ] $,
     $ \int_0^t Z_s dB_s  $ and
     $   \int_{(0,t]}  \n   \int_\cX  \n   U_s (x) \tnp (ds,dx) $, $ t  \ins  [0,T]  $
     leads  to  \eqref{eq:a235}.

   Next, let $(Z',U') \in  \hZ^{2,p} \ti \hU^p $ be another pair satisfying \eqref{eq:a235},
   so one has  that  \pas ~
 \beas
 \int_0^t    (Z_s \-  Z'_s) dB_s  +  \int_{(0,t]}   \int_\cX  \big( U_s (x)  \-  U'_s (x) \big) \tnp (ds,dx)
   =  0  , \q t  \ins  [0,T] .
 \eeas
 Clearly, the quadratic variation of the above process is
 $ \int_0^t    |Z_s \-  Z'_s|^2 ds  +  \int_{(0,t]}   \int_\cX  \big| U_s (x)  \-  U'_s (x) \big|^2 N_\fp (ds,dx)
   \=  0  $, $  t  \ins  [0,T] $,
  which implies that $Z_t (\o) \=  Z'_t (\o) $ for \dtp ~ $(t,\o) \ins [0,T] \ti \O$,
 and $ U(t,\o,x) \= U'(t,\o,x) $  for $dt \ti dP  \ti  \nu (dx)-$a.s.
  $(t,\o,x)  \ins  [0,T]  \ti  \O  \ti  \cX$.  \qed

  \no {\bf Proof of Theorem \ref{thm_BSDEJ2}: }
 By  the Burkholder-Davis-Gundy inequality, there exists a constant $\wt{\k} \> 0$ such that
 \beas
 E [ M_* ] \le \wt{\k} E \big\{  [M,M]^{1/2}_T \big\}
 \eeas
 for any 
 c\`adl\`ag local martingale $M $.
 Set constants $\wp_1  \df  \big( \frac{q}{4} \big)^{\frac{1}{q}}$,
 $\wp_2  \df   ( \nu(\cX)T  )^{\frac{1}{p}}  2^{\frac{2}{p}+1} \wt{\k}  p^{1-\frac{1}{p}} \wp^{-1}_1   $,
 $\wp_3 \df T \+ 1   -   \frac{p}{2}     \+  \wp^{-p}_2  \nu(\cX) T $,
 $ \wp_4 \df \frac{2}{p - 1} \big(16 \wt{\k}^2 p^2 \wp_3 \+  \frac{p}{2} \big)  $,
 $\wp_5 \df \frac{2}{p - 1} (1 \+ \wp_2)^{2-p} \big(  2^{2+p } \wt{\k}^p p^{p-1} \wp^{-p}_1   \wp_3 \+ 1 \big)  $
 and   $ \wp_6 \df  ( 4 p \wp_3  \+  \wp_4   \+ \wp_5 )^{-\frac{1}{p}}  $.
 Then we set $ \wt{C} \df
   \big\| \int_0^T \n  \big( \wt{\beta}^q_t \ve \L^2_t \big)  dt \big\|_{L^\infty_+ (\cF_T) } $ and
 define processes
        \beas
        a_t := \frac{2 \wt{\beta}^q_t  \+  \wt{C}^{\frac{q}{2}-1} \L^2_t
         }{q \wp^{p+q}_6 (\wp_4   \+ \wp_5)}
  \q \hb{and} \q           A_t  \df  p \n \int_0^t  a_s   ds  , \q
    t  \ins  [0,T] \, .
 \eeas
 Clearly, $ \dis C_A  \df  \|A_T\|_{L^\infty_+ (\cF_T) }
  \ls  \frac{(p \- 1)(4 p \wp_3  \+  \wp_4   \+ \wp_5)}{  \wp^q_6 (\wp_4   \+ \wp_5)}
  \big( 2 \wt{C} \+  \wt{C}^{\frac{q}{2}}   \big)   $.

   Let us introduce the following norm on $\hS^p$:
  \beas
  \q  \big\| (Y,Z,U) \big\|_\sharp \df \bigg\{ E \bigg[ \, T \underset{t \in [0,T]}{\sup} \big( e^{A_t } |Y_t|^p \big)
 \+   \Big( \int_0^T \n e^{ \frac{2}{p} A_t } | Z_t |^2   \, dt  \Big)^{\frac{p}{2}}
 \+ \int_0^T  \n \int_\cX \n e^{A_t }  |U_t (x) |^p \nu (dx) dt
  \bigg] \bigg\}^{\frac{1}{p}} , ~ \fa (Y,Z,U)  \ins  \hS^p .
 \eeas

   Fix $(Y,Z,U)  \ins  \hS^p$.
 The $\sP  \oti  \sB  ( \hR^l    )   \oti  \sB  ( \hR^{l \times d}  )
  \oti  \sB\big(L^p_\nu   \big)/\sB(\hR^l   ) -$measurability of generator $f$,
  the $\bF-$predictability of processes $(Y,Z)$ as well as the
  $\wh{\sP}    \oti    \cF_\cX -$measurability of the  random field  $ U $
  implies that the   process $ \ff_t \df f(t, Y_t, Z_t,U_t) $, $t \in [0,T] $ is $\bF-$progressively measurable.
    Also, \eqref{eq:f173}, Lemma \ref{lem_lp_esti} and H\"older's inequality show that
  \beas
   E \bigg[ \Big(\int_0^T |\ff_t| dt \Big)^p \bigg]
 & \tn \le  & \tn   E \bigg[ \Big( \int_0^T \n \big(\big|f(t, 0, 0,0)\big|  \+  \wt{\beta}_t |Y_t|
    \+  \L_t |Z_t|  \+  \wt{\beta}_t \|U_t\|_{L^p_\nu} \big)  dt \Big)^p \bigg] \nonumber \\
   & \tn    \le  & \tn   4^{p-1} E \bigg[     \Big( \n   \int_0^T  \n \big|f(t, 0, 0,0)\big| dt \Big)^p
    \dn + \n  \Big( \n  \int_0^T \wt{\beta}^q_t dt \Big)^{\frac{p}{q}}   \int_0^T \n \big( |Y_t|^p
    \+ \|U_t\|^p_{L^p_\nu} \big) dt
   \+  
   \Big(    \n \int_0^T \n  \L^2_t dt \Big)^{\frac{p}{2}}   \Big(    \n \int_0^T \n  |Z_t|^2 dt \Big)^{\frac{p}{2}}
      \bigg]  \nonumber  \\
    & \tn   \le  & \tn  4^{p-1} E \bigg[     \Big( \n   \int_0^T  \n \big|f(t, 0, 0,0)\big| dt \Big)^p \bigg]
   \+ 4^{p-1}  \Big( \wt{C}^{p-1} T \|Y\|^p_{\hD^p}
   \+    \wt{C}^{\frac{p}{2}}  \|Z\|^p_{\hZ^{2,p}}
   \+ \wt{C}^{p-1} \|U\|^p_{\hU^p}     \Big) \< \infty .
  \eeas
    Clearly, $\ff$ satisfies (H1)$-$(H6), so   Theorem \ref{thm_BSDEJ1} shows that
  the BSDEJ\,$(\xi,\ff)$   admits a unique solution $(\cY,\cZ,\cU )
    \ins     \hS^p    $.

  We  set $\Psi (Y,Z,U) \df   (\cY,\cZ,\cU  )$. To see that
 $\Psi$ defines a contraction map on $\hS^p $
 under  the  norm $\|\cd\|_\sharp$,  let $ \big( \wt{Y},\wt{Z},\wt{U}\big) $
 be another triplet in $   \hS^p $
 and let $ \big(\wt{\cY},\wt{\cZ},\wt{\cU} \big) $ be the unique solution to
 the BSDEJ\,$\big(\xi, \wt{\ff} \, \big)$ with
 $ \wt{\ff}_t \df f \big( t, \wt{Y}_t, \wt{Z}_t,\wt{U}_t \big) $, $t \in [0,T] $,
 so $\Psi \big( \wt{Y},\wt{Z},\wt{U}\big) \=  \big(\wt{\cY},\wt{\cZ},\wt{\cU} \big)$.
   For simplicity, we denote $(\sY \n , \sZ \n ,\sU) \df \big(\cY \- \wt{\cY},
 \cZ \- \wt{\cZ}, \cU \- \wt{\cU} \big)$.

 The structure of the following proof is similar to that of Lemma \ref{lem-a-priori},
 please make a comparison.

  \no  {\bf 1)}
 Let $ (t,\e)  \ins  [0 , T ]    \ti  (0,1]$.
  For the  function $ \vf_\e  $ defined in \eqref{eqn-b150},  applying It\^o's formula
  to process $e^{A_s} \vf^p_\e ( \sY_s )$ over the interval $[t,    T]$,
  we see from \eqref{Y_jump} that \pas
 \bea
   && \hspace{-1.2 cm}
   e^{A_t} \vf^p_\e (   \sY_t)
   \+  \frac{1}{2} \n \int_t^T \n e^{A_s}
   \hb{trace} \big( \sZ_s \sZ^T_s D^2   \vf^p_\e ( \sY_s)   \big)  ds
     \+  \underset{s \in (t,T] }{\sum}
   e^{A_s}\Big( \vf^p_\e ( \sY_s) \- \vf^p_\e ( \sY_{s-}) \-
      \big\lan D  \vf^p_\e ( \sY_{s-}) , \D \sY_s \big\ran  \Big)   \nonumber \\
  &&   = e^{A_T} \e^{\frac{p}{2}} \+ p \n  \int_t^T  \n   e^{A_s}   \big[   \vf^{p-2}_\e ( \sY_s)
  \big\lan \sY_s ,     \ff_s \- \wt{\ff}_s  \big\ran     \-  a_s \vf^p_\e ( \sY_s)\big]  ds
  \-  p \big( \sM_T \- \sM_t  \+   \fM_T  \-   \fM_t \big)  , \q \qq
  \label{eqn-b190c}
 \eea
   where  $   \sM_s   \df   \sM^\e_s \=  \int_0^s
    \n   e^{A_r}  \vf^{p-2}_\e ( \sY_{r-}) \big\lan \sY_{r-} ,  \sZ_r   dB_r \big\ran $
       and $  \fM_s   \df   \fM^\e_s \=     \int_{(0, s]} \n
       \int_\cX e^{A_r} \vf^{p-2}_\e ( \sY_{r-})   \big\lan   \sY_{r-} , \sU_r(x)  \big\ran    \tnp(dr,dx) $,
           $\fa  s \ins [0,T]$. 
  Similar to \eqref{eq:a114},
  we can deduce from   Taylor's Expansion Theorem and \eqref{eqn-b150}  that
  \beas
    \vf^p_\e \big( \sY_s  \big) \- \vf^p_\e \big( \sY_{s-}  \big) \- \big\lan D
 \vf^p_\e \big( \sY_{s-} \big), \D \sY_s   \big\ran
 \gs       p  (p \- 1)
 |\D \sY_s   |^2 \n \int_0^1 \n (1 \- \a) \vf^{p-2}_\e  ( \sY_{s-} \+ \a \D \sY_s  )  d \a .
 \eeas
   When $|\sY_{s-}| \ls \wp_2 |\D \sY_s|$, one has $ \vf^{p-2}_\e( \sY_{s-} \+ \a \D \sY_s )
 \gs \big( ( |\sY_{s-}|   +      \a  |\D \sY_s| )^2    +    \e   \big)^{\frac{p}{2}-1}
 \gs \big( (1   +   \wp_2)^2  |\D \sY_s|^2    +    \e   \big)^{\frac{p}{2}-1}
  \gs  (1   +   \wp_2)^{p-2}
   \big(    |\D \sY_s|^2    +    \e   \big)^{\frac{p}{2}-1}  $, $\fa \a \ins [0,1]$.
  So an analogy to \eqref{eq:a451} and \eqref{eq:b617} show that   \pas
    \beas
       &      &   \hspace{-1.5cm}   \underset{s \in (  t,T] }{\sum}
  e^{A_s } \Big( \vf^p_\e  ( \sY_s ) \- \vf^p_\e (  \sY_{s-} ) \- \big\lan D
 \vf^p_\e (  \sY_{s-} ), \D \sY_s   \big\ran\Big)  \nonumber   \\
 & & \ge    \frac12 (1 \+ \wp_2)^{p-2}  p  (p \- 1)   \n \int_{( t,T ]}
 \n \int_\cX \n   \b1_{\{ |\sY_{s-}|  \le \wp_2 | \sU_s(x)|        \}}
  e^{A_s}   \big| \sU_s (x)   \big|^2
  \big(   | \sU_s(x) |^2  \+  \e   \big)^{\frac{p}{2}-1} N_\fp (ds,dx)    .
       \eeas

  Plugging it and an analogy to \eqref{eqn-b198} back into \eqref{eqn-b190c},
  we can deduce from \eqref{eqn-b150}   that \pas
  \bea
   && \hspace{-2 cm}
   e^{A_t} \vf^p_\e (   \sY_t)
   \+  \frac{p(p \- 1)}{2}     \dn  \int_t^T \n  e^{A_s}   \vf^{p-2}_\e ( \sY_s) |\sZ_s|^2  ds
     \+  \frac{p(p \- 1)}{2}   (1 \+ \wp_2)^{p-2} \dn \int_{(t , T]}  \n \int_\cX \n   \b1_{\{ |\sY_{s-}|  \le \wp_2 |\sU_s(x)|  \}}
  e^{A_s}   \big|\sU_s (x)   \big|^2
  \big(   |\sU_s(x)|^2  \+  \e   \big)^{\frac{p}{2}-1} N_\fp (ds,dx)  \nonumber \\
  &&   \le   e^{C_A} \e^{\frac{p}{2}} \+ p \eta^\e_t
  \-  p \big( \sM_T \- \sM_t  \+   \fM_T  \-   \fM_t \big)  ,
  \label{eqn-b190d}
 \eea
 where $\eta^\e_t \df  \int_t^T  \n   e^{A_s}   \big[   \vf^{p-1}_\e ( \sY_s)
        | \ff_s \- \wt{\ff}_s  |    \-  a_s \vf^p_\e ( \sY_s)\big]  ds $.

  By Lemma \ref{lem_lp_esti}, the random variable
    $\vr_\e    \df  \underset{t \in [   0 ,  T  ]}{\sup}
   \big( e^{A_t} \vf^p_\e(  \sY_t ) \big) $  satisfies
 $   E  [   \vr_\e    ]    \ls
   e^{C_A} E \big[      \sY^p_*  \+   \e^{\frac{p}{2}} \big]
    \=  e^{C_A} \big(  \|\sY\|^p_{\hD^p}  \n +     \e^{\frac{p}{2}} \big)   \<  \infty  $.
    Since an analogy to \eqref{eq:a017}  shows that
    \bea \label{eq:f021}
    E \bigg[ \Big( \int_{(0, T]}  \n  \int_\cX
          |\sU_s (x) |^2     N_\fp (  ds, dx) \Big)^{\n \frac{p}{2}}   \bigg]
    \le E   \int_0^T  \n  \int_\cX         |\sU_s (x) |^p  \nu(dx) ds  ,
    \eea
 one can deduce from the Burkholder-Davis-Gundy inequality and Young's inequality    that
   \beas
   && \hspace{-1.5cm}   E   \bigg[ \, \underset{s \in [0, T]  }{\sup}|\sM_s|
    \+ \underset{s \in [0, T]}{\sup}|\fM_s|    \bigg]
    \ls  \wt{\k}   E \n \left[
  (\vr_\e)^{\frac{p-1}{p}   } \bigg( \int_0^T  \n
      e^{  \frac{2}{p} A_s } |\sZ_s|^2    ds \bigg)^{\n \frac12}
   \+  (\vr_\e)^{\frac{p-1}{p}  } \bigg( \int_{(0, T]}  \n  \int_\cX
  \n   e^{  \frac{2}{p} A_s }      |\sU_s (x) |^2     N_\fp (  ds, dx) \bigg)^{\n \frac12}   \right]  \nonumber \\
 &&   \ls  \frac{\wt{\k}}{p}
    E   \bigg[  2(p \- 1)  \vr_\e  \+  e^{C_A} \bigg( \int_0^T  \n
       |\sZ_s|^2    ds \bigg)^{\n \frac{p}{2}}    \+   e^{C_A}  \n  \int_0^T    \n  \int_\cX
         |\sU_s (x) |^p  \nu(dx) ds      \bigg]   \<  \infty .  
 \eeas
 So both $ \sM $ and    $ \fM  $ are   uniformly integrable    martingales.

  Taking $t \= 0$ and  taking expectation in \eqref{eqn-b190d}  yields that
 \bea
  &&  \hspace{-2cm}   E \n  \int_0^T \n e^{A_s}   \vf^{p-2}_\e ( \sY_s) |\sZ_s|^2  ds
     \+  (1 \+ \wp_2)^{p-2}      E  \n \int_0^T  \n \int_\cX \n   \b1_{\{ |\sY_{s-}|  \le \wp_2 |\sU_s(x)|  \}}
   e^{A_s}   \big|\sU_s (x)   \big|^2
  \big(   |\sU_s(x)|^2  \+  \e   \big)^{\frac{p}{2}-1} \nu(dx) ds \nonumber \\
 &&   \le \n \frac{ 2  }{p (p \- 1)} \big( e^{C_A} \e^{\frac{p}{2}} \+ p E[\eta^\e_0] \big) . \q \label{eq:f024}
 \eea

  \no {\bf 2)}  Clearly, $ \lmtu{\e \to 0} |\sU (s,\o,x)|^2 \big(   |\sU (s,\o,x)|^2   +   \e   \big)^{\frac{p}{2}-1}
 \= |\sU (s,\o,x)|^p $, $\fa (s,\o,x)  \ins  [0,T]  \ti  \O  \ti  \cX$,
 so the monotone convergence theorem  implies  that
 \beas
  \lmtu{\e \to 0}  E \n \int_0^T \n \int_\cX \n   \b1_{\{|\sY_{s-}| \le \wp_2 |\sU_s(x)|      \}}
  e^{A_s}    |\sU_s (x)   |^2
  \big(   |\sU_s(x)|^2   \+   \e   \big)^{\frac{p}{2}-1} \nu (dx) ds
  \=  E \n \int_0^T  \n \int_\cX \n   \b1_{\{|\sY_{s-}| \le \wp_2 |\sU_s(x)|      \}}
  e^{A_s}    |\sU_s (x)   |^p     \nu (dx) ds  .
 \eeas
 On the other hand, one has    $ \lmt{\e \to 0} \big( \vf^{p-1}_\e ( \sY_s)
  | \ff_s \- \wt{\ff}_s  |    \-  a_s \vf^p_\e ( \sY_s) \big)
   \=    |\sY_s|^{p-1}   | \ff_s \- \wt{\ff}_s  |    \-  a_s |\sY_s|^p $, $ s  \ins  [0,T]$
   and   $|\eta^\e_t| \ls
   \wt{\eta} \df    \int_0^T  \n  e^{A_s}  \big[   \vf^{p-1}_1 ( \sY_s)
        | \ff_s \- \wt{\ff}_s  |    \+  a_s \vf^p_1 ( \sY_s)\big]  ds $, $\fa \e \ins (0,1]$.
   Since \eqref{eq:f173}, Young's inequality and H\"older's inequality imply that
        \beas
        \hspace{-3mm}
     E \big[ \wt{\eta} \big] & \tn  \tn \le  & \tn  \tn   E \n \int_0^T  \n   e^{A_s}     \Big[   \vf^{p-1}_1 ( \sY_s)
      \big( \wt{\beta}_s \big|Y_s  \- \wt{Y}_s\big|   \+  \L_s \big|Z_s  \- \wt{Z}_s\big|
     \+   \wt{\beta}_s \big\|U_s  \- \wt{U}_s\big\|_{L^p_\nu} \big)     \+  a_s \vf^p_1 ( \sY_s) \Big]  ds \\
        & \tn  \tn  \le  & \tn  \tn   E \n \int_0^T  \n   e^{A_s}    \bigg[
       \frac{1}{p} \big|Y_s  \- \wt{Y}_s\big|^p \+ \frac{1}{p} \big\|U_s  \- \wt{U}_s\big\|^p_{L^p_\nu}
        \+ \Big(\frac{2}{q} \wt{\beta}^q_s
        \+ a_s \Big) \vf^p_1 ( \sY_s) \bigg]   ds
        \+    E \Bigg[ \vr_1^{1/q } \bigg(\n \int_0^T \n   \L^2_s ds \bigg)^{\frac{1}{2}}
        \bigg( \n  \int_0^T \n e^{\frac{2}{p} A_s} \big|Z_s  \- \wt{Z}_s\big|^2 ds \bigg)^{\frac12} \Bigg] \qq \\
         & \tn \tn  \le & \tn  \tn  \frac{1}{p} \big\| \big(Y \- \wt{Y},Z \- \wt{Z},U \- \wt{U}\big) \big\|^p_\sharp         \+  \Big(\frac{2}{q} \wt{C}
         \+  \frac{1}{p} C_A  \+  \frac{1}{q}   \wt{C}^{\frac{q}{2}}   \Big)
        E[\vr_1]    \< \infty ,
        \eeas
  an application of the dominated convergence theorem shows that
  \bea \label{eq:f041}
  \lmt{\e \to 0} E[\eta^\e_0] \= E[\eta ]
  \eea
    with   $\eta  \df \int_0^T  \n   e^{A_s}   \big[   |\sY_s|^{p-1}   | \ff_s \- \wt{\ff}_s  |
    \-  a_s |\sY_s|^p \big]  ds$. Letting $\e \n \to \n  0$ in \eqref{eq:f024} then yields that
    \bea \label{eq:f039}
   (1 \+ \wp_2)^{p-2}    E  \n \int_0^T  \n \int_\cX \n   \b1_{\{ |\sY_{s-}| \le \wp_2 |\sU_s(x)|  \}}
  e^{A_s}   \big|\sU_s (x)   \big|^p \nu(dx) ds   \ls  \frac{2 }{p \- 1 } E[\eta ] .
    \eea

    Now, fix $\e \ins (0, 1]$ again.  We also see from \eqref{eqn-b190d} that
  \bea \label{eq:f033}
  E \big[  \vr_\e  \big]  \ls e^{C_A} \e^{\frac{p}{2}} \+ p E \n \int_0^T  \n   e^{A_s}     \vf^{p-1}_\e ( \sY_s)
        | \ff_s \- \wt{\ff}_s  |   ds
        \+ 2p E \bigg[ \, \underset{s \in [0,T]}{\sup} |\sM_s| \+ \underset{s \in [0,T]}{\sup} |\fM_s|  \bigg] .
  \eea
   The Burkholder-Davis-Gundy inequality, Young's inequality,
 \eqref{Y_jump} as well as  an analogy to \eqref{eq:f021} implies that
 \beas
    && \hspace{-1.5 cm}  2p E\bigg[ \, \underset{s \in [0,T]}{\sup}| \sM_s|
   \+ \underset{s \in [0,T]}{\sup}|\fM_s |  \bigg]
 \ls   2 \wt{\k}  p      E \Bigg[ \vr_\e^{1/2  }
     \bigg(  \n \int_0^T \n e^{A_s}
   \vf^{p-2}_\e ( \sY_{s-})  |\sZ_s|^2 ds  \bigg)^{\n \frac12} \+
    (\vr_\e)^{\frac{p-1}{p}}   \bigg(  \n  \int_{(0 , T]} \int_\cX
   \n  e^{  \frac{2}{p} A_s }  |\sU_s(x)|^2 N_\fp(ds, dx) \bigg)^{\n \frac12}  \Bigg] \nonumber \\
&\tn \le & \tn \dn \Big( \frac14 \+ \frac{1}{q} \wp^q_1 \Big) E   \big[  \vr_\e    \big]
  \+  4 \wt{\k}^2 p^2      E \n \int_0^T   \n e^{A_s}  \vf^{p-2}_\e ( \sY_s )  |\sZ_s|^2 ds
     \+  2^p \wt{\k}^p p^{p-1} \wp^{-p}_1
      E \n   \int_0^T   \n  \int_\cX \n e^{A_s}     |\sU_s(x)|^p \nu(dx) ds   .
 \eeas
 One can also deduce that
  \bea
 && \hspace{-1cm} E \n   \int_0^T  \n \int_\cX \n e^{A_s }  | \sU_s (x) |^p \nu (dx) ds
    \le  E \n    \int_0^T  \n \int_\cX \n e^{A_s }
 \Big( \b1_{\{   |\sY_{s-} | > \wp_2 |\sU_s (x)  |    \}}  \wp^{-p}_2  \vf^p_\e ( \sY_{s-} )
  \+   \b1_{\{     |\sY_{s-} | \le  \wp_2  |\sU_s (x)  |    \}}   | \sU_t (x) |^p
  \Big) \nu (dx) ds    \nonumber \\
&& \hspace{-4mm} \le   \wp^{-p}_2  \nu(\cX) T   E[\vr_\e]
      \+  E \n \int_0^T  \n \int_\cX \n
   \b1_{\{ |\sY_{s-} | \le  \wp_2  |\sU_s (x)  |    \}}  e^{A_s }   | \sU_t (x) |^p    \nu (dx) ds . \label{eq:f035}
 \eea

 \no {\bf 3)}  As $ \frac{1}{q} \wp^q_1 \=   2^p \wt{\k}^p p^{p-1} \wp^{-p}_1   \wp^{-p}_2  \nu(\cX) T  \= \frac14 $,
 putting these two inequalities into \eqref{eq:f033} yields    that
 \bea
  E \big[  \vr_\e   \big]
   & \tn  \le  & \tn 4 e^{C_A} \e^{\frac{p}{2}} \+ 4 p E \n \int_0^T  \n   e^{A_s}     \vf^{p-1}_\e ( \sY_s)
        | \ff_s \- \wt{\ff}_s  |   ds
  \+  16 \wt{\k}^2 p^2      E \n \int_0^T   \n e^{A_s}  \vf^{p-2}_\e ( \sY_s )  |\sZ_s|^2 ds \nonumber \\
   & \tn  & \tn   +    2^{2+p } \wt{\k}^p p^{p-1} \wp^{-p}_1
      E \n   \int_0^T   \n  \int_\cX \n  \b1_{\{ |\sY_{s-} | \le  \wp_2  |\sU_s (x)  |    \}}
       e^{A_s}     |\sU_s(x)|^p \nu(dx) ds    .   \label{eq:f037}
 \eea
 Similar to \eqref{eqn-b338},  Young's inequality shows that
 \beas
  E \bigg[   \Big( \int_0^T \n e^{ \frac{2}{p} A_s } | \sZ_s |^2   \, ds  \Big)^{\frac{p}{2}}   \bigg]
  \ls  E \bigg[   \Big( \big( \vr_\e \big)^{\frac{2-p}{p}} \dn
  \int_0^T \n e^{   A_s } \vf^{p-2}_\e ( \sY_s ) | \sZ_s |^2   \, ds  \Big)^{\frac{p}{2}}    \bigg]
   \ls  \frac{2 \- p}{2}   E[\vr_\e]
   \+  \frac{p}{2} E \n \int_0^T \n e^{   A_s } \vf^{p-2}_\e ( \sY_s ) | \sZ_s |^2   \, ds .
 \eeas
 Then we see from \eqref{eq:f035}, \eqref{eq:f037}, \eqref{eq:f024} and \eqref{eq:f039}  that
 \beas
 && \hspace{-8mm} \|(\sY,\sZ,\sU)\|^p_\sharp \ls
 E \bigg[ T \vr_\e  \+  \Big( \int_0^T \n e^{ \frac{2}{p} A_s } | \sZ_s |^2   \, ds  \Big)^{\frac{p}{2}}
 \+ \int_0^T  \n \int_\cX \n e^{A_s }  | \sU_s (x) |^p \nu (dx) ds
  \bigg] \\
&&   \le \big( T \+ 1 \- \frac{p}{2}     \+  \wp^{-p}_2  \nu(\cX) T \big) E[\vr_\e]
   \+  \frac{p}{2} E \n \int_0^T \n e^{   A_s } \vf^{p-2}_\e ( \sY_s ) | \sZ_t |^2   \, ds
      \+  E \n \int_0^T  \n \int_\cX \n
   \b1_{\{ |\sY_{s-} | \le \wp_2  |\sU_s (x)  |    \}}  e^{A_s }   | \sU_t (x) |^p    \nu (dx) ds \\
&&   \le \Big(4  \wp_3 \+ \frac{1}{p} \wp_4 \Big) e^{C_A} \e^{\frac{p}{2}} \+
 4 p \wp_3  E \n \int_0^T  \n   e^{A_s}     \vf^{p-1}_\e ( \sY_s)
        | \ff_s \- \wt{\ff}_s  |   ds    \+ \wp_4 E[\eta^\e_0]        \+ \wp_5 E[\eta] \\
&&        =  \Big(4  \wp_3 \+ \frac{1}{p} \wp_4 \Big) e^{C_A} \e^{\frac{p}{2}}
 \+  4 p \wp_3 E \n \int_0^T  \n  a_s e^{A_s}     \vf^p_\e ( \sY_s)  ds
 \+ ( 4 p \wp_3  \+  \wp_4 ) E[\eta^\e_0]        \+ \wp_5 E[\eta]    .
 \eeas
 Since $E \n \int_0^T  \n  a_s e^{A_s}     \vf^p_1 ( \sY_s)  ds
 \ls  \frac{1}{p} E [ \vr_1 A_T ]  \ls  \frac{C_A}{p} E [ \vr_1 ]  \<  \infty $,
 letting $\e \n \to \n 0$, we can deduce from the dominated convergence theorem and \eqref{eq:f041} that
 \bea
 \|(\sY,\sZ,\sU)\|^p_\sharp & \tn \le  & \tn   4 p \wp_3 E \n \int_0^T  \n  a_s e^{A_s}     | \sY_s |^p  ds
 \+ ( 4 p \wp_3  \+  \wp_4   \+ \wp_5 ) E[\eta] \nonumber \\
  & \tn = & \tn  ( 4 p \wp_3  \+  \wp_4   \+ \wp_5 )
 E \n \int_0^T  \n   e^{A_s}   |\sY_s|^{p-1}   | \ff_s \- \wt{\ff}_s  | ds
 \- (  \wp_4   \+ \wp_5 ) E \n \int_0^T  \n  a_s e^{A_s}     | \sY_s |^p  ds .   \label{eq:f043}
 \eea

  By \eqref{eq:f173}, Young's inequality and H\"older's inequality again,
        \beas
        \hspace{-5mm}
     E \n \int_0^T  \n   e^{A_s}   |\sY_s|^{p-1}   | \ff_s \- \wt{\ff}_s  | ds
      & \tn  \tn \le  & \tn  \tn   E \n \int_0^T  \n   e^{A_s}   |\sY_s|^{p-1}
      \big( \wt{\beta}_s \big|Y_s  \- \wt{Y}_s\big|   \+ \L_s  \big|Z_s  \- \wt{Z}_s\big|
     \+  \wt{\beta}_s \big\|U_s  \- \wt{U}_s\big\|_{L^p_\nu} \big)      ds \\
        & \tn  \tn  \le  & \tn  \tn   E \n \int_0^T  \n   e^{A_s}    \bigg[
       \frac{\wp_6^p}{p} \big|Y_s  \- \wt{Y}_s\big|^p
       \+ \frac{\wp_6^p}{p} \big\|U_s  \- \wt{U}_s\big\|^p_{L^p_\nu}
        \+ 2  \frac{\wp_6^{-q}}{q} \wt{\beta}^q_s   |\sY_s|^p \bigg]   ds \\
      &&  \n  +        E \Bigg[ \bigg(   \int_0^T \n  \L^2_s   ds \bigg)^{1-\frac{1}{q}-\frac12}
     \bigg(   \int_0^T \n  \L^2_s e^{A_s}     | \sY_s |^p ds \bigg)^{\frac{1}{q}}
        \bigg(   \int_0^T \n e^{\frac{2}{p} A_s} \big|Z_s  \- \wt{Z}_s\big|^2 ds \bigg)^{\frac12} \Bigg] \qq \\
         & \tn \tn  \le & \tn  \tn  \frac{1}{p} \wp_6^p
         \big\| \big(Y \- \wt{Y},Z \- \wt{Z},U \- \wt{U}\big) \big\|^p_\sharp
        \+ \frac{\wp_6^{-q}}{q}  E \n \int_0^T  \n   e^{A_s}
           \big(2 \wt{\beta}^q_s  \+  \wt{C}^{\frac{q}{2}-1} \L^2_s \big)   |\sY_s|^p   ds   .
        \eeas
 Plugging this inequality back into \eqref{eq:f043} yields that
 $ \big\| \big(\cY \- \wt{\cY},\cZ \- \wt{\cZ},\cU \- \wt{\cU}\big) \big\|^p_\sharp
   \le  \frac{1}{p}  \big\| \big(Y \- \wt{Y},Z \- \wt{Z},U \- \wt{U}\big) \big\|^p_\sharp $.
  Therefore, $\Psi$ is a contraction mapping on $\hS^p $
 under  the  norm $\|\cd\|_\sharp$. Then the unique fixed point
 $(Y,Z,U) \ins \hS^p$ of $\Psi$
   forms a unique solution of BSDEJ\,$(\xi,f)$ in $\hS^p$. \qed

\no {\bf Proof of Corollary \ref{cor_f_tau}:}
 Let $\tau \ins \cT$ and $\xi \ins L^p (\cF_\tau)$. Since
   the $p-$generator $f_\tau$ also satisfies  either  \(H1\)$-$\(H6\) or \eqref{eq:f173},
   Theorem \ref{thm_BSDEJ1} and Theorem \ref{thm_BSDEJ2}
  show  that BSDEJ\;$(\xi,f_\tau)$ admits a unique solution
   $    (Y,Z,U )   \ins  \hS^p $.

By  \eqref{eq:f347} and  \eqref{def_Poisson_integr},
$M_t \df \int_0^t \n Z_s dB_s  \+ \int_{(0,t]} \n \int_\cX \n U_s(x) \tnp (ds,dx) $ is a uniformly integrable martingale.
  Since  $Y_\tau \= \xi \+ \int_\tau^T \n f_\tau (s,Y_s,Z_s,U_s) ds \- M_T \+ M_\tau
 \= \xi \- M_T \+ M_\tau   $, \pas ~
 taking conditional expectation $E[~|\cF_\tau]$ shows that $ Y_\tau \= \xi $, \pas ~
 Then   processes $\Big\{(\cY_t,\cZ_t,\cU_t) \df \big( Y_{\tau \land t}, \b1_{\{t \le \tau\}} Z_t ,
 \b1_{\{t \le \tau\}} U_t \big) \Big\}_{t \in [0,T]} \ins \hS^p$   satisfy that \pas
  \beas
  \cY_t & \tn =& \tn  Y_{\tau \land t} = Y_\tau + \int_{\tau \land t}^\tau \n   f_\tau (s,  Y_s,  Z_s,  U_s)   ds
  \- \int_{\tau \land t}^\tau \n Z_s dB_s
  \- \int_{(\tau \land t,\tau]} \n \int_\cX \n  U_s (x) \tnp (ds,dx) \\
  & \tn =& \tn  \xi + \int_t^T   \n   f_\tau  \big(s,  \cY_s, \cZ_s,   \cU_s \big)   ds
  \- \int_t^T \n \cZ_s dB_s
  \- \int_{(t,T]} \n  \int_\cX \n \cU_s (x) \tnp (ds,dx) , \q t \in [0,T] ,
 \eeas
 so $(\cY ,\cZ ,\cU)$ also solves BSDEJ\;$(\xi,f_\tau)$. It follows that
 $P\{Y_t \= \cY_t \= Y_{\tau \land t}, \; t \ins [0,T] \} \= 1$ and that
 $ (Z_t,U_t) \= (\cZ_t, \cU_t) \= \b1_{\{t \le \tau\}} (Z_t,U_t) $, \dtp \qed

\no {\bf Proof of Proposition \ref{gen_repre}: } 
Set   $ C_{\wt{\beta}} \df
   \big\| \int_0^T \n    \wt{\beta}^q_t    dt \big\|_{L^\infty_+ (\cF_T) } $
   and let $\cC$ denote a generic constant    depending on $T$, $\nu(\cX)$, $p$,   $ \fc $  and $C_{\wt{\beta}} $,
   whose form may vary from line to line.
 Fix   $ (  t,y,z,u)  \ins    [0,T) \ti  \hR^l \ti  \hR^{l \times d} \ti L^p_\nu$
 such that \eqref{eq:x231} holds.

 \no {\bf 1)} Let $s \ins (t,T]$.    For any $A \ins \cF_t$,
  the Burkholder-Davis-Gundy inequality, \eqref{eqn-d011} and an analogy to \eqref{eq:a017} imply that
 \beas
  E \bigg[ \b1_A \underset{t' \in  (t,s]}{\sup} \big| V (t,t',z,u) \big|^p \bigg]
     & \tn \dn  \ls   & \tn \dn  2^{p-1}  E \bigg[   \underset{t' \in  (t,s]}{\sup}
     \Big| \int_t^{t'} \b1_A  z d B_s \Big|^p \+   \underset{t' \in  (t,s]}{\sup}
     \Big| \int_{s \in (t,t']} \n \int_\cX \n \b1_A  u(x) \tnp (ds,dx) \Big|^p  \bigg]     \nonumber \\
           & \tn \dn  \ls   & \tn \dn    c_{p,l}     E   \bigg[    \Big( \int_t^s  \n
        \b1_A |z|^2     dr  \Big)^{  \frac{p}{2}}
  \+  \int_t^s  \n  \int_\cX \n  \b1_A |u (x)|^p  \nu(dx) dr    \bigg]
   \ls  c_{p,l}  \Big( (s \- t)^{\frac{p}{2}}|z|^p \+ (s\-t) \|u\|^p_{L^p_\nu} \Big) P(A)      .
 \eeas
 As $A $ vary over $\cF_t$, we obtain
 \bea
  E \bigg[   \underset{t' \in  (t,s]}{\sup} \big| V (t,t',z,u) \big|^p \Big| \cF_t \bigg]
  \le c_{p,l} \Big( (s \- t)^{\frac{p}{2}}|z|^p \+ (s\-t) \|u\|^p_{L^p_\nu} \Big) .
 \label{eq:x201}
 \eea
 In particular, taking expectation implies that $  V (t,s,z,u) \ins L^p(\cF_s)  $.
 Since the $p-$generator $f_s$   also  satisfies   \eqref{eq:f173c},
      Theorem \ref{thm_BSDEJ2} and Corollary \ref{cor_f_tau}
    show   that the BSDEJ\,$ \big(y \+ V (t,s,z,u), f_s  \big)$
  admits a unique solution $ (Y^{s,y,z,u},Z^{s,y,z,u},U^{s,y,z,u})  \ins \hS^p $
  satisfying that
   $P \big\{Y^{s,y,z,u}_r \= Y^{s,y,z,u}_s, \, r \ins [s,T] \big\} \= 1$ and   that
  $\big(Z^{s,y,z,u}_r,U^{s,y,z,u}_r\big) \= 0$, $dr \ti dP-$a.s.   on $(s,T] \ti \O$.

  Fix $\e \ins \big( 0 , \d  \land 1 \big]  $.
  We simply denote $ (Y^{t+\e,y,z,u},Z^{t+\e,y,z,u},U^{t+\e,y,z,u})$
  by $(Y^\e,Z^\e,U^\e)$.
 For any $    (s,\o, y',z',u')\in [0,T]\ti  \O \ti  \hR^l \ti  \hR^{l \times d} \ti L^p_\nu$,  define
 \beas
  \ff^\e(s, \o, y',z',u') \df   \b1_{\{s <  t+\e \}}
  f \big(s, \o, y'\+ y \+ \b1_{\{s > t\}}    V (t,s,z,u) (\o)  ,z'\+\b1_{\{s > t\}} z,u'\+ \b1_{\{s > t\}}  u\big) ,
  \eeas
  which    satisfies     \eqref{eq:f173c} again.
  Since $ \big\{ \b1_{\{s > t\}} V (t,s,z,u) \big\}_{s \in  [0,T]}$ is an $\bF-$progressively measurable process,
  the $\sP     \oti  \sB  ( \hR^l )   \oti  \sB  ( \hR^{l \times d}  )   \oti  \sB\big(L^p_\nu   \big)/\sB(\hR^l )
  -$measurability   of   $  f    $   implies that of $\ff^\e$.
     H\"older's inequality and \eqref{eq:x201} show that
  \beas
 && \hspace{-1cm} E \bigg[ \Big(\int_0^T \n | \ff^\e(s,0,0,0) |  ds \Big)^p \bigg]
  \= E \bigg[ \Big(\int_0^{t+\e} \n |f(s,y\+ \b1_{\{s > t\}}    V (t,s,z,u)   ,
  \b1_{\{s > t\}} z,  \b1_{\{s > t\}}  u) |  ds \Big)^p \bigg] \\
 &&    \ls    E \bigg[ \Big(\int_0^{t+\e} \n   |f(s,0  , 0,  0) | ds
   \+  \int_0^{t+\e} \n  \wt{\beta}_s  ds  \Big(  |y| \+ \underset{s \in (t,t+\e]}{\sup} |V(t,s,z,u)| \Big)
   \+ \fc \int_t^{t+\e} \n \big(|z|\+\|u\|_{L^p_\nu}\big) ds    \Big)^p \bigg] \\
 &&    \ls 4^{p-1} E \bigg[ \Big(\int_0^{t+\e} \n   |f(s,0  , 0,  0) | ds \Big)^p
 \+ (t\+\e) \Big(\int_0^{t+\e} \n  \wt{\beta}^q_s   ds \Big)^{p-1}
   \Big(  |y|^p \+ \underset{s \in (t,t+\e]}{\sup} |V(t,s,z,u)|^p \Big)
     \bigg] \+ 4^{p-1} \fc^p \e^p  \big(|z|\+\|u\|_{L^p_\nu}\big)^p \\
 &&    \ls 4^{p-1} E \bigg[ \Big(\int_0^T \n   |f(s,0  , 0,  0) | ds \Big)^p     \bigg]
 \+ 4^{p-1} T  C_{\wt{\beta}}^{p-1}   \Big(  |y|^p \+ c_{p,l}  \big( \e^{\frac{p}{2}}|z|^p \+ \e \|u\|^p_{L^p_\nu} \big) \Big)
      \+ 4^{p-1} \fc^p \e^p  \big(|z|\+\|u\|_{L^p_\nu}\big)^p      \<\infty .
  \eeas
  By \eqref{eq:x201} again, one can deduce that
  \beas
  \big(\cY^\e_s, \cZ^\e_s , \cU^\e_s \big) \df
  \big( Y^\e_s \- y \- \b1_{\{ s >  t \}} V \big(t,s \ld (t \+ \e) ,z,u \big) ,
    Z^\e_s \- \b1_{\{ t < s \le t+\e \}} z  ,   U^\e_s \- \b1_{\{ t < s \le t+\e \}} u \big) ,
    \q \fa s \in [0,T]
 \eeas
 belongs to $\hS^p$ and satisfies BSDEJ\,$(0,\ff^\e)$.

 \no {\bf 2)} Next, let $ A \ins  \cF_t$.
  Applying Corollary \ref{cor_martingale} with $\xi \=    \b1_A \cY^\e_t \ins L^p(\cF_t)$
 shows that there exists   a unique pair
 $ \big(\sZ^\e, \sU^\e\big) \ins    \hZ^{2,p}      \ti  \hU^p $  such that \pas
 \beas
 \sY^\e_s \df  E \big[\b1_A \cY^\e_t|\cF_s\big]   \=    E \big[ \b1_A \cY^\e_t \big]    \+    \int_0^s    \sZ^\e_r dB_r
   \+    \int_{(0,s]}   \int_\cX   \sU^\e_r (x) \tnp (dr,dx) , \q s  \ins  [0,T] .
 \eeas
   Define    $\ol{Y}^\e_s \df \b1_{\{s<t  \}} \sY^\e_s  \+ \b1_{\{s \ge t  \}} \b1_A \cY^\e_s $,
    $ \big( \ol{Z}^\e_s,\ol{U}^\e_s \big) \df \b1_{\{s \le t  \}} \big( \sZ^\e_s, \sU^\e_s \big)
 \+ \b1_{\{s > t  \}} \b1_A (\cZ^\e_s, \cU^\e_s)$, $\fa s \ins [0,T]$.
 One can deduce from Doob's martingale inequality that  $\big(\ol{Y}^\e,\ol{Z}^\e,\ol{U}^\e\big) $ belong to $  \hS^p $.
  Since $\big\{ \b1_A \b1_{  \{s \ge t     \} } \big\}_{s \in [0,T]}$ is an $\bF-$adapted c\`adl\`ag process,
   the mapping
   \beas
     \ff^\e_A (s, \o, y',z',u') \df  \b1_A \b1_{\{s \ge  t    \}} \ff^\e (s, \o, y',z',u') ,
 \q \fa    (s,\o, y',z',  u) \in [0,T]\ti  \O \ti  \hR^l \ti  \hR^{l \times d} \ti L^p_\nu
 \eeas
  is also a $p-$generator that satisfies \eqref{eq:f173c} and $\int_0^T \n |\ff^\e_A (t,0,0,0)| dt \ins L^p_+ (\cF_T)$.

 For any $s \ins [t,T]$, multiplying $\b1_A$ to the BSDEJ\,$(0,\ff^\e)$ over period $[s, T]$ yields that
\beas
 \b1_A \cY^\e_s & \tn = & \tn   \int_s^T \n \b1_A
 \ff^\e  (r, \cY^\e_r,\cZ^\e_r, \cU^\e_r)dr
 \-   \int_s^T \n \b1_A \cZ^\e_r dB_r  \-   \int_{(s,T]} \n \int_\cX \n  \b1_A  \cU^\e_r (x) \tnp (dr,dx) \nonumber \\
 & \tn = & \tn   \int_s^T \n \b1_A \b1_{\{r \ge t \}} \ff^\e   (r, \b1_A \cY^\e_r, \b1_A \cZ^\e_r, \b1_A \cU^\e_r)dr
 \-   \int_s^T \n \ol{Z}^\e_r dB_r  \-   \int_{(s,T]} \n \int_\cX \n  \ol{U}^\e_r (x) \tnp (dr,dx) \nonumber \\
  & \tn = & \tn    \int_s^T \n \ff^\e_A  \big(r, \ol{Y}^\e_r  , \ol{Z}^\e_r, \ol{U}^\e_r \big) dr
 \-   \int_s^T \n \ol{Z}^\e_r dB_r  \-   \int_{(s,T]} \n \int_\cX \n    \ol{U}^\e_r (x) \tnp (dr,dx) ,
 \q \pas
\eeas
    The right continuity of process $\cY^\e$  shows that \pas
\bea \label{eq:x361}
 \b1_A \cY^\e_s  \=     \int_s^T \n \ff^\e_A   \big(r, \ol{Y}^\e_r  , \ol{Z}^\e_r, \ol{U}^\e_r \big) dr
 \-   \int_s^T \n \ol{Z}^\e_r dB_r  \-   \int_{(s,T]} \n \int_\cX \n    \ol{U}^\e_r (x) \tnp (dr,dx) ,
 \q s \in [t,T] .
\eea
 On the other hand,
 for any $s \ins [0,t)$, since $ \sY^\e_t \= E \big[ \b1_A \cY^\e_t |\cF_t \big] \= \b1_A \cY^\e_t $, \pas ~
   taking $s \= t $ in \eqref{eq:x361} yields that
 \beas
 \sY^\e_s & \tn \=  & \tn  \sY^\e_t    \- \int_s^t \n   \sZ^\e_r dB_r
 \- \int_{(s,t]} \n \int_\cX \n   \sU^\e_r (x) \tnp (dr,dx)
   \=   \b1_A \cY^\e_t  \- \int_s^t \n    \ol{Z}^\e_r dB_r
\- \int_{(s,t]} \n \int_\cX \n    \ol{U}^\e_r (x) \tnp (dr,dx)   \\
 & \tn \=  & \tn   \int_s^T  \ff^\e_A (r,\ol{Y}^\e_r  , \ol{Z}^\e_r, \ol{U}^\e_r ) dr
\- \int_s^T \n    \ol{Z}^\e_r dB_r
\- \int_{(s,T]} \n \int_\cX \n    \ol{U}^\e_r (x) \tnp (dr,dx) , \q \pas
\eeas

This equality together with   \eqref{eq:x361} and the right continuity of   $\ol{Y}^\e$ shows  that
\pas
 \beas
 \ol{Y}^\e_s = \b1_{\{s<t  \}} \sY^\e_s   \+ \b1_{\{s \ge t  \}} \b1_A \cY^\e_s
 =  \int_s^T   \ff^\e_A  (s,\ol{Y}^\e_r  , \ol{Z}^\e_r, \ol{U}^\e_r ) dr
\- \int_s^T \n    \ol{Z}^\e_r dB_r
\- \int_{(s,T]} \n \int_\cX \n    \ol{U}^\e_r (x) \tnp (dr,dx)  , \q s \ins [0,T] .
 \eeas
 Hence,    $\big(\ol{Y}^\e,\ol{Z}^\e,\ol{U}^\e\big)$   solve  BSDEJ\,$( 0, \ff^\e_A )$.

 By Proposition \ref{prop-a-priori}, one has
 $ \big\| \ol{Y}^\e \big\|^p_{\hD^p} \+ \big\| \ol{Z}^\e \big\|^p_{\hZ^{2,p}}
 \+ \big\| \ol{U}^\e \big\|^p_{\hU^p}
   \ls  \cC  E \Big[    \big( \int_0^T \n |\ff^\e_A(s,0,0,0)| ds \big)^p \, \Big] $.
 It then follows from \eqref{eq:f173c},  \eqref{eqn-d011} and H\"older's inequality   that
  \beas
 && \hspace{-1.2cm}  E \bigg[  \b1_A \bigg( \,  \underset{s \in [t,t+\e]}{\sup} \big|\cY^\e_s\big|^p
 \+   \Big( \int_t^{t+\e} \n |\cZ^\e_s|^2 ds \Big)^{\frac{p}{2}}
 \+   \int_t^{t+\e} \n \|\cU^\e_s\|^p_{L^p_\nu} ds \bigg) \bigg] \=
  E \bigg[ \,  \underset{s \in [t,t+\e]}{\sup} \big|\ol{Y}^\e_s\big|^p
 \+   \Big( \int_t^{t+\e} \n |\ol{Z}^\e_s|^2 ds \Big) ^{\frac{p}{2}}
 \+   \int_t^{t+\e} \n \|\ol{U}^\e_s\|^p_{L^p_\nu} ds \bigg] \nonumber  \\
 && \hspace{-0.6cm} \ls  \big\| \ol{Y}^\e \big\|^p_{\hD^p} \+ \big\| \ol{Z}^\e \big\|^p_{\hZ^{2,p}}
 \+ \big\| \ol{U}^\e \big\|^p_{\hU^p}
   \ls  \cC  E \bigg[    \Big( \int_0^T \n |\ff^\e_A(s,0,0,0)| ds \Big)^p \, \bigg]
   \= \cC  E \bigg[    \Big( \int_t^{t+\e} \n \b1_A  |f(s,y\+   V (t,s,z,u)  ,  z,   u)| ds \Big)^p \, \bigg] \nonumber \\
 &&  \hspace{-0.6cm}   \ls \cC  E \bigg[ \, \b1_A   \Big( \int_t^{t+\e} \n \big[ |f(s,y,0,0)|
     \+  \wt{\beta}_s    |V (t,s,z,u)| \+ \fc |z|
  \+ \fc \|u\|_{L^p_\nu} \big]  ds \Big)^p \, \bigg] \nonumber \\
 &&  \hspace{-0.6cm}    \ls \cC  E \bigg[  \e^p \b1_A  \underset{s \in [t,t+\e]}{\sup}   |f(s,y,0,0)|^p
    \+   \e   \b1_A    \underset{s \in  (t,t+\e]}{\sup} \big| V (t,s,z,u) \big|^p
    \Big( \int_t^{t+\e} \n \wt{\beta}^q_s ds   \Big)^{p-1}
 \+ \e^p \fc^p \, \b1_A  \big( |z|^p \+ \|u\|^p_{L^p_\nu} \big)   \bigg]    .
  \eeas
 As $A $ varies over $\cF_t$, \eqref{eq:x201} shows that
\bea
&& \hspace{-1.5 cm} E \bigg[    \,  \underset{s \in [t,t+\e]}{\sup} \big|\cY^\e_s\big|^p
 \+   \Big( \int_t^{t+\e} \n |\cZ^\e_s|^2 ds \Big)^{\frac{p}{2}}
 \+   \int_t^{t+\e} \n \|\cU^\e_s\|^p_{L^p_\nu} ds \Big| \cF_t \bigg] \nonumber \\
 &&  \ls \cC   E \bigg[  \e^p \underset{s \in [t,t+\e]}{\sup}   |f(s,y,0,0)|^p
     \+  \e C_{\wt{\beta}}^{p-1}  \underset{s \in  (t,t+\e]}{\sup} \big| V (t,s,z,u) \big|^p    \Big| \cF_t \bigg]
     \+ \cC \e^p    \big( |z|^p \+ \|u\|^p_{L^p_\nu} \big)    \nonumber \\
     &&  \ls \cC  \e^p    E \bigg[  \underset{s \in [t,t+\e]}{\sup}   |f(s,y,0,0)|^p
        \Big| \cF_t \bigg]
        \+   \cC  \e^p    \big(|z|^p \+ \|u\|^p_{L^p_\nu} \big) .
     \label{eq:x211}
\eea

\no {\bf 3)}
Since $ \cY^\e_{t+\e} \= Y^\e_{t+\e} \- y \-  V \big(t, t \+ \e  ,z,u \big)
\= Y^\e_T \- y \-  V \big(t, t \+ \e  ,z,u \big) \= 0 $, \pas, we see from the BSDEJ\,$(0,\ff^\e)$ that
\beas
    \cY^\e_t   \=    \int_t^{t+\e}
 \n \ff^\e (s, \cY^\e_s,\cZ^\e_s , \cU^\e_s )   ds
\- \int_t^{t+\e} \n \cZ^\e_s dB_s \- \int_{(t,t+\e]} \n \int_\cX \n \cU^\e_s (x) \tnp (ds,dx)    , ~ \pas
\eeas
   Taking the conditional expectation $E[~|\cF_t]$ yields that
 \beas
  \frac{1}{\e} \cY^\e_t \-f(t,y,z,u) \=     \frac{1}{\e}   E\bigg[  \int_t^{t+\e} \n
   \big(f(s,\cY^\e_s \+ y + V(t,s,z,u) ,\cZ^\e_s\+z,\cU^\e_s \+ u) \-f(t,y  ,z,u) \big)
   ds\Big|\cF_t \bigg] , \q \pas
 \eeas
   H\"older's inequality and \eqref{eq:f173c}   imply that \pas
 \beas
  && \hspace{-1.3cm} \Big| \frac{1}{\e} \cY^\e_t \-f(t,y,z,u) \Big|
    \ls   \frac{1}{\e} E\bigg[  \int_t^{t+\e} \n
 \Big(  \wt{\beta}_s \big(|\cY^\e_s| \+   |V(t,s,z,u)| \big) \+  \fc |\cZ^\e_s| \+  \fc \|\cU^\e_s\|_{L^p_\nu}
   \+ \big|f(s,y  ,z,u) \-f(t,y  ,z,u)\big| \Big) ds\Big|\cF_t \bigg]   \\
 && \hspace{-0.8cm}  \ls \frac{1}{\e} E\bigg[ \Big( \n \int_t^{t+\e} \n \wt{\beta}^q_s ds \Big)^{\frac{1}{q}}
  \Big\{ \Big( \n \int_t^{t+\e} \n   |\cY^\e_s|^p ds \Big)^{\frac{1}{p}}
  \+ \Big( \int_t^{t+\e} \n  |V(t,s,z,u)|^p ds \Big)^{\frac{1}{p}} \Big\}
  \+ \fc    \e^{\frac12} \Big( \int_t^{t+\e} \n |\cZ^\e_s|^2 ds \Big) ^{\frac{1}{2}}
  \+    \fc    \e^{\frac{1}{q}}   \Big(\int_t^{t+\e} \n \|\cU^\e_s\|^p_{L^p_\nu} ds \Big)^{\frac{1}{p}}        \\
  && \hspace{-0.8cm} \q  +    \int_t^{t+\e} \n |f(s,y ,z,u  ) \-f(t,y,z,u)  |  ds \Big|\cF_t \bigg] .
 \eeas
 Then one can deduce from  Jensen's inequality, \eqref{eqn-d011}, \eqref{eq:x201}    and \eqref{eq:x211}   that
 \bea
 && \hspace{-1.2cm} 5^{1-p} \Big| \frac{1}{\e} \cY^\e_t \- \n f(t,y,z,u) \Big|^p
    \n \ls  \n       \frac{1}{\e^p}   E \bigg[ \e C_{\wt{\beta}}^{p-1}
    \underset{s \in [t,t+\e]}{\sup} \big( |\cY^\e_s|^p  \+   |V(t,s,z,u)|^p \big)
 \n \+ \fc^p  \e^{\frac{p}{2}} \Big( \int_t^{t+\e} \n |\cZ^\e_s|^2 ds \Big) ^{\frac{p}{2}}
 \n \+   \fc^p  \e^{p-1} \dn
  \int_t^{t+\e} \n \|\cU^\e_s\|^p_{L^p_\nu} ds
     \Big| \cF_t  \bigg] \n  \+ \vth_\e \nonumber \\
 &&  \hspace{-0.7cm}  \ls   \frac{C_{\wt{\beta}}^{p-1} \n  \+ \fc^p}{\e }    E \bigg[ \,  \underset{s \in [t,t+\e]}{\sup} |\cY^\e_s|^p
 \+   \Big( \int_t^{t+\e} \n |\cZ^\e_s|^2 ds \Big) ^{\frac{p}{2}}
 \+   \int_t^{t+\e} \n \|\cU^\e_s\|^p_{L^p_\nu} ds  \Big| \cF_t  \bigg]
 \+ \frac{C_{\wt{\beta}}^{p-1}}{\e^{p-1} }    E \bigg[ \, \underset{s \in [t,t+\e]}{\sup} |V(t,s,z,u)|^p \Big| \cF_t  \bigg]
  \+ \vth_\e    \nonumber \\
   &&  \hspace{-0.7cm}    \ls    \cC \e^{p-1}
     E \bigg[  \underset{s \in [t,t+\e]}{\sup}   |f(s,y,0,0)|^p
        \Big| \cF_t \bigg]   \+ \cC (\e^{p-1} \+ \e^{1- \frac{p}{2}} )   \big(|z|^p \+ \|u\|^p_{L^p_\nu} \big)  \+ \vth_\e ,
        \qq    \label{eq:x219}
 \eea
  where $\vth_\e \df  E\Big[  \big( \frac{1}{\e} \int_t^{t+\e} \n | f (s,y ,z,u  ) \- f (t,y,z,u)  |  ds
  \big)^p \Big|\cF_t \Big] $.

  The right continuity of process $ \{f(s,y,z,u)\}_{s \in [0,T]}$ at $t$ implies that
  $ \lmt{\e \to 0+} \frac{1}{\e} \int_t^{t+\e} \n |f(s,y ,z,u  ) \-f(t,y,z,u)  |  ds \= 0 $, \pas~
  Since \eqref{eq:f173c} shows that $ \frac{1}{\e} \int_t^{t+\e} \n |f(s,y ,z,u  ) \-f(t,y,z,u)  |  ds
  \ls 2 \underset{s \in [t,t+\d]}{\sup} |f(s,y,z,u)|
  \ls 2 \underset{s \in [t,t+\d]}{\sup} |f(s,y,0,0)|  \+ \n 2 \fc \big(|z|  \+ \|u\|_{L^p_\nu} \big) $, \pas ~
  for any $\e \ins \big( 0 , \d \ld 1 \big] $
  and since $E \Big[ \, \underset{s \in [t,t+\d]}{\sup} |f(s,y,0,0)|^p \Big] \< \infty$,
  a conditional-expectation version of  the  dominated convergence theorem shows that
  $ \lmt{\e \to 0+} \vth_\e \= 0$. Thus, letting $\e \to 0 $ in \eqref{eq:x219} yields that
 \beas
 \hspace{3cm}
 f(t,y,z,u) = \lmt{\e \to 0+} \frac{1}{\e } \cY^{\e }_t  \=
 \lmt{\e \to 0+} \frac{1}{\e } (Y^{\e }_t \- y )
 =  \lmt{\e \to 0+} \,  \frac{1}{\e } \big(Y^{t+\e,y,z,u}_t    \- y \big) , \q \pas
 \hspace{3cm} \hb{\qed}
  \eeas

\appendix
\renewcommand{\thesection}{A}
\refstepcounter{section}
\makeatletter
\renewcommand{\theequation}{\thesection.\@arabic\c@equation}
\makeatother

\section{Appendix}

\begin{lem} \label{lem_lp_esti}
 Let  $\{a_i\}_{i \in \hN} \subset [0, \infty)$. For any $p \in (0,\infty)$ and  $n \in \hN$ with $n \ge 2$, we have
   \bea   \label{eqn-d011}
  \big( 1 \land n^{p-1} \big) \sum_{i=1}^n a_i^p \le
\left(\sum_{i=1}^n a_i \right)^p \le \big( 1 \vee n^{p-1} \big) \sum_{i=1}^n a_i^p.
\eea
\end{lem}

 \no {\bf Proof:} Suppose that  $p \in (1, \infty)$ first. For any $0 \le b \le c <    \infty$, one can deduce that
\bea \label{eqn-d015}
   (b+c)^p-c^p 
   = p \int_c^{b+c}   t^{p-1} dt  \ge  p \int_c^{b+c}  b^{p-1} dt =p b^p
  \ge  b^p,  \;\hb{ or equivalent, }\;  (b+c)^p \ge  b^p+c^p   .
\eea
Thus, $(a_1+a_2)^p \ge  a^p_1+a^p_2 $.   When $n \ge 3$, applying \eqref{eqn-d015} consecutively, we obtain
 \bea \label{eqn-d019}
  \left(\sum_{i=1}^n a_i \right)^p \ge a^p_1+\left(\sum_{i=2}^n a_i \right)^p  \ge a^p_1+a^p_2+\left(\sum_{i=3}^n a_i \right)^p
\ge \cds \ge  \sum_{i=1}^{n-2} a_i^p  +\left(\sum_{i=n-1}^n a_i \right)^p \ge \sum_{i=1}^n a_i^p.
 \eea

  Now, let $\fm_n$ be the counting probability measure on $ S_n=\{1,\cds, n\}$ with $\fm_n(i)=\frac{1}{n}$ for each $i \in S_n $.  Jensen's inequality implies that
 \beas
   \left(\sum_{i=1}^n \frac{a_i}{n} \right)^p=\left(\int_{S_n} a_i \fm_n(di)\right)^p \le \int_{S_n} a^p_i \fm_n(di) =\sum_{i=1}^n \frac{a^p_i}{n} .
 \eeas
  Multiplying $n^p $ to both sides, we see from  \eqref{eqn-d019}  that
 \bea  \label{eqn-d025}
 \sum_{i=1}^n  a^p_i \le \left(\sum_{i=1}^n a_i \right)^p \le n^{p-1} \sum_{i=1}^n  a^p_i.
 \eea

Clearly, the case\;``$p=1$" is trivial. So it remains to show \eqref{eqn-d011} for $p \in (0,1)$:  Applying \eqref{eqn-d025} with $\wt{p}=\frac{1}{p}$ and $\wt{a}_i=a^p_i$, $i \in S_n $ yields that
 \beas
   \sum_{i=1}^n  a_i  =\sum_{i=1}^n  \wt{a}^{\wt{p}}_i  \le   \left(\sum_{i=1}^n \wt{a}_i \right)^{\wt{p}}=\left(\sum_{i=1}^n a^p_i \right)^{\frac{1}{p}}
    \le n^{\wt{p}-1} \sum_{i=1}^n  \wt{a}^{\wt{p}}_i  = n^{\frac{1}{p}-1} \sum_{i=1}^n  a_i.
 \eeas
Taking $p-$th power on both inequalities above, we obtain
 \beas
\hspace{6cm}  n^{ p-1}  \sum_{i=1}^n a^p_i  \le \left(\sum_{i=1}^n  a_i  \right)^p \le  \sum_{i=1}^n a^p_i.   \hspace{6cm} \hb{\qed}
 \eeas

\begin{lem} \label{lem_a02}
     For any $b, c \in [0, \infty)$, we have
  \bea   \label{eqn-d031}
 \big|b^p-c^p\big| \le \left\{ \ba{ll}
 |b-c|^p,   & \hb{if } p \in (0,1], \vspace{1mm}\\
   p(b \vee c)^{p-1}|b-c| , \qq \q & \hb{if } p \in (1, \infty) .
\ea
  \right.
   \eea
\end{lem}

 \no {\bf Proof:} It is trivial when $b=c$. Since $b$ and $c$ take the symmetric roles in \eqref{eqn-d031},  we only need to assume $b < c$ without loss of generality.

 \no $\bullet$ When $p \in (0,1]$, applying Lemma \ref{lem_lp_esti} with $a_1=b$ and $a_2 = c-b$ yields that
 $c^p = (a_1+a_2)^p \le  a^p_1+a^p_2 =  b^p + (c-b)^p   $,
  which implies that $  \big|b^p-c^p\big|= c^p-b^p \le (c-b)^p = |b-c|^p$;

 \no $\bullet$ When $p \in (1, \infty)$, one can deduce that $  c^p-b^p
   = p \int_b^c   t^{p-1} dt \le p \int_b^c   c^{p-1} dt = p c^{p-1} (c-b)$, which leads to that
 $  \big|b^p-c^p\big|= c^p-b^p \le p \, c^{p-1} (c-b) = p(b \vee c)^{p-1}|b-c|$.      \qed

\begin{lem}  \label{lem:theta-fcn} \(Bihari's inequality\)
Let $\th  : [0, \infty) \to [0, \infty)$ and
$  \z, \chi: [0, T] \to [0, \infty)$ be  three   functions    such that

 \no (\,\,i) either $\th \equiv 0$ or \,$\th(x)>0$ for any $x>0$;

 \no (\,ii) $\th$ is increasing and satisfies $\int_{0+}^1 \frac{1}{ \th(x)} dx= \infty$;

 \no (iii) $\z$ is integrable and $\chi$ is bounded.

  \no  If   $ \chi(t) \le \int_t^T   \th\big(\chi(s)\big) \z(s)  ds$ for any  $ t \in [0,T]$, then $\chi \equiv 0$.
 \end{lem}

  \no {\bf Proof:} The case\;``$\th \equiv 0$"  is trivial. So we only  assume
 that $\th(x)>0$ for any $x>0$ by (i).  It follows from (iii) that
 $ \int_0^T \th\big(\chi(s)\big) \z(s) ds \le \th \Big( \underset{s \in [0,T]}{\sup} \chi(s)  \Big)
 \int_0^T   \z(s) ds < \infty$. Thus
 \beas
h(t) := \int_t^T \th\big(\chi(s)\big) \z(s) ds \in [0,\infty),  \q  \fa t  \in [0,T]   .
\eeas
defines a continuous and decreasing function. Since $\th$ is increasing,  differentiating function $\f$ yields that
  \bea   \label{eqn-d041}
  h'(t)=-\th\big(\chi(t)\big) \z(t) \ge -\th\big(h(t)\big) \z(t), \q \fa t \in [0,T] .
 \eea

  Assume $h(0)>0$.
      Then $\wh{T} := \inf\{t\in (0,T]: h(t)=0\} \in (0,T]$ and it is clear that  $\lmtd{t \to \wh{T}-} h(t) = 0$.
  As the continuous and decreasing function $\f$ images $[0,\wh{T})$ onto $(0 ,h(0)]$,
Changing of variable  gives that
  \bea  \label{eqn-d045}
   \int_{0+}^{h(0)} \frac{1}{\th(x)} dx = -\int_0^{\wh{T}-} \frac{1}{\th\big(h(t)\big)} d h(t)
= -\int_0^{\wh{T}-} \frac{h'(t)}{\th\big(h(t)\big)} dt \le  \int_0^{\wh{T}-}  \z(t) dt   \le \int_0^T \z(t) dt < \infty ,
 \eea
where we used  \eqref{eqn-d041} and (iii).  For any $0<a<b< \infty$, one can deduce from the monotonicity of function $\th$  that
 $   \int_a^b \frac{1}{\th(x)} dx \le \int_a^b \frac{1}{\th(a)} dx = \frac{b-a}{\th(a)} < \infty$, which together with \eqref{eqn-d045} implies that
  $\int_{0+}^1 \frac{1}{\th(x)} dx  < \infty$. This results in a contradiction to assumption (ii). Therefore,
  $h(0)=0$, which forces    $h (\cdot) \equiv 0$. As a consequence,    $\chi(\cdot)    \equiv    0$.    \qed

      For the next three lemmas, we consider   a generic vector space $\hE$  with norm $\|\cd\|$.

\begin{lem} \label{lem_pi}
 Let $\hE$ be  a   vector space   with   inner product $\lan \cd, \cd \ran$ and norm $\|\cd\|$.
 For any $x,y \in \hE$, we have
 \bea \label{eqn-d051}
\big\|\pi_r(x)-\pi_r(y)\big\| \le \|x-y\|, \q \fa  r \in (0, \infty).
 \eea
  Consequently,
     \bea   \label{eqn-d055}
 \|x-y\|  \ge  \big(\|x\|\land \|y\|\big) \big\|  \sD(x) - \sD(y)  \big\|.
   \eea
\end{lem}

 \no {\bf Proof:} Without loss of generality, we assume that  $\|x\| \ls \|y\|$ in the whole proof.

 \ss  To see \eqref{eqn-d051}, let us  discuss by three cases:

  \no (1) When $r \>  \|y\|$: Since  $\pi_r(x) \= x $ and $\pi_r(y) \= y $,
 one simply has  $ \big\|\pi_r(x) \- \pi_r(y)\big\|  \=  \|x \- y\|$;

  \no (2) When $ \|x\|  \ls   r  \ls  \|y\|$:
 Let us set $c  \df  \big\lan x, \sD(y) \big\ran$ and $\wh{y}  \df  c  \sD(y)$.
 Since $\big\lan x \- \wh{y}, \sD(y) \big\ran  \=  0 $,
it holds for any $\a  \ins  \hR$ that
 \beas
\big\| x  -\a \sD(y) \big\|^2= \big\| x-\wh{y} -(\a \-  c) \sD(y) \big\|^2= \big\|  x-\wh{y}\big\|^2 + \big\|  (\a \-  c) \sD(y) \big\|^2 = \big\|  x-\wh{y}\big\|^2 +  (\a \-  c)^2   .
\eeas
Hence, it follows that
 \beas
  \big\|\pi_r(x)-\pi_r(y)\big\|^2 = \big\|x-r \sD(y)  \big\|^2 = \big\|  x-\wh{y}\big\|^2 +  (r \-  c)^2  \le  \big\|  x-\wh{y}\big\|^2 +  \big( \|y\| \-  c\big)^2 =  \big\|x-y \big\|^2 ,
\eeas
   where we used the fact that $c  \ls  \big| \big\lan x , \sD(y) \big\ran \big|  \ls  \|x\|
    \ls   r  \ls  \|y\|$    by the Schwarz inequality.

  \no (3) When $r < \|x\|  $:   We know from (2) that
 \beas
   \q      \|x \- y\| \ge \big\|\pi_{\|x\|}(x) \- \pi_{\|x\|}(y)\big\| =   \big\|x \-  \|x\| \sD(y)\big\|
   =\|x\|   \big\|\sD(x) \-  \sD(y)\big\|       \ge  r \big\|\sD(x)  \-  \sD(y)\big\| = \big\|\pi_r(x) \- \pi_r(y)\big\|.
 \eeas

If $x \= {\bf 0}$, \eqref{eqn-d055} holds trivially. Otherwise,
since $\|x\| \ls \|y\|$, applying \eqref{eqn-d051} with $r \= \|x\|$ gives rise to \eqref{eqn-d055}. \qed

\begin{lem} \label{lem_pi2}
 Let $\hE$ be  a   vector space   with   norm $\|\cd\|$ only.
 For any $x,y \in \hE$, we have
 \beas
\big\|\pi_r(x)-\pi_r(y)\big\| \le 2 \|x-y\|, \q \fa  r \in (0, \infty).
 \eeas

\end{lem}

 \no {\bf Proof:} Let $x,y \in \hE$. Since $ |a \ve b \- a \ve c | \ls |b\-c|$ holds for any $a,b,c \ins \hR $,
   the triangle inequality implies that
 \beas
 \hspace{1.5cm} \big\|\pi_r(x) \- \pi_r(y)\big\|
 & \tn \= & \tn  \Big\|\frac{r}{r \ve \|x\|}x \- \frac{r}{r \ve \|y\|}y \Big\|
 \ls \frac{r}{r \ve \|x\|} \|x \- y\| \+ \Big| \frac{r}{r \ve \|x\|} \- \frac{r}{r \ve \|y\|} \Big| \|y\| \\
  & \tn \ls & \tn  \|x \- y\| \+ \frac{r\|y\|}{(r \ve \|x\|)(r \ve \|y\|)} \big| r \ve \|x\| \- r \ve \|y\| \big|
 \ls \|x \- y\| \+ \big| \|x\| \-   \|y\| \big|
 \ls 2 \|x \- y\| .  \hspace{1.5cm} \hb{\qed}
 \eeas

\begin{lem} \label{lem_a04}
  Let $\hE$ be  a   vector space   with   inner product $\lan \cd, \cd \ran$ and norm $\|\cd\|$.
  For any $p \ins (0,1]$ and $x,y \ins \hE$, we have
  $\big\| \|x\|^p\sD(x) \- \|y\|^p\sD(y) \big\| \ls  (1 \+ 2^p)\|x \- y\|^p$.
\end{lem}

  \no {\bf Proof:} The case\;``$p \= 1$" is trivial since
  $\big\| \|x\| \sD(x) \- \|y\| \sD(y) \big\| 
 \=  \|x \- y\|   $. For  $p  \ins  (0,1)$, we assume without loss of generality
 that $\|x\|  \ls  \|y\|$ and  discuss by three cases:

 \no 1) When $x \= 0$: $\big\| \|y\|^p\sD(y) \big\|  \= \|y\|^p $;

 \no 2) When $0  \<  \|x\|  \ls  \|x-y\|$: $\big\| \|x\|^p\sD(x) \- \|y\|^p\sD(y) \big\|
   \ls  \big\| \|x\|^p\sD(x) \big\| \+  \big\| \|y\|^p\sD(y) \big\|  \=  \|x\|^p \+  \|y\|^p
   \ls  \|x\|^p \+  \big( \|x\| \+ \|x \- y\| \big)^p  \ls   (1 \+ 2^p)\|x \- y\|^p$;

 \no 3) When  $ \|x\|  \>  \|x-y\|$: As $\|x\|  \ls  \|y\|$, \eqref{eqn-d055} and Lemma \ref{lem_a02} show that
$\big\| \|x\|^p\sD(x) \- \|y\|^p\sD(y) \big\|
   \ls  \|x\|^p \big\| \sD(x) \-  \sD(y) \big\|  \+  \big| \|x\|^p  \-  \|y\|^p \big|
    \ls   \|x\|^{p-1} \| x \-  y \|
     \+  \big| \|x\|  \-  \|y\|  \big|^p    \<   2\|x \- y\|^p $.    \qed

\bibliographystyle{siam}
\bibliography{lp_sol}

\end{document}